\documentclass{article}

\usepackage{natbib}
 \bibpunct[, ]{(}{)}{,}{a}{}{,}%
\usepackage[text={6.5in,9in},centering]{geometry}

\usepackage{setspace}
\usepackage{amssymb}
\usepackage{amsmath}
\usepackage{amsthm}
\usepackage{epsfig}
\usepackage{graphics}
\usepackage{graphicx}
\usepackage{verbatim}
\usepackage{latexsym}
\usepackage{multirow}
\usepackage{rotating}
\usepackage{hyperref}
\usepackage{url,colordvi}
\usepackage{color}
\usepackage{bbm}
\usepackage{dsfont}
\usepackage{tabularx, booktabs}
\usepackage{float}
\usepackage{subfigure}
\usepackage{mathrsfs}
\usepackage{enumitem,authblk}

\newcommand {\otoprule }{\midrule[\heavyrulewidth]}

\newcommand{\bfb}{\mathbf{b}}
\newcommand{\bfe}{\mathbf{e}}
\newcommand{\bfm}{\mathbf{m}}

\newcommand{\bfS}{\mathbf{S}}
\newcommand\1{\mathds{1}}
\newcommand{\indi}[1]{\1_{\{#1\}}}

\newcommand{\Jcal}{\mathcal{J}}
\newcommand{\Ical}{\mathcal{I}}
\newcommand{\Kcal}{\mathcal{K}}

\newcommand{\mbf}{\mathbf{m}}

\newcommand{\dbf}{\mathbf{d}}
\newcommand{\ybf}{\mathbf{y}}
\newcommand{\tildeV}{\tilde{V}}

\DeclareMathOperator*{\argmax}{arg\,max}

\newtheorem{theorem}{Theorem}
\newtheorem{remark}{Remark}
\newtheorem{lemma}{Lemma}
\newtheorem{corollary}{Corollary}
\newtheorem{proposition}{Proposition}

\newtheorem{definition}{Definition}

\newcommand{\bo}{\textcolor{black}}       
\newcommand{\nan}{\textcolor{black}}       
\newcommand{\nanr}{\textcolor{black}}  
\definecolor{olive}{RGB}{107,142,35}
\newcommand{\peterr}{\textcolor{black}}  

\newcommand\veclambda{\boldsymbol \lambda}

\title{Managing Appointment Booking under Customer Choices}

\author[1]{Nan Liu}
\author[2]{Peter M. van de Ven}
\author[3]{Bo Zhang}

\affil[1]{Operations Management Department, Carroll
School of Management, Boston College} 
\affil[2]{Centrum Wiskunde \& Informatica
(CWI), Amsterdam, The Netherlands}
\affil[3]{IBM Research AI}

\date{}

\begin{document}

\maketitle

\begin{abstract}
\nan{Motivated by the increasing use of online appointment booking
platforms, we study how to offer appointment slots to customers in order to maximize the total number of slots booked.  We develop two models, \textit{non-sequential} offering and \textit{sequential} offering, to capture different types of interactions between customers and the scheduling system. In these two models, the scheduler offers either a single set of appointment slots for the arriving customer to choose from, or multiple sets in sequence, respectively.} 
%
%
\nan{For the non-sequential model, we identify a static randomized policy which is asymptotically optimal when the system demand and capacity increase simultaneously, and we further \nanr{show that offering all available slots at all times has a constant} factor of 2 performance guarantee. For the sequential model, we derive a closed-form optimal policy for a large class of instances and develop a simple, effective heuristic for those instances without an explicit optimal policy.} 
\nanr{By comparing these two models, our study generates useful operational insights for improving the current appointment booking processes.}
\nan{In particular, our analysis reveals an interesting equivalence between the sequential offering model and the non-sequential offering model with perfect customer preference information. This equivalence allows us to apply sequential offering in a wide range of interactive scheduling contexts. Our extensive numerical study shows that sequential offering can significantly improve the slot fill rate (6-8\% on average and up to 18\% in our testing cases) compared to non-sequential offering.}

\end{abstract}

%
%


{\bf keywords} service operations management; customer choice;
appointment scheduling; Markov decision process; asymptotically optimal
policy



%



\section{Introduction}\label{sec.introduction}



\nan{Appointment scheduling is a common tool used by service firms (e.g., tech support, beauty services and healthcare providers) to \nanr{match their service capacity with uncertain customer demand.} 
With the widespread use of Internet and smartphones, customers often resort to online channels when searching for information and reserving services. To keep up with customers' preferences and needs, many service organizations have developed online appointment scheduling portals. For instance, TIAA allows its clients to book appointments with their financial consultants online. There are also a rising number of online service reservation companies that offer online appointment booking software or apps as a service for (small) businesses. Examples include zocdoc.com for medical appointments, opentable.com for dinner reservations, mindbodyonline.com for fitness classes, booker.com for spa services, and salonultimate.com for haircuts.}

\nan{The interfaces of these online appointment booking systems
vary.}  
%
%
%
\nanr{Some are more towards one-shot offering, i.e., a single list
of available appointments are shown on a single screen for customers
to choose from. Others offer a small number of options to start, and
customers must press ``more'' or ``next'' to view additional
appointments that are available. This way of scheduling resembles
the traditional telephone-based scheduling process, in which the
scheduling agent may reveal availability of appointment slots in a
sequential manner. Such a sequential way of displaying options is
often seen on mobile devices with a small screen as well. 
}

\nan{Our research is motivated by these various ways of appointment booking, and we seek to understand how a service provider can best use these (online) appointment booking systems. In scheduling practice, service providers first predetermine for each day an appointment template, which specifies the total number of slots, the length of each slot, and characteristics of customers (e.g., nature of the visit) to be scheduled for each slot. For instance, in a gym setting one has to determine the number of classes and their capacity, and in healthcare the service provider first determines the number of patients a clinician will see that day and at what times.
%
%
With an appointment template in place, service providers then decide how to assign incoming customer
requests to the available slots -- nowadays this process is often
done via online appointment scheduling as mentioned above.} 
The relevant performance metric for this process is the \emph{fill
rate}, i.e., the fraction of slots in a template booked before the
scheduling process closes. While the fill rate is not equivalent to
the eventual capacity utilization due to various post-scheduling
factors (e.g., cancellations, no-shows and walk-ins),
it is the first, and in many cases, the most important step to
achieving a high utilization (and thus a high revenue), and it is
the objective of the research presented in this paper.

Our focus is on modeling the scheduling process, and developing
stochastic dynamic optimization models to inform appointment
scheduling decisions in the presence of customer choice behavior.
Notwithstanding the surge of interest in 
service operations management in the past decade, basic single-day,
choice-based dynamic decision models are absent for a broad class of
real-world scheduling systems.
To our knowledge, the existing operations research and management
literature on this type of dynamic appointment scheduling is very
limited; \peterr{most, if not all, related research} 
%
%
%
assumes that customers reveal their preferences first and the
scheduler decides to accept or reject; see, e.g., \cite{GuptaWang08}
and \cite{wang2011adaptive}. 
%
%
%
%
\nan{However, as discussed above, in many real-world scheduling platforms the system (i.e. the scheduler) offers its
availability to customers to choose from either in a one-shot
format or in a sequential manner, with no explicit knowledge on customer preferences. Customers interact with the
scheduler in ways that have not been fully explored in the literature.}
%
%
This paper fills a gap in the literature by proposing the first
choice-based dynamic optimization models for making scheduling
decisions in systems where customers are allowed to choose among
offered appointment slots from an established appointment template.
%
%
We demonstrate how the current appointment booking processes can be improved by
developing optimality results, heuristics and managerial insights in
the context of the proposed models.

%

We propose and study two models for the interaction between customers and the service provider. The first one is referred to as
the \emph{non-sequential offering} model. In this model, the
scheduler offers a single set of appointment slots to each customer.
If some of the offered slots are acceptable to the customer, she
chooses one from them; otherwise, she does not book an appointment.
\nan{This simple, one-time interaction resembles the mechanism of many online appointment systems which provide one-shot
offerings, and our results on this model have direct
implications on how to manage these systems.} Our second model is a
\emph{sequential offering} model, in which the scheduler may offer
several sets of appointment choices in a sequential manner. This is
motivated by 1) web-based appointment applications designed to
reveal only a small number of appointment options, \nan{one web page
at a time (e.g., mobile-based appointment applications)};
%
%
and 2) the traditional telephone-based scheduling process, in which
the scheduler offers appointment slots sequentially. This second
model is stylized in the sense that it does not incorporate customer
recall behavior (i.e., a customer choosing a previously offered slot
after viewing more offers), which is allowed in both online and
phone-based scheduling. Our goal here is to glean insights on how the fill rate can be improved by ``smarter'' sequencing when
sequential offering is part of the scheduling process.

For both cases we are interested in which slots to offer in order to
improve and maximize the fill rate. We answer this question by
investigating the optimal offering policy using Markov decision
processes (MDPs), as well as by discussing heuristics. \nan{Intuitively,
sequential offering should lead to a higher fill rate than
non-sequential offering, because sequential offering gives the
scheduler more control over the service capacity. We are also
interested in how much improvement a service provider can get by
switching from non-sequential scheduling to sequential scheduling.
We answer this question by comparing the fill rates resulting from these two
models, and the gap in the fill rates represents the ``value'' of
sequential offering.}

\nan{We make the following main contributions to the literature.
\begin{itemize}[listparindent=2em]
\item To the best of our knowledge, our paper is the first to study
and compare two main scheduling paradigms, non-sequential (online)
and sequential (mobile- or telephone-based), used in the service
industries.
\item For the non-sequential offering model, we characterize the
optimal policy for a few special instances, and demonstrate that the
optimal policy can be highly complex in general. We then identify a
static randomized policy (arising from solving a single linear
program) which is \emph{asymptotically optimal} when the system
demand and capacity increase by the same factor. 
%
We further show that the
offering-all policy (i.e., offering all available capacity
throughout) has a \emph{constant factor of 2 performance guarantee}.
\item For the sequential offering model, we show that there exists an optimal policy that offers slot types one at a time based on their marginal
values. We are able to determine these values for a broad class of
model instances, which leads to a \emph{closed-form}
optimal policy in these cases. For model instances without an explicit optimal policy we develop a simple, effective heuristic.  
\item We show that a sequential offering model is
equivalent to a non-sequential offering model with perfect customer
preference information. This equivalence \nanr{ensures that sequential offering can be optimally applied in various interactive scheduling contexts, in particular when customer-scheduler interaction can (partially) reveal customer preference information during the appointment booking process.}
%
\item Via extensive numerical experiments, we demonstrate that the
offering-all policy and the heuristic developed for sequential
offering work remarkably well in their respective settings, and thus can
serve as effective approximate scheduling policies for practical
use. We also show that by switching from non-sequential to
sequential offering, the slot fill rate can be significantly improved
(6-8\% on average and up to 18\% in our testing cases).
\end{itemize}
}

The remainder of the paper is organized as follows.
Section~\ref{sec.lit.review} briefly reviews the relevant
literature.
Section~\ref{sec:demandmodel} introduces the common capacity and
demand model that will be used in both the non-sequential and
sequential settings. Sections~\ref{sec.non.seq.model}
and~\ref{sec.seq.model} discuss the non-sequential offering case and
the sequential offering case, respectively. Section
\ref{sec.num.results} presents an extensive numerical study that
complements our analytic work.
%
%
In Section~\ref{sec.conclusion}, we make concluding remarks. All proofs of our technical
results can be found in the Online Appendix.

\subsection{Literature Review} \label{sec.lit.review}

From an application perspective, our work is related and complementary to the literature on appointment template design, a
topic that has been studied extensively \citep{CayirliVeral03,GuptaDenton08}.  
%
\nan{Our
work departs from this literature in that we start from an
established template, and then study how to manage the interaction
between the customers and the scheduler in order to best direct
customers to various slots.}
%
%
Among the existing work on dynamic appointment scheduling,
\cite{feldman2014appointment} is the only study, other than the few
papers mentioned in the previous section, that explicitly models
customer choice behavior. However, \cite{feldman2014appointment}
focus on customer choices across different days and use a newsvendor
model to capture the use of daily capacity; this aggregate daily
capacity model does not allow them to
%
%
consider (allocating customers into) detailed appointment time slots
within a daily template.

From a modeling perspective,~\cite{zhang2005revenue} looks at a
similar choice model to ours, in the context of revenue management
for parallel flights. In contrast to the present paper, their
approach focuses on deriving bounds on the value function of the
underlying MDP, and using them to construct heuristics. Three recent
studies on assortment optimization are particularly relevant to our
paper: \cite{bkx2015}, \cite{golrezaei2014real} and
\cite{gallego2016online}. \cite{bkx2015} study a dynamic assortment
customization problem, mathematically similar to our non-sequential
appointment offering problem, assuming multiple types of customers,
each of which has a multinomial logit choice behavior over all
product types.
%
%
They assume that the customer type is observable to the seller
(corresponding to our scheduler), which differs from our setting.
\cite{golrezaei2014real} adopt a general choice model 
%
and also allow an arbitrary customer arrival process. 
\cite{gallego2016online} extend the work by \cite{golrezaei2014real}
to allow rewards that depend on both the customer type and product
type.
%
The last two studies assume that the customer type is known to the
seller, and their focus is on developing control policies
competitive with respect to an offline optimum, a different type of
research question from ours. The other distinguishing feature of our
research from all previous work is that we consider sequential
offering, an offering paradigm which has not been studied before.
%

%

Finally, our work is related to two other branches of literature.
The first on online bipartite matching \citep{mehtareview}, 
%
and the
second on general stochastic dynamic optimization, in particular
stochastic depletion problems (e.g., \citealt{ChanFarias}) and submodular optimization (e.g., \citealt{golovin}). These two lines of
research mainly aim to obtain performance guarantee results with
respect to offline optimums, which is not our research goal.

\section{Capacity and Demand Model} \label{sec:demandmodel}
We consider a single day in the future that has just opened for
appointment booking. The day has an established appointment
template, but none of the slots are filled yet. We divide the
appointment scheduling window, i.e., the time between when the day is first opened
for booking and the end time of this booking process, into $N$ small
periods. Specifically, we consider a discrete-time $N$-period
dynamic optimization model with $I$ customer types (that may come)
and $J$ appointment slot types (in the template), where customer types are
characterized by their set of \emph{acceptable} slot
types. Denote by $\Omega_{ij}$ the 0-1 indicator of whether slot
type $j$ is acceptable by customer type $i$, so the $I \times J$
\emph{choice matrix} $\Omega:=[\Omega_{ij}]$ consists of distinct
row vectors, each representing a unique customer type. \nanr{Such a customer type structure is similar to those in the literature that model customer segments characterized by different product preferences (e.g., \citealt{bkx2015}).}


%
%
%
%

We now present the details of our customer arrival and choice model.
%
%
In each period at most one customer arrives. The customer is type $i$
with probability $\lambda_i>0$, and with probability $\lambda_0:= 1-
\sum_{i=1}^I \lambda_i$ no customer arrives.
Upon a customer arrival, the scheduler offers her a set
$S\subseteq\{1,...,J\}$ of slot types, \emph{without} knowledge of
the customer type. When {\em offer set} $S$ contains one or more
acceptable slot types, the customer chooses one uniformly at random.
If no type in $S$ is acceptable to this customer, we distinguish two
possibilities. Either we use a {\em non-sequential} model where the
scheduler can only offer a single set, and the customer immediately
leaves \nanr{if none of the offered slots are acceptable} 
(Section
\ref{sec.non.seq.model}), or we use a {\em sequential} model where the
scheduler may offer any number of sets sequentially, until the
customer either encounters an acceptable slot, or the customer finds
no acceptable slots in any offer set and leaves without booking a
slot (Section \ref{sec.seq.model}). We start from an initial
capacity of $b_j$ slots of type~$j$ at the beginning of the
reservation process, and denote $\mathbf{b}:=(b_1,...,b_J)$. Every
time a customer selects a slot, the remaining slots of this
type are reduced by 1.
The scheduler aims to maximize the fill rate at the end of the
reservation process by deciding on the offer set(s) in each
period. This is also equivalent to maximizing the \emph{fill
	count}, i.e., the total number of slots reserved at the end of the booking process, because
the initial capacity $\mathbf{b}$ is fixed.

Our capacity and demand model generalizes that of
\cite{wang2011adaptive} in the following sense. Our notion of
`slot type' can be viewed as an abstraction of the service provider and time block combination in their model, and thus we allow a
generalization of using other attributes of a slot that may affect
its acceptability to customers, such as duration.
%
%
\cite{wang2011adaptive} consider distinct customer panels, each
characterized by a possibly different acceptance probability
distribution over all possible combinations of service providers and time
blocks and a set of revenue parameters. In contrast, we define the
notion of customer type and identify it with a unique set of
acceptable slot types. Their arrival rate (probability) parameters
are associated with each customer panel, while we directly have the
demand rate for each of the $I$  customer types as model
primitives.

\nan{Our choice model assumes that for a particular customer type,
	slot types are either ``acceptable'' or ``unacceptable''.} \bo{This dichotomized classification of slots closely mimics
	the decision process on whether a time slot works for one's daily schedule. For instance, such a slot-choosing process is seen at the popular polling website 
	www.doodle.com, where each participant responds to a poll by indicating whether a particular time works (i.e., is acceptable) by him or her.}
	\peterr{This relatively parsimonious choice model enables a tractable analysis of the interplay between appointment booking and customer choice.}
	\nanr{Its parameters may for instance be estimated by conducting a market survey on customers' acceptance on various slot types.}



As discussed earlier, the distinction
between the non-sequential and sequential customer-scheduler
interactions reflects the differences present in various real-life
appointment scheduling systems. The non-sequential model is best
suited for web-based appointment scheduling systems such as
\url{www.zocdoc.com}. In such systems the customer 
%
is presented with a list of time slots to choose from, which
corresponds to a single offer set. In contrast, sequential
scheduling reflects the iterative nature of, for instance,
telephone-based appointment scheduling. Here the scheduler 
may propose one or more slots initially, and
may present more if these are rejected by the customer. \bo{While allowing an unlimited number of offer sets in sequence does not conform with many real-world systems, the sequential model is a valuable object of study because the scheduler in this setting enjoys the greatest flexibility and hence the resulting optimal fill rate serves as an upper bound for that in both the non-sequential model and some intermediate paradigms such as those allowing a limited number of offer sets or with customer reneging}.

The assumption on the unobservability of the customer type is unique
in our work, and is present in all real-world systems that we
consider. \bo{Users of web and mobile-based appointment scheduling systems often prefer a simple interface soliciting no or minimal personal information before displaying availabilities}. Many telephone-based schedulers only know some
basic information of the customers. 
%
%
Even if these collected data are useful in predicting customer
preferences, \bo{many service firms 
may lack the necessary resources} (e.g., human, technology and software) to make such predictions and then
use them in scheduling decisions. This is another important
motivation why we choose to assume exact customer type is unknown to the scheduler in our
models.

Our objective is to maximize the fill rate (or equivalently,
fill count), thereby assuming that each customer
contributes to the objective equally. We choose this objective for
a few reasons. First, 
%
fill rate is a widely-used reporting metric by service firms for their operational and financial performance. The simplicity of
this metric also makes it more tractable for analysis. Second,
fairness may carry more weight than profitability in the
vision of a service firm, e.g., a healthcare delivery organization.
Third, while different customers may bring different rewards (e.g.,
revenues) to the service firm, how to associate such rewards with
customer (preference) types is not well understood in
the literature. \nanr{In the present study}, we choose a straightforward objective
instead, without guessing a complicated reward
structure lacking empirical support.

Finally, our discrete-time customer arrival model with at most one
arrival per period is widely accepted and used by many operations
management studies, including those on healthcare scheduling (e.g.,
\citealt{green2006managing}) and on revenue management (e.g., \citealt{talluri2004revenue, bkx2015}).
%
%
\nanr{One could set $N$, the total number of time periods, sufficiently large so that the probability of multiple customers arriving during a single period is negligible (and thus as is the probability of more than $N$ customers arriving in total).} This demand model can be used to approximate an inhomogeneous Poisson arrival process
\citep{subramanian1999airline}.

\nan{In the following sections, we focus on analyzing the models described above. We acknowledge that our models
	do not explicitly capture the rolling-horizon feature of the
	appointment scheduling practice, in which customers may book
	appointments in future days and unused capacity in a day is wasted when the day is past. However, the rolling-horizon
	multi-day scheduling model is known for its intractability \citep{LiuZiyaKulkarni08, feldman2014appointment}.
	The single-day model is more tractable and often used in the
	literature to generate useful managerial insights (e.g.,
	\citealt{GuptaWang08,wang2011adaptive}). Indeed, in Section~\ref{sec.multi.day.num} we will
	numerically demonstrate how our single-day models can inform
	decision making in a rolling-horizon multi-day setting. }

\section{Non-sequential Offering}
\label{sec.non.seq.model} We first consider the non-sequential
offering model, in which only one offer set $S$ is presented to each
arriving customer. Denote by $\mathbf{m}\leq\mathbf{b}$ a
$J$-dimensional, non-negative integer vector that represents the
current number of remaining slots of each type, and by
$\mathbf{e}_j$ the $J$-dimensional unit vector with its $j$th entry
being 1 and all others zero. Define
$\bar{S}(\mathbf{m}):=\{j=1,...,J:m_j>0\}$, the set of slot types
with positive capacity, and $V_n(\mathbf{m})$ as the expected
maximum number of appointment slots that can be booked from period
$n$ to period 1 with $\mathbf{m}$ slots available at the beginning
of period $n$. Note that we count time backwards.

Further, denote by $q_{ij}(S)$ the probability that slot type $j$ is
chosen conditional on a type-$i$ customer arrival and an offer set
$S \in \bar{S}(\mathbf{m})$. We have, for any $j$,
\begin{equation} \label{eqn.qij}
q_{ij}(S) = \left\{
\begin{array}{cc}
\frac{\Omega_{ij}}{\sum_{k\in S}\Omega_{ik}},&  ~\textrm{if
	$\sum_{k\in S}\Omega_{ik}>0$}, \\
0, & ~\textrm{otherwise.}
\end{array}
\right.
\end{equation}
%
%
Then, the probability that slot type $j$ is chosen when offer
set $S$ is given is
\begin{equation}\label{eqn.dist.ln}
q_j(S)= \sum_{i=1}^I \lambda_i q_{ij}(S),
\end{equation}
and the no-booking probability is $q_0(S)=1-\sum_{j=1}^Jq_j(S)$.
The optimality equation is
\begin{equation}\label{eqn.bellman.obj}
V_{n}(\mathbf{m}) =   \max_{S\subseteq\bar{S}(\mathbf{m})}\left[\sum_{j\in S}q_j(S)(1-\Delta_{n-1}^j(\mathbf{m}))\right]+V_{n-1}(\mathbf{m}),\mbox{~~for~}n=N,N-1,...,1,
\end{equation}
where $V_0(\cdot)=0$ and
$\Delta_{n-1}^j(\mathbf{m}):=V_{n-1}(\mathbf{m})-V_{n-1}(\mathbf{m}-\mathbf{e}_j)$
denotes the marginal benefit due to the $m_j$th unit of slot type
$j$ at period $n-1$.

\nan{We first analyze the non-sequential offering model for a few
	specific instances, and demonstrate that in general the
	optimal non-sequential offering policy seems to have no appealing
	structural properties.
	Thus, characterizing the optimal policy for general, large-scale
	non-sequential offering models is very challenging, if not
	impossible. \nanr{We then focus our efforts on constructing simple scheduling policies that have performance guarantees and may perform well in practice. We first} 
	%
	%
	consider a limited class of policies (called static randomized offering policies), and identify one such
	policy which is \emph{asymptotically optimal} \nanr{when we increase the system demand and capacity simultaneously}. \bo{We further show} that 
	\nanr{a simple} policy that offers all available
	slots at all times has a \textit{constant ratio of 2 performance guarantee}, independent of all model parameters. In Section \ref{sec.num.results},
	we show via extensive numerical instances that this offering-all
	policy significantly outperforms its theoretical
	bound. 
	It may thus serve as a simple, effective heuristic
	offering rule for \bo{many} practitioners in the non-sequential offering
	context.}



\subsection{Results for Specific Model Instances}

\begin{figure}
	\centering
	\subfigure[N model]{\includegraphics[width=0.14\textwidth]{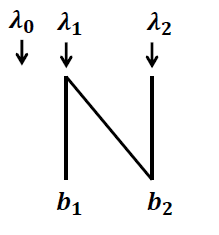}}\hspace{0.5in}
	\subfigure[W model]{\includegraphics[width=0.14\textwidth]{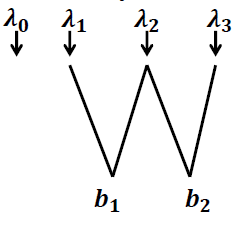}}\hspace{0.5in}
	\subfigure[M model]{\includegraphics[width=0.14\textwidth]{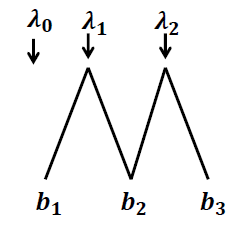}}\hspace{0.5in}
	\subfigure[M+1 model]{\includegraphics[width=0.16\textwidth]{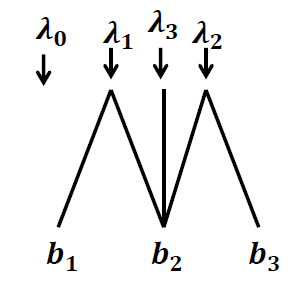}}
	\caption{The N, W, M, and M+1 model instances.}
	\label{fig:three.models}
\end{figure}

When there are $J=2$ slot types, the choice matrix $\Omega$ has two
possible non-trivial values:
$$
\Omega =
\left(
\begin{array}{ll}
1 & 1 \\
0 & 1 \\
\end{array}
\right),
\quad \textrm{and} \quad
\Omega =
\left(
\begin{array}{ll}
1 & 0 \\
1 & 1 \\
0 & 1 \\
\end{array}
\right).
$$
These we refer to as the N model instance (see Figure
\ref{fig:three.models}(a)) and the W model instance (see Figure
\ref{fig:three.models}(b)), respectively. These two model instances
are, for example, applicable to the popular Chinese scheduling system
\url{www.guahao.com.cn}, which allows customers to book either a
morning or an afternoon (medical) appointment for a certain day without
providing more granular time interval options. In both model instances, we show that it is optimal to offer all available slots at all times
\nan{(which we call the \emph{offering-all} policy in the rest of
	this article)}, as not doing so would unnecessarily risk sending
away certain customers. This is formalized in the following result.
\begin{proposition}
	\label{lemma.W.opt.policy} For the N and W model instances, the
	\nan{offering-all} policy is optimal.
\end{proposition}

When there are $J=3$ slot types, the simplest nontrivial choice
matrix is the M model instance in Figure \ref{fig:three.models}(c) with 
$$
\Omega =
\left(
\begin{array}{lll}
1 & 1 & 0\\
0 & 1 & 1 \\
\end{array}
\right).
$$
It turns out that in this case, the offering-all policy is not
always optimal;
%
%
rather, rationing of the versatile type-2 slot is needed. We define
policy $\pi_1$ according to its offer set:
\begin{equation}\label{eqn.opt.M}
S^{\pi_1}(\mathbf{m}):= \left\{
\begin{tabular}{ll}
$\{1,3\}$ & if $m_1>0$ and $m_3 > 0$,\\
$\bar{S}(\mathbf{m})$   &   otherwise.
\end{tabular}
\right.
\end{equation}
So policy $\pi_1$ proposes to hold back on offering type-2 slots
until either type-1 or type-3 slots are used up. We now formalize
that one cannot do better than this.
\begin{proposition}
	\label{lemma.M.opt.policy} For the M model instance, $\pi_1$ is optimal.
\end{proposition}
The intuition behind Proposition~\ref{lemma.M.opt.policy} is that
blocking slot type 2 does not lead to any immediate loss of customer
demand compared to offering it, while forcing early customers into
\emph{less popular} slot types (types 1 and 3). This 
preserves the \emph{popular (or, versatile)} slots (type 2) for
later arrivals, when slots run low. \nan{For convenience of
	discussion, we say a slot type is more popular (or, versatile) if
	this slot type is accepted by a superset of customer types compared to its counterpart.}

Following from Proposition \ref{lemma.M.opt.policy}, we know that a
versatile type 2 slot is at least as valuable as one of the other
two less popular slot types \nan{at all times, for otherwise it
	would be better to offer type 2 slots but not offering the more
	valuable, less popular slot type}. \nan{To be more specific, we have
	the following corollary.}
%
\begin{corollary}\label{mcorollary}
	In the M model instance, for either $j=1$ or 3 or both,
	\begin{equation}\label{m.morevaluable}
	V_n(\mathbf{m}-\bfe_2)\leq V_n(\mathbf{m}-\bfe_j), ~ \forall
	\mathbf{m}>0, ~ n \in \{1,\dots,N\}.
	\end{equation}
\end{corollary}


\nan{However, it is important to note that} one of the two less
popular slot types (1 and 3) may be \emph{strictly more valuable}
than the popular type 2. For example, for $\lambda_1=\lambda_2=0.5$,
it is easy to verify that $V_2(2,1,0)=1.625<1.75=V_2(2,0,1)$. The
reason here is the following. With $m_1=2$ and $n=2$, sufficient
capacity is available for potential type 1 customer demand (i.e., at
most 2 units). If $(m_2,m_3)=(1,0)$, the one unit of type 2 slot has
a positive probability of being taken by a type 1 customer (which
would be a waste); in contrast, if $(m_2,m_3)=(0,1)$, the one unit
of type 3 slot can only be exclusively offered to type 2 customers
(for whom no sufficient capacity is available), ant thus this is
more efficient. This simple example shows that because of customers'
ability to (randomly) choose from their offer set, \emph{less
	popular slots may be more valuable than versatile slots due to
	resource imbalance}. \nan{This observation implies that the (future)
	value of keeping a slot type cannot be viewed solely based on the
	number of accepting customer types, irrespective of the arrival
	probabilities or slot capacities. This complication renders the
	optimal policy for a general model instance quite complex, as we demonstrate now.}

The next model instance that we focus on is the M+1 model instance
shown in Figure \ref{fig:three.models}(d), with choice matrix
$$
\Omega = \left(
\begin{tabular}{lll}
1   &   1   &   0\\
0   &   1   &   1\\
0   &   1   &   0
\end{tabular}
\right).
$$
\nan{Note that the only difference between the M+1 and M model
	instances is the additional customer type 3 that only accepts type 2 slots.} It turns out that the
simple, elegant form of the optimal policies in the previous 
cases does not carry over to the M+1 model instance.
%
%
%

\begin{figure}
	\centering
	\subfigure[$\lambda_1=\lambda_2=0.1$, $\lambda_3=0.8$]{\includegraphics[width=0.4\textwidth]{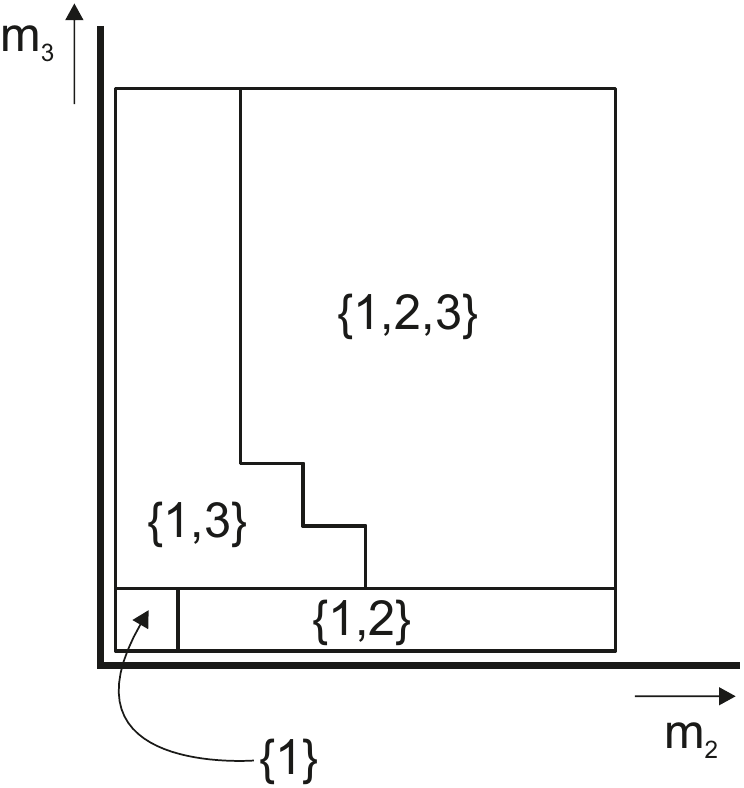}}\label{fig:M+1_1}\hspace{0.7in}
	\subfigure[$\lambda_1=\lambda_2=0.475$, $\lambda_3=0.05$]{\includegraphics[width=0.4\textwidth]{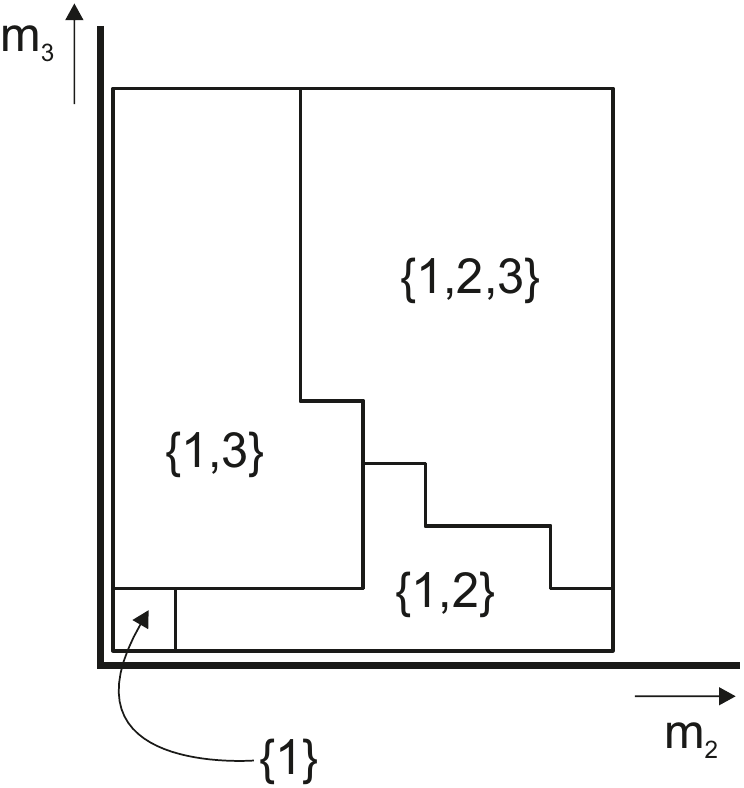}}\label{fig:M+1_2}
	\caption{The optimal policy for an M+1 model instance, with $m_1=4$, $n=5$.}
	\label{fig:M+1}
\end{figure}

To illustrate the complexity of the M+1 model instance, consider
the case with $m_1=4$ and $n=5$. Figure~\ref{fig:M+1} shows the
\emph{unique} optimal offer set, identified with $S \subset
\{1,2,3\}$, as a function of $m_2$ and $m_3$. (For instance, if
$S=\{1,3\}$, it means offering slot types 1 and 3 but not slot
type 2.) 
Consider $\lambda_1=\lambda_2=0.1$, $\lambda_3=0.8$
(Figure~\ref{fig:M+1}(a)) or $\lambda_1=\lambda_2=0.475$,
$\lambda_3=0.05$ (Figure~\ref{fig:M+1}(b)).
%
%
\nan{As discussed earlier, resource imbalance can make a less
	popular slot more valuable than a more popular one, which naturally
	would suggest an action of saving the less popular slot by only
	offering the versatile slot. Indeed, in Figure \ref{fig:M+1}(b), we
	see that action $\{1,2\}$ can be the \emph{unique} optimal action
	even when $\mathbf{m}>0$. This is true because $m_3$ is relatively
	small (equal to 1 or 2 in this case), while $m_1=4$ is ample given
	$n=5$ and the symmetric arrival rates of type 1 and type 2 customers.
	Blocking type 3 and offering the versatile type 2 earlier rather
	than later can help to resolve the resource imbalance by maximizing
	the total expected amount of type 2 slots taken by type
	2 customers \nanr{(and thus saving type 3 slots that can only serve type 2 customers for the future)}.} 
	

\nan{In addition, we see that the arrival rate now has a strong
	impact on the optimal policy, in contrast to the
	other cases we discussed so far: when $\lambda_3$ is large it is
	often optimal to include type 2 slots in the offer set, while for
	$\lambda_3$ small this is not the case. The reasoning here is that
	for $\lambda_3$ small the model is very close to the M model
	instance, for which we know it is optimal to save versatile
	type 2 slots for later in the booking process. However, not
	offering type 2 slots also implies turning away all type 3 customers,
	which explains why this slot type should be offered when $\lambda_3$
	is large.}

\nan{These observations we make on the M+1 model instance suggest
	that the optimal policy depends on the customer preference profiles,
	arrival rates and available slot capacity of the specific model
	instance under consideration. The optimal policy for a general model
	can be quite complex and have no straightforward structural
	properties. Thus, we shall focus our efforts on identifying
	\nanr{simple} 
	scheduling policies that have performance guarantees \nanr{and perform well in practice}.}

	\subsection{Asymptotically Optimal Policy}
	\label{sec.asymp} In this section, we \bo{construct} a
	\emph{static randomized policy} that is \emph{asymptotically
		optimal} when we increase the system demand and capacity simultaneously.
	We first introduce the class of static randomized policies. At any
	decision epoch, there are altogether $2^J$ possible actions in terms
	of which slot types to offer. Here we use a binary
	vector to denote the offer set, with  a 1 at position $j$ meaning that slot
	type $j$ is offered, and 0 otherwise. For example, we denote by the action
	of closing all slots as $\mathbf{w}^1:=(0,...,0)$, the
	$J$-dimensional zero vector, and the action of opening all slots as
	$\mathbf{w}^{2^J}:=(1,...,1)$, the $J$-dimensional one vector. We
	call the set of all $2^J$ $J$-dimensional 0-1 vectors as set
	$\mathcal{W}:=\{0,1\}^J$ and name the elements of the set as
	$\mathbf{w}^1,\mathbf{w}^2...,\mathbf{w}^{2^J}$. Define
	$\Kcal:=\{1,...,2^J\}$ as the action index set and so $\Kcal$ and
	$\mathcal{W}$ have the same cardinality.  
	
	A policy $\pi^{p}$ is a static randomized policy if
	$\pi^p$ offers $\textbf{w}^k$ with some fixed probability $p_k$,
	independent of the system state and the time period.\footnote{\nan{Even
			if some of the slot types are unavailable, $\pi^p$ would still offer
			these slot types according to the probability vector $\mathbf{p}$.
			However, customers only consider those slot types that are available in the booking process.}} The class of static randomized policies contains all
	$\pi^p$'s such that the vector $\textbf{p}=\{p_k\}_{k=1}^{2^J}$ is a
	probability vector. For instance, the offering-all policy is a
	special case in this class with $p_{2^J}=1$ and $p_k=0$ for all
	$k\neq 2^J$.
	
	We show that there exists a vector $\textbf{p}^*$ such that
	$\pi^{p^*}$ is asymptotically optimal when the demand and capacity
	are scaled up simultaneously. The choice of $p^*_k$ relies on the
	fluid model corresponding to the stochastic model
	\eqref{eqn.bellman.obj} considered above, in which we can readily determine the optimal offering policy. We choose $p^*_k$ such that it represents the fraction of the time in which the action
	$\textbf{w}^k$ is used in \bo{a} fluid model under optimal control. Below we \bo{construct} this asymptotically optimal policy
	$\pi^{p^*}$ and defer \bo{more} technical details to the Online Appendix.
	
	\subsubsection{Fluid Model}
	We first introduce our fluid model. To differentiate from the
	notation in the stochastic model formulation
	\eqref{eqn.bellman.obj}, we shall put the time index $n$ in
	parentheses, instead of as a subscript. Instead of discrete customers arriving in each slot, we represent a customer by a unit of fluid. In total one unit of demand arrives in each time period, a fraction $\lambda_i$ of which corresponds to customer type $i$. \nanr{This fluid is distributed evenly among all available slots that are offered and accepted by the corresponding customer type.}
	
	
	For each $n=1,...,N$, the
	decision vector in the fluid model is
	$\mathbf{z}(n)=(z_1(n),...,z_{2^J}(n))$, which is a $2^J$
	dimensional vector, each component $z_k(n) \in [0,1]$ representing
	the time during which action $k \in \Kcal$ is being used in period
	$n$. Note that each action can be used for any fractional unit of
	time. Thus we require that
	\begin{equation}
	0 \leq z_k(n) \leq 1, \forall k\in\Kcal, n=1,...N; \label{zzeroone}
	\end{equation}
	\begin{equation}
	\sum_{k\in\Kcal}z_k(n)=1, \forall n=1,...N.\label{zsum}
	\end{equation}
	Constraint \eqref{zsum} ensures that the total time spent on all
	possible actions (including the one that closes all slot types) in
	one period adds up to one.

	Let $\boldsymbol\tau(n)=[\tau_{k,j}(n)]$ be a $2^J\times J$ matrix,
	each row of which corresponds to one of the $2^J$ possible actions.
	We use $\tau_{k,j}(n)$ to indicate the amount of time for which type
	$j$ slots are offered during the time when the $k$th action is taken in
	period $n$. We have that
	\begin{equation}\label{opentime}
	\tau_{k,j}(n)= z_k(n) \mathbf{w}^k_j,~~ \forall k\in\Kcal, ~
	j=1,2,\dots,J, ~n=1,...N,
	\end{equation}
	where $\mathbf{w}^k_l$ denotes the $l$th coordinate of vector
	$\mathbf{w}^k$. Constraint \eqref{opentime} is presented mainly to
	make the formulation clearer and easier to understand. It ensures
	that slot type $j$ can be open when action $k$ is chosen only if
	action $k$ offers slot type $j$. If action $k$ does not offer slot
	type $j$, then $\mathbf{w}^k_j=0$ and $\tau_{k,j}(n)$ is zero by \eqref{opentime}. Let
	$\Jcal_k = \{j: \mathbf{w}^k_j=1,~ j=1,2,\dots,J.\}$ be the full set
	of slot types offered by action $k$. Note that \eqref{opentime}
	implies that
	\begin{equation*} 
	\tau_{k,j_1}(n) = \tau_{k,j_2}(n) , ~~\forall k\in\Kcal, ~n=1,...N,
	~j_1,j_2 \in \Jcal_k.
	\end{equation*}
	%
	That is, if an action $k$ offers multiple slot types, the offering
	durations of these slot types are the same.

	Let $y_{i,j}(n)$ denote the amount of type $j$ slot's capacity
	filled by type $i$ customers during period $n$ and $\Kcal_j = \{s:
	\mathbf{w}^s_j=1, s \in \Kcal \}$ be the index set of actions that
	offer type $j$ slots. If $\Omega_{i,j}=1$,
	\begin{equation}\label{yexp1}
	y_{i,j}(n)=\sum_{k\in \Kcal_j
	}\tau_{k,j}(n)\cdot\frac{\lambda_i}{\sum_{l=1}^J \min\{\Omega_{i,l},\mathbf{w}^k_l\}},~
	i=1,2,\dots,I,~j=1,2,\dots,J,~n=1,...,N;
	\end{equation}
	and otherwise if $\Omega_{i,j}=0$, then
	\begin{equation}\label{yexp2}
	y_{i,j}(n)=0,~ i=1,2,\dots,I,~j=1,2,\dots,J,~
	n=1,...,N.
	\end{equation}
	Note that all terms in \eqref{yexp1} except $\tau_{k,j}(n)$'s are
	constants and therefore \eqref{yexp1} as a set of constraints for
	the optimization problem is linear in the decision variables
	$\tau_{k,j}(n)$.
	
	Let $M_j(t)$ be the amount of type $j$ slots left with $t$ time
	periods to go and let $Z_N(\mathbf{m})$ be the optimal amount of
	(fluid) customers served with initial capacity vector $\mathbf{m}$
	and $N$ periods to go. The goal is to choose $z_k(n)$ (and
	$\tau_{k,j}(n)$) in order to solve for
	\begin{align}
	\tag{P1}
	& Z_N(\mathbf{m}) = \max  \sum_{n=1}^N  \sum_{j=1}^J \sum_{i=1}^I y_{i,j}(n), & \label{lp.obj} \\
	\mbox{s.t.~~}&\mbox{\eqref{zzeroone}, \eqref{zsum}, \eqref{opentime}, 
		\eqref{yexp1}, \eqref{yexp2}, and~}\notag\\
	& M_j(N)=m_j,  & j=1,2,\dots,J, \label{lp.initial}\\
	& M_j(n-1)=M_j(n)-\sum_{i=1}^I y_{i,j}(n), & j=1,2,\dots,J, ~n=1,...,N, \label{lp.flow}\\
	& M_j(n)\geq 0,  & j=1,2,\dots,J,~ n=0,1,...,N-1. \label{lp.nonneg}
	\end{align}
	In (P1), constraint \eqref{lp.initial} specifies the initial
	capacity vector, \eqref{lp.flow} updates the capacity vector for each period, and
	\eqref{lp.nonneg} ensures that all slot types have nonnegative
	capacity throughout. We remark that in our formulation, control can
	be exerted anytime continuously throughout the horizon but the
	system is observed only at discrete time epochs $0,1,2,\dots,N$ to
	match the stochastic model formulation \eqref{eqn.bellman.obj}.

	\subsubsection{Choice of $\mathbf{p}^*$}
	Let $p^*_k$ be the fraction of the time in which the optimal policy
	chooses action $k$ in the fluid model (P1). That is,
	\begin{equation}
	\label{eq:p*def} p^*_k = \frac{\sum_{n=1}^N z^*_k(n)}{N},
	\end{equation}
	where $z^*_k(n)$ is the optimal solution to (P1). We now translate this optimal policy for the fluid model to our original discrete and stochastic setting by defining a
	policy $\pi^{p^*}$ such that in each period $n$, this policy offers
	$\textbf{w}^k$ with probability $p^*_k$, independent of everything else.
	
	The intuition behind choosing $\mathbf{p}^*$ as the offering
	probability vector is that if we scale up the system demand (i.e.,
	$N$) and capacity (i.e, $\bfm$) in the stochastic model, using
	$\pi^{p^*}$ makes the proportion of total customer demand going to
	each slot type in the stochastic model approximately matches that in
	the fluid model. Thus, the total fill counts in the stochastic model
	is similar to that of the fluid model. Because the fluid model is a
	deterministic model which provides an upper bound on the objective
	value of the stochastic model (more on this below), we know that
	$\mathbf{p}^*$ is (close to) optimum in the stochastic model as the
	system becomes large. We formalize this intuition in the next
	section.

	\subsubsection{Main Result} Consider a sequence of problems indexed by
	$K=1,2,3\dots$. The problems in this sequence are identical except
	that for the $K$th problem, the number of total periods is $NK$ and
	the capacity vector is $\textbf{m} K$. We call the problem instance
	with $K=1$ as the base problem instance. Let $V_n^{\pi^p}(\cdot)$ be
	the total expected number of slots filled under a policy $\pi^p$ with the offering probabiity vector $\mathbf{p}$ in
	the stochastic model. The main result is shown in the following
	theorem.
	\begin{theorem} \label{prop.asym}
		\mbox{}\\
		(i) $K^{-1} V_{NK}(\textbf{m} K) \leq K^{-1} Z_{NK}(\bfm K)
		= Z_N(\bfm),~\forall \bfm
		\geq 0,~ K=1,2,3,\dots$; \\
		(ii) $\lim_{K \rightarrow \infty} K^{-1}
		V^{\pi^{p^*}}_{NK}(\textbf{m} K) = Z_N(\bfm)$.
	\end{theorem}
	Recall that $V_{n}(\cdot)$ is the optimal value of the stochastic
	model defined in \eqref{eqn.bellman.obj}. Thus Theorem 
	\ref{prop.asym}(i) says that the ``normalized'' optimal value of the
	non-sequential offering stochastic model (i.e., the original value
	divided by $K$) is \bo{bounded from above} by that of the corresponding fluid
	model, and that the normalized objective value of the fluid model is
	the same as the objective value of the base fluid model with $K=1$.
	Theorem \ref{prop.asym}(ii) states that as the system grows
	large, the normalized objective value in the stochastic system under
	policy $\pi^{p^*}$ converges to this constant upper bound, and thus
	$\pi^{p^*}$ is asymptotically optimal.
	
	The proof of \bo{Theorem} \ref{prop.asym} entails a few key steps
	which are outlined below (full details can be found in the Online
	Appendix). We first show that the optimal objective value of the
	fluid model is an upper bound of the optimal value of the stochastic
	model, i.e., $Z_N(\bfm) \geq V_N(\bfm)$, for any given set of model
	parameters. Then, based on any static randomized policy $\pi^p$, we
	construct a lower bound for $V_n^{\pi^p}(\cdot)$, and this lower
	bound is naturally a lower bound for the optimal value of the
	stochastic model $V_N(\bfm)$ (because $\pi^p$ is not necessarily
	optimal). Finally, we show that when $p$ is chosen as $p^*$ defined
	in \eqref{eq:p*def}, the normalized lower bound converges to
	$Z_N(\bfm)$ when the system grows in both demand and capacity. Thus,
	$\pi^{p^*}$ is asymptotically optimal.
	
	Our findings build upon the early classic results in the revenue
	management literature, which show that allocation policies arising
	from a single linear program make the normalized total expected
	revenue converge to an upper bound on the optimal value
	\citep{cooper2002asymptotic}. Our results are different and new in
	several important aspects. The model in \cite{cooper2002asymptotic}
	can designate/allocate a particular product (slot) type upon a
	customer arrival (because he assumes that customer preference
	is known upon arrival), while our model offers multiple product (slot)
	types for customers to choose from (because the customer preference is not known). Using the offer set as a decision in the
	model creates significant new challenges. 
	First of all, our
	fluid model formulation needs to explicitly take care of customer
	choice processes and is much more complicated than that in
	\cite{cooper2002asymptotic}. Leveraging the fluid model
	formulation, the asymptotically optimal policy in
	\cite{cooper2002asymptotic} accepts customer requests up to some
	customer type-specific thresholds, because the optimal solution in
	Cooper's fluid model prescribes such thresholds for each customer type. As a result, 
	Cooper's asymptotic policy leads to a closed-form expression for each type of the customer demand served, allowing him to directly show
	that the normalized demand served converges in distribution
	to a constant which matches the optimal fluid model decision.
	However, due to customer choices, our fluid model cannot give rise to
	such a simple policy. Our fluid model informs the optimal duration
	in which a particular offer set is used, and we use this information
	to construct our policy which has a completely different form compared to Cooper's policy. 
	%
	%
	%
	%
	%
	As we cannot control the exact product (slot) type in an offer set that
	will be chosen by an arriving customer, we do not have a closed-form
	expression as in \cite{cooper2002asymptotic} for the total demand that goes
	into each product (slot) type and eventually gets served. To deal with this difficulty, we construct a (very) tight lower bound on the objective value and show that this lower bound, after being normalized, converges to the optimal objective value of the fluid model. The idea of our proof may be useful to identify
	effective approximate policies in other capacity management contexts
	when the manager cannot directly control the product a customer may
	pick.

	\subsection{Constant Performance Guarantee of the Offering-all Policy}
	In this section, we \bo{focus on} a simple scheduling policy: the \textit{offering-all policy}. Let $\pi_0$
	represent this policy, so the offer set under $\pi_0$ is the full
	set of all slot types, irrespective of the period $n$. Note that the
	effective offer set at state $\mathbf{m}$ is $\bar{S}(\mathbf{m})$.
	That is, when customers arrive, they only consider those slot types
	with positive capacity when making a choice. We denote by
	$V_{n}^{\pi_0}(\mathbf{m})$ the expected fill count attained by
	applying the offering-all policy $\pi_0$ throughout.
	%
	%
	Indeed, this simple policy has a constant performance guarantee that
	states that for any set of parameters, the offering-all policy
	$\pi_0$ achieves at least half of the
	optimal fill count. 
	
	\begin{theorem}\label{2approximation}
		For any $\Omega$, $n$, and $\mathbf{m}$, $V_n(\mathbf{m})\leq 2
		V_{n}^{\pi_0}(\mathbf{m})$.
	\end{theorem}
	
	It is worth noting that Theorem~\ref{2approximation} in fact holds
	more broadly for all so-called \emph{myopic policies}, which at each
	period offer a set maximizing the expected number of filled slots
	for that period. Myopic policies, however, do not have to offer all
	slot types in all periods. For instance, offering slot types 1 and 3
	in the M model instance would constitute a myopic policy.

	Performance guarantee results on myopic policies exist in various
	dynamic optimization settings, and a ratio of 2 is often the best
	provable performance bound; see, e.g.,
	\cite{mehtareview,ChanFarias}. While this performance bound may seem
	a little loose, we shall see empirically in Section \ref{sec.offer.all.num} that the offering-all policy performs very well and much better than this lower bound; \nanr{in finite regimes, the offering-all policy also seems to perform better than the asymptotically optimal policy constructed in Section \ref{sec.asymp}.} 

\section{Sequential Offering}
\label{sec.seq.model}

 We now present our second scheduling paradigm, which allows the scheduler to offer multiple sets of slots sequentially. Recall that this way of offering slots may represent for instance web-based scheduling where available slots are not revealed simultaneously, as well as telephone-based scheduling. Intuitively, having the scheduler offer slots sequentially instead of all at once will be able to steer customers into selecting more favorable slots from the perspective of system optimization. The question we address in this section is then how many and what sets of slots to offer in order to maximize the fill rate. We start by introducing the sequential offering model next. 


\subsection{Model Outline}

Upon customer arrival, the scheduler chooses a $K$, $1 \le K \le
J$, and sequentially presents the customer with $K$ mutually
exclusive subsets $S_1$, $S_2$,$\dots$, $S_K \subseteq
\bar{S}(\mathbf{m})$. We denote this action as
$\mathbf{S}:=S_1-S_2-\dots-S_K$.
Denote by
$\mathscr{S}(\mathbf{m})$ the set of all possible such actions at
state $\mathbf{m}$, and by $I_k(\mathbf{S}):=\{i: \sum_{j\in S_k}
\Omega_{ij} \geq 1, i \notin \cup_{l=1}^{k-1} I_l(\mathbf{S})\}$,
$k=1,...,K$, the set of customer types who do not accept any slot
from the first $(k-1)$ offer sets but encounter at least one
acceptable slot in $S_k$. \nan{So $I_k(\mathbf{S})$ represents the set of customers who, given sequential offering $\mathbf{S}$, accept some slot upon arrival into the system. Moreover, the slot chosen by these customers belongs to the $k$th offer set $S_k$.} The probability that slot type $j$ is
chosen under action $\mathbf{S}$ may then be written as
$q_j(\mathbf{S}):= \sum_{k=1}^K \sum_{i \in I_k(\mathbf{S})}
\lambda_i q_{ij}(S_k)$, with $q_{ij}(\cdot)$ as
in \eqref{eqn.qij}. \nan{The assumption that the sets $S_1,S_2,\dots,S_K$ are mutually exclusive is made from a practical rather than mathematical standpoint: there is simply no reason to offer the same slot type in two or more sets, because the customer will book a slot as soon as she is offered a set with at least one acceptable slot. Thus only the first set in which such a slot is included is relevant.}

For ease of presentation, we still use $V_n(\bfm)$ to denote the
expected maximum number of slots that can be booked with $\bfm$
slots available and $n$ periods to go in this section. For an action
$\mathbf{S}$, we let $\bigcup \mathbf{S} := \bigcup_{i=1}^{K} S_i$
denote the set of all slot types offered throughout action $\mathbf{S}$. Then, for the sequential offering model, we have
\begin{equation}\label{eqn.bellman.obj.seq}
V_{n}(\mathbf{m}) = \max_{\mathbf{S}\in
	\mathscr{S}(\mathbf{m})}\left[\sum_{j\in \bigcup
	\mathbf{S}}q_j(\mathbf{S})(1-\Delta_{n-1}^j(\mathbf{m}))\right]+V_{n-1}(\mathbf{m}),\mbox{~~for~}n=N,N-1,...,1,
\end{equation}
where $V_0(\cdot)=0$ and
$\Delta_{n-1}^j(\mathbf{m}):=V_{n-1}(\mathbf{m})-V_{n-1}(\mathbf{m}-\mathbf{e}_j)$
denotes the marginal benefit due to the $m_j$th unit of slot type
$j$ at period $n-1$. We observe that both the transition
probability $q_j(\mathbf{S})$ and the set of feasible actions
$\mathscr{S}(\mathbf{m})$ are much more complicated than their
counterparts in the non-sequential model.

The sequential offering setting can be viewed as a generalization of
non-sequential scheduling to any number $K \ge 1$ of offer sets.
Consequently, it stands to reason that the offering-all policy will
not perform well in the sequential setting, as this would limit the
scheduler to a single offer set ($K=1$). We indeed numerically
confirm this conjecture in Section~\ref{sec.seq.num}. Note that, in
contrast to the non-sequential case, an offering-all policy is
unlikely to be used in a practical setting such as 
telephone scheduling (because it would take too much time for the scheduler to
go over every possible appointment option). 
\nan{In the online setting, there is a way to take advantage of sequential
	offerings by redesigning the customer interface that releases information sequentially.} 

\nan{To provide a roadmap of analyzing the sequential model, we summarize our key findings in this section as follows. 
	\begin{itemize}
		\item We first consider a general setting and derive various structural results that provide more insights; in particular, we show that it is optimal to offer slot types one by one. 
		\item For a large class of problem instances with nested preference structures (to be discussed later), we derive a \textit{closed-form} optimal sequential offering policy.
		\item For problem instances not in this class, we develop a simple and highly effective heuristic based on the idea of balanced resource use and fluid models. 
		\item We prove that the optimal sequential offering does as well as in the \nanr{non-sequential} case where the scheduler has full information on the customer type upon arrival; we argue that this equivalence allows us to apply the idea of sequential offering in various interactive scheduling contexts. 
\end{itemize}}

\subsection{Results for the General Sequential Offering Model}

We now present some properties of the sequential model with general
choice matrices. First, we derive some structural properties of the
value function.

\begin{lemma}
	\label{lemma.Vn.monotone.seq}
	The value function $V_n(\mathbf{m})$ satisfies: \\
	(i) $0 \leq V_{n+1}(\mathbf{m}) - V_n(\mathbf{m}) \leq 1$, $\forall \bfm \geq 0$, $\forall n=1,2,\dots,N-1$; and \\
	(ii) $0 \leq V_n(\mathbf{m} + \bfe_j) - V_n(\mathbf{m}) \leq 1$,
	$\forall \bfm \geq 0$, $\forall n=1,2,\dots,N$.
\end{lemma}

Part (ii) of Lemma~\ref{lemma.Vn.monotone.seq} implies that
$V_n(\mathbf{m} + \bfe_j) \leq V_n(\mathbf{m}) + 1$, i.e., it is
better to have a slot booked now rather than saving it for future.
Therefore, in the context of sequential offering, it is better to
keep offering slots if none has been taken so far.
This is formalized in the following result, which shows that there exists an optimal sequential
offering policy that exhausts all available slot types in each
period.

%

\begin{lemma}
	\label{lemma.fullset.seq} 
	For any $\Omega$, $\bfm$, and $n$, there exists an optimal action $\mathbf{S}^*$  such that $\bigcup \mathbf{S}^* =
	\bar{S}(\bfm)$.
	%
\end{lemma}
%
%

Building upon Lemma \ref{lemma.fullset.seq}, we are able to characterize the
structure of an optimal sequential offering policy, described in
the theorem below.


\begin{theorem}
	\label{lemma.one.by.one.seq} Let $\mathbf{m}>0$ be the system state
	at period $n \geq 1$, and let $j_1, j_2, ..., j_J$ be a permutation
	of $1,2,...,J$ such that $V_{n-1}(\mathbf{m} - \bfe_{j_k}) \geq
	V_{n-1}(\mathbf{m} - \bfe_{j_{k+1}})$, $k=1,2,\dots,J-1$. Then the
	action $\{j_1\}-\dots-\{j_J\}$ is optimal.
\end{theorem}
Theorem \ref{lemma.one.by.one.seq} implies that there exists an
optimal policy that offers one slot type at a time. 
%
More importantly, this result shows a specific optimal offer
sequence based on the value function to go.  To understand this,
recall that $V_{n-1}(\mathbf{m}) - V_{n-1}(\mathbf{m} - \bfe_{j})$
can be viewed as the value of keeping the $m_j$th type $j$ slot from
period $n-1$ onwards. As all customers bring in the same amount of
reward, it benefits the system the most if an arrival customer can be
booked for the slot type with the least value to keep, i.e., the
slot type with the largest $V_{n-1}(\mathbf{m} - \bfe_{j})$.

Even if the scheduler does not know the exact customer type,
following the optimal offer sequence described in Theorem
\ref{lemma.one.by.one.seq} ensures that the arriving customer takes
the ``least valuable'' slot (as long as there is at least one
acceptable slot remaining). Indeed, matching customers with slots in
this way would be the best choice for the scheduler, even if she had
perfect information about customer type, i.e., she knew exactly the
customer type upon arrival. Following this rationale, our next result
shows an interesting and important correspondence between (i) the sequential
offering without customer type information and (ii) the
non-sequential offering with \emph{perfect} customer type
information. To distinguish these two settings, we let
$V^s_n(\mathbf{m})$ and $V^f_n(\mathbf{m})$ represent the value
functions for settings (i) and (ii), respectively, in the next
theorem.

\begin{theorem}
	\label{theorem.equivalence} $V^s_n(\mathbf{m}) =
	V^f_n(\mathbf{m})$, $\forall~\mathbf{m} \geq 0,~n=0,1,2,\dots,N$.
\end{theorem}

Theorem \ref{theorem.equivalence} \nan{suggests that the optimal sequential offering can fully exploit the value of customer type information;} 
%
however, it does
\emph{not} imply that it 
can fully
elicit customer type. Specifically, optimal sequential offering
happens to result in the same system state changes as if the
scheduler had full information about customer type, but does not let
the scheduler know exactly the customer type (see Remark
\ref{remark.W.seq} in Section \ref{sec.beyond}). 
Theorem \ref{theorem.equivalence} suggests that sequential offering
is a useful operational mechanism to improve the scheduling
efficiency in the absence of customer type information. Our
numerical experiments in Section \ref{sec.num.results} confirm and
quantify such efficiency gains.

\subsection{Optimal Sequential Offering Policies}
\label{sec.opt.seq}


In this section we fully characterize the optimal sequential
offering policy for a large class of choice matrix instances, which
include the N, M and M+1 model instances (see Figure \ref{fig:three.models}).
%
%
%
To this end, let $I(j)$ be the set of customer types who accept slot
type $j$, i.e., $I(j)=\{i=1,2,\dots,I : \Omega_{ij}=1\}$, $\forall j
=1,2,\dots,J$.
It makes intuitive sense that if $I(j_1) \subset I(j_2)$, then slot
type $j_2$ is more valuable than $j_1$, and thus slot type $j_1$
should be offered first. Combining this observation with
Theorem~\ref{lemma.one.by.one.seq} could then help us to design an
optimal policy. Let us first introduce a specific class of
model instances.

\begin{definition}\label{def:nested}
	We say that a model instance characterized by $\Omega$ is {\em
		nested} if for all $j_1, j_2=1,2,\dots,J$ and $j_1 \neq j_2$, one
	of the following three conditions holds: (i) $I(j_1) \cap I(j_2) =
	\emptyset$, (ii) $I(j_1) \subset I(j_2)$, or (iii) $I(j_1) \supset
	I(j_2)$.
\end{definition}

Note that not all model instances are nested. One simple example
is the W model instance from Figure~\ref{fig:three.models}(b),
where $I(1) = \{1,2\}$ and $I(2) = \{2,3\}$. None of the
conditions (i)-(iii) from Definition~\ref{def:nested} hold in this
case for $j_1 = 1$ and $j_2 = 2$. However, it is readily verified
that the N, M and M+1 model instances are all nested.

\begin{remark}
	The concept of a nested model instance is related to the \emph{star
		structure} considered in the previous literature on flexibility
	design; see, e.g., \cite{akccay2010dynamic}. Consider a system with
	a certain number of resource types (corresponding to slot types in
	our context), which can be used to do jobs of certain types
	(customer types in our context). A star flexibility structure is one
	such that there are specialized resource types, one for each job
	type, plus a versatile resource type that can perform all job types.
	The nested structure generalizes the star structure.
\end{remark}

It turns out that we can fully characterize an optimal policy for
nested model instances as follows.

\begin{theorem}
	\label{thm.vn.order}
	%
	Suppose $\Omega$ is nested, any policy that offers slot type $j_1$
	before offering slot type $j_2$ for any $j_1,~j_2$ such that
	$I(j_1) \subset I(j_2)$ is optimal.
\end{theorem}
Theorem~\ref{thm.vn.order} proposes to offer nested slot types in an
increasing order of the accepting customer types. Note that when two
slot types are mutually exclusive (i.e., $I(j_1) \cap I(j_2) =
\emptyset$), the order in which they are offered is irrelevant,
since customers that would select a slot from one set could never
from the other. \nan{To give some specific examples,} we can fully
characterize the optimal policy for the N, M and M+1 model instances
using Theorem~\ref{thm.vn.order}.




\begin{corollary}
	\label{cor:N.seq} For the N model instance and any $n$ and $\bfm$, an optimal sequential offering
	policy is to offer $\mathbf{S} = \{1\}-\{2\}$.
\end{corollary}

\begin{corollary}
	\label{cor:MandMplus1.seq} For the M and M+1 model instances and any $n$, an optimal
	sequential offering policy is to offer
	$$
	\mathbf{S} = \left\{
	\begin{array}{ll}
	\{1,3\}-\{2\},  & {\rm ~if~} m_1, m_2, m_3 \ge 1,\\
	\{1\}-\{2\},    & {\rm ~if~} m_3=0,\\
	\{3\}-\{2\},    & {\rm ~if~} m_1=0.
	\end{array}
	\right.
	$$
\end{corollary}


\subsection{Beyond Nested Model Instances}
\label{sec.beyond}

While Theorem \ref{thm.vn.order} solves a large class of the
sequential model instances, not all instances have a nested
structure. In this section, we analyze the W
model instance (see Figure \ref{fig:three.models}) to glean some insights
into the instances which are not nested.

To analyze the W model instance, one can formulate an MDP
with three possible actions: $\{1,2\}$, $\{1\}-\{2\}$, and
$\{2\}-\{1\}$ (and the corresponding actions at the boundaries).
However, there exist no straightforward offering orders for slot
types, and the optimal sequential policy turns out to be state
dependent. Specifically, we find that the optimal policy is a
\emph{switching curve} policy: with the availability of one type of
slots held fixed, it is optimal to offer the other type of slots
first as long as there is a sufficiently large amount of such slots
left.

\begin{figure}[!h]
	\centering
	\includegraphics[width=0.7\textwidth]{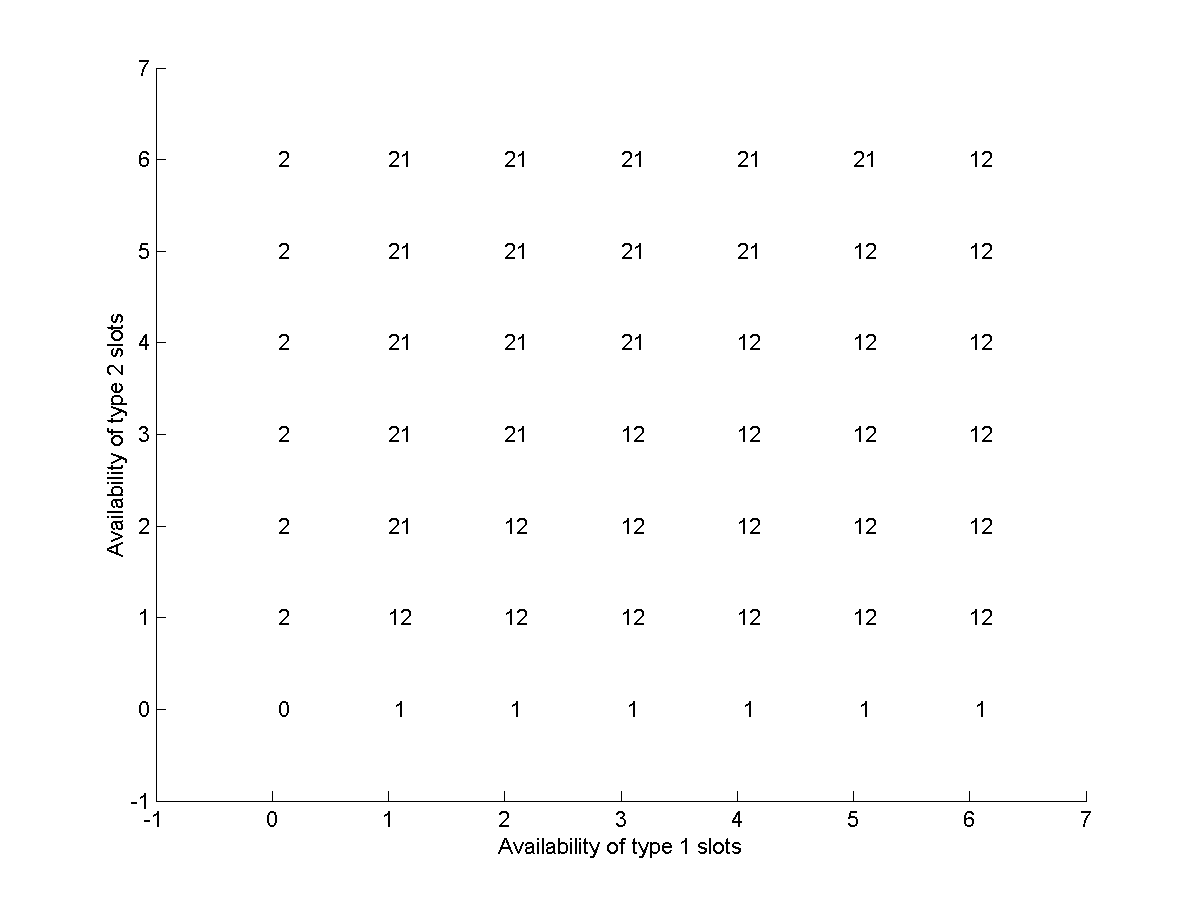}
	\caption{Structure of the optimal policy under W model instance with sequential offers.}
	\label{fig:WseqPolicy}
\end{figure}

Figure \ref{fig:WseqPolicy} illustrates the optimal actions for the
W model instance at different system states with $\mathbf{\lambda} =
(0.2, 0.5, 0,3)$ and $n=6$. The symbols ``0'', ``1'', ``2'', ``12''
and ``21'' correspond to the actions of offering nothing, offering
type 1 slots only, offering type 2 slots only, offering type 1 slots
and then type 2 slots, and offering type 2 slots and then type 1
slots, respectively. The optimal actions at boundary are obvious. In
the interior region of the system states, we can clearly see the
switching curve structure. For instance, when the system state is
$(3,3)$, it is optimal to offer $\{1\}-\{2\}$. When the number of
type 2 slots increase to 4, then it is optimal to offer
$\{2\}-\{1\}$.

The intuition behind this is different from that of the model instances
considered above where customer preferences are nested (e.g., the
N, M and M+1 model instances). In the W model instance, type 1 (3 resp.) customers
only accept type 1 (2 resp.) slots; but type 2 customers accept
both types of slots.
If there are relatively more type 1 slots than type 2 slots, then it
makes more sense to ``divert'' type 2 customers to choose type 1
slots, thus saving type 2 slots only for type 3 customers.
Accordingly, the switching curve policy stipulates that type 1 slots
to be offered first, ensuring that type 2 customers if any will pick
type 1 slots. The intuition above is formalized in the proposition
below.
\begin{proposition}
	\label{theorem.W.seq.opt.structure} Consider the W model instance with
	sequential offers. Given $m_2$, if there exists an $m_1^*$ such
	that the optimal action at state $(m_1^*,m_2)$ is $\{1\}-\{2\}$,
	then $\forall \bfm \in \{(m_1,m_2), m_1 \geq m_1^*\}$, the optimal
	action is $\{1\}-\{2\}$. Similarly, given $m_1$, if there exists
	an $m_2^*$ such that the optimal action at state $(m_1,m_2^*)$ is
	$\{2\}-\{1\}$, then $\forall \bfm \in \{(m_1,m_2), m_2 \geq
	m_2^*\}$, the optimal action is $\{2\}-\{1\}$.
\end{proposition}

\begin{remark}
	\label{remark.W.seq} In Section \ref{sec.opt.seq}, we state
	that sequential offering may not fully reveal exact customer types,
	but allows the system to evolve in the same optimal way as if the
	scheduler knew exactly the customer type. We use the W model instance
	to illustrate this point. Consider the W model instance with
	non-sequential offers and the scheduler knows the exact type of
	arriving customers. Suppose the optimal action is to offer $\{1\}$
	when type 1 or type 2 customers arrive; and to offer $\{2\}$ when
	type 3 customer arrives. Now, in a sequential offering model where
	the scheduler does not know the exact type of arriving customers, the
	scheduler would have offered $\{1\}$-$\{2\}$ to any arriving
	customer. If we encountered type 1 or 2 customers, type 1 slot would
	be taken, but we do not know the exact type of this customer (we know
	she must be either type 1 or type 2 though); if type 3 customer
	arrived, she would reject type 1 slot, but take type 2 slot. In this
	way, the system evolves as if the scheduler had perfect information
	on customer type.
\end{remark}


The structural properties of the optimal policy described in
Proposition \ref{theorem.W.seq.opt.structure} are likely the best we
can obtain for the W model instance; the exact form of the optimal
policy depends on model parameters and the system state, much like
with the M+1 model instance in the non-sequential case.
If customer preference structures become more
complicated, it is very difficult, if not impossible, to develop
structural properties for the optimal sequential offering policy.
Thus, for model instances that do not satisfy the conditions of Theorem
\ref{thm.vn.order}, we propose an effective heuristic below.

\subsection{The ``Drain'' Heuristic}
\label{sec:drain}

If customer preferences are not nested, the analysis of the W model
instance suggests that the optimal policy is to offer slots with
more capacity relative to its customer demand. Inspired by this
observation and using the idea of fluid models, we propose the
following heuristic algorithm which aims to ``drain'' the abundant
slot type first followed by less abundant ones. This heuristic aims
to have all slot types emptied simultaneously, thus maximizing the
fill rate. \nan{That is, this heuristic tries to ``balance'' the
	resource use.} Specifically, the drain algorithm works in the
following simple way. At period $n$ and for each slot type $j \in
\bar{S}(\bfm)$, we calculate
\begin{equation}
\label{eqn:drain} I_j := \frac{m_j}{n \sum_{i=1}^I \lambda_i
	\frac{\Omega_{ij}}{\sum_{k \in \bar{S}(\bfm)} \Omega_{ik}}}.
\end{equation}

Note that
$\frac{\Omega_{ij}}{\sum_{k=1}^{J} \Omega_{ik}}$ represents the
share of type $i$ customers who will choose type $j$ slots, assuming all available slot types are offered simultaneously. Taking
expectation with respect to the customer type distribution and
multiplying by $n$, the number of customers to come, the
denominator of \eqref{eqn:drain} can be viewed as the expected
load on type $j$ slots in the next $n$ periods. As a result, the
index $I_j$ can be regarded as the ratio between capacity left and
``expected'' load.

The drain algorithm is then to calculate all $I_j$s at the beginning
of each period, and to offer slots in decreasing order of the $I_j$.
The algorithm calls for offering slot types with larger $I_j$ first,
as these slot types have relatively more capacity compared to
demand. In other words, a slot type with a larger $I_j$ is likely to
have a smaller marginal value to keep, and thus can be offered
earlier. \nan{We could of course safely remove $n$ in the definition of $I_j$, and obtain the exact same order of slots. However, we leave $n$ in the denominator of \eqref{eqn:drain} because this allows us to interpret $I_j$ as the ratio between capacity left and ``expected'' number of requests. Based on this interpretation, it is clear that this heuristic aims to have all slot types emptied
	simultaneously, thus maximizing the fill rate.} We will test the
performance of this algorithm in
Section~\ref{sec:performance_drain}.

\subsection{Applications to Interactive Scheduling}
\label{sec.interactive}

\nanr{In Sections \ref{sec.non.seq.model} and \ref{sec.seq.model} we discuss two different models of customer-scheduler interactions %
	%
%
in the appointment booking practice. In one model, the scheduler makes a one-shot offering, and in the other, the scheduler enjoys the full flexibility of sequential offering.} \peterr{The appointment booking process, however, can fall in between these two models in terms of the degree to which the customer preference information is collected and used during the interaction between the scheduler and each customer.} 
%
%
%
%
\nan{Such interactions may be present both in a traditional setting with human interaction (e.g., a
	customer, after being offered an appointment at 8am by a receptionist, may
	indicate that none of the morning slots are acceptable) or fully digital (e.g., the Partners HealthCare Patient
	Gateway online booking website allows patients to indicate their
	acceptable time slots upfront).}

	
\nanr{When additional customer preference information is gathered during the appointment booking process, the scheduler can still follow the optimal list of slot types $\{j_1\}-\dots-\{j_J\}$ obtained from Theorem~\ref{lemma.one.by.one.seq}, but simply skip all slots known to be unacceptable, either up front or dynamically as additional information is collected. This offering strategy is still optimal because it would end up with the same system state compared to not skipping} \peterr{those slots indicated as unacceptable before or during the booking process}
\nanr{(e.g., directly declared by the customer) and thus give the exact same fill count that can be obtained if the scheduler had full information about the customer type (see Theorem \ref{theorem.equivalence}).}
	%
	%
%
\nan{Recall that the order of $\{j_1\}-\dots-\{j_J\}$
	can be readily obtained with nested customer preferences (Theorem \ref{thm.vn.order}), or otherwise an approximate order can
	be easily formed by the drain heuristic \eqref{eqn:drain}.}

\nan{Although outside the scope of this paper, these considerations on interactive scheduling raise various issues related to the tradeoff between obtaining the best fill rate and providing a convenient experience to the customer. For instance, the scheduler may want to limit the number of sets offered to the customer to provide a smooth user experience. In light of Theorem \ref{lemma.one.by.one.seq}, one potential idea for future study is to group slot types based on the order of $\{j_1\}-\dots-\{j_J\}$.}



%
%
%

\section{Numerical Results}
\label{sec.num.results}

\nan{In the last two sections, we consider non-sequential offering and sequential offering. For each setting, we derive optimal or near-optimal booking policies. 
In this section, we	run extensive numerical experiments to test and compare these policies and the two scheduling paradigms.}


\nan{We organize this section as follows. Section
	\ref{sec.offer.all.num} and Section \ref{sec:performance_drain} discuss the performance of the offering-all policy in the
	non-sequential model and of the drain heuristic in
	the sequential model, respectively. We demonstrate that
	these two algorithms obtain fill rates that are remarkably close to that of the respective optimal policies, and therefore can
	serve as simple, effective heuristics for practical use. Section \ref{sec.seq.num} compares the
	differences in the expected fill rate under the non-sequential and
	sequential offering models, where this difference represents the
	``value'' of sequential offering. Specifically, we evaluate the
	differences between the optimal policies and those between the
	heuristics. The former represent the ``theoretical''
	value of sequential offering compared to non-sequential offering,
	while the latter can be thought of as the ``practical'' value if
	practitioners adopt the heuristics mentioned above for
	each setting. Finally, Section \ref{sec.multi.day.num} extends
	our scheduling policies to a multi-day rolling horizon setting and demonstrates that the insights obtained from our analysis remain valid in this setting.}

\subsection{Performance of the Offering-all Policy}
\label{sec.offer.all.num}

We start our evaluation of the offering-all policy in two specific
model instances considered above: M and M+1 model instances. (We
need not to evaluate the offering-all policy in the N and W model
instances because the offering-all policy is optimal there.)
%
%
Here we use backward induction to determine the expected performance
of the optimal policy $\pi_1$, and compare it through simulation to
that of the offering-all policy $\pi_0$. To this end we simulate the offering-all policy for 1000 days. The performance metric of interest is the
percentage optimality gap defined as $(u_g-u_o)/u_o \times 100\%$,
%
%
where $u_o$ is the expected fill count of the optimal policy and
$u_g$ is the average fill count over 1000 simulated days under the
offering-all policy.

Table~\ref{tab:M} summarizes the statistics on the optimality gap of
the offering-all policy in the M model instance. For each
$N=20,30,40,50$, we evaluate the maximum, average and median
optimality gap over all possible initial capacity vectors
$(b_1,b_2,b_3) \in \mathbb{Z}_+^3$ such that $b_j \geq 0.2
N,~\forall j$ and $b_1+b_2+b_3=N$. The number of initial capacity
vectors considered for each $N$ is shown as the number of scenarios
in the second column of Table~\ref{tab:M}. In general, the
optimality gap of the offering-all policy in the M model instance is
relatively small ($\approx 3-4\%$) and is not sensitive to model
parameters.



\begin{table}[!ht]
	\caption{\small 
		Optimality gap of the offering-all policy in the M model instance.}
	\label{tab:M}
	\begin{center}
		{\small
			\begin{tabular}{c | c | ccc | ccc | ccc}
				\toprule
				\multirow{2}{*}{$N$}     &                   {\footnotesize $\#$ of}   & \multicolumn{3}{c|}{\footnotesize $(\lambda_1,\lambda_2) = (1/2, 1/2)$} & \multicolumn{3}{c|}{\footnotesize $(\lambda_1,\lambda_2) = (1/3, 2/3)$} & \multicolumn{3}{c}{\footnotesize $(\lambda_1,\lambda_2) = (1/4, 3/4)$} \\
				&           {\footnotesize Scenarios}            & \footnotesize Max & \footnotesize Average  & \footnotesize Median   & \footnotesize Max   &  \footnotesize Average  & \footnotesize Median   & \footnotesize Max &  \footnotesize Average  & \footnotesize Median  \\
				\otoprule
				$20$                                          & 45                     &  -4.4\% &  -3.6\%  &  -3.6\% &   -4.1\%  &  -3.3\% &  -3.3\%     & -3.6\%  &  -2.9\%  & -3.0\%                     \\
				$30$                                          & 91                     &  -4.8\% &  -3.7\%  &  -3.8\% &  -4.5\%   &  -3.5\%  & -3.5\%     & -3.8\%  &  -3.1\%  & -3.2\%             \\
				$40$                                          & 153                    &  -5.1\% &  -3.8\%  &  -3.8\% &  -4.7\%   &  -3.6\%   & -3.6\%     & -4.0\%  & -3.2\%  & -3.3\%         \\
				$50$                                           & 231                   &  -5.3\% &  -3.8\%  & -3.8\%  &  -4.8\%   &  -3.7\%   & -3.7\%     & -4.1\%  & -3.3\%  & -3.4\%       \\
				\bottomrule
			\end{tabular}
		}
	\end{center}
\end{table}


Table~\ref{tab:M+1nonseq} shows the optimality gap statistics for
the M+1 model instance, and the setup of this table is similar to
Table \ref{tab:M}.
%
When $\lambda_3$ is small, the M+1 model instance is very similar to
the M model and thus the optimality gaps of the offering-all policy
are similar to those observed in Table \ref{tab:M}.
%
%
As $\lambda_3$ increases the performance of the offering-all policy
improves, since offering-all becomes more likely to be optimal.

\begin{table}[!ht]
	\caption{\small Optimality gap of the offering-all policy in the M+1
		model instance.} \label{tab:M+1nonseq}
	\begin{center}
		{\small
			\begin{tabular}{c | c | ccc | ccc | ccc}
				\toprule
				\multirow{2}{*}{$N$}     &                   {\footnotesize $\#$ of}   & \multicolumn{3}{c|}{\footnotesize $(\lambda_1,\lambda_2,\lambda_3) = (9/20,9/20, 1/10)$} & \multicolumn{3}{c|}{\footnotesize $(\lambda_1,\lambda_2,\lambda_3) = (2/5,2/5, 1/5)$} & \multicolumn{3}{c}{\footnotesize $(\lambda_1,\lambda_2,\lambda_3) = (3/10,3/10, 2/5)$} \\
				&           {\footnotesize Scenarios}            & \footnotesize Max & \footnotesize Average  & \footnotesize Median   & \footnotesize Max   &  \footnotesize Average  & \footnotesize Median   & \footnotesize Max &  \footnotesize Average  & \footnotesize Median  \\
				\otoprule
				$20$                                          & 45                     &  -3.1\% &  -2.0\%  &  -1.9\% &   -2.0\%  &  -1.1\% &  -0.9\%     & -0.7\%  &  -0.3\%  & -0.2\%                     \\
				$30$                                          & 91                     &  -3.4\% &  -2.1\%  &  -2.0\% &  -2.3\%   &  -1.1\%  & -1.0\%     & -0.8\%  & -0.3\%  & -0.2\%             \\
				$40$                                          & 153                    &  -3.7\% &  -2.1\%  &  -2.0\% &  -2.5\%   &  -1.2\%   & -1.0\%     & -0.8\%  & -0.2\%  & -0.1\%         \\
				$50$                                           & 231                   &  -3.9\% &  -2.2\%  & -2.0\%  &  -2.6\%   & -1.2\%   & -1.0\%     & -0.8\%  & -0.2\%  & -0.1\%       \\
				\bottomrule
			\end{tabular}
		}
	\end{center}
\end{table}

	To evaluate the performance of the offering-all policy in settings beyond these two simple instances, we carry out an extensive numerical study using
	randomly generated customer preference matrices. Fixing the number of slot types $J$, there are $2^J$ different possible customer types, including those that accept no slots at all. By allowing any possible combination of these customer types, there could be $2^{2^J}-1$ possible preference matrices (excluding the empty matrix). 
	In order to test the performance of offering-all in a robust and yet computationally tractable manner, we compare its performance among many randomly generated such preference matrices. 
	
	We also vary the arrival probability vector $\veclambda = (\lambda_1,\dots,\lambda_I)$ for each preference matrix. In particular, we test three possible vectors: $\veclambda^{(1)}$ such that $\lambda^{(1)}_i = 1/I$, $\veclambda^{(2)}$ such that $\lambda^{(2)}_i = 2(I+i-2)/(3I^2-3I)$ and  $\veclambda^{(3)}$ with $\lambda^{(3)}_i = 2 (I+3i-4)/(5I^2-5I)$. In all three cases the $\lambda_i$'s add up to one. For $\veclambda^{(2)}$ and $\veclambda^{(3)}$, $\lambda_1$ is the largest, and each successive $\lambda_i$ is smaller by a factor 2 or 4, respectively. Note that the value of $I$ depends on the randomly generated preference matrix, and may vary from $I=1$ (since we exclude the empty matrix) to the maximum number of customer types.
	
	Our results are summarized in Table~\ref{tab:random}, where we show the optimality gap of the offering-all policy. We compute the performance of offering-all through simulation as before, and the performance of the optimal policy through backward induction. We fix $J$ and $N$, and then generate the number of random instances indicated in the table (`number of instances'). For each instance we also vary the initial capacity vectors similar to what was done for Tables~\ref{tab:M} and~\ref{tab:M+1nonseq} (`number of scenarios'). Fixing the structure of the arrival rate vector, we then report the maximum, average and median optimality gap over all instances and scenarios. It is clear from this table that the offering-all policy continues to do very well, and the average gap with the optimal policy is around $0.5\%$ throughout, independent of the size of the matrix and the arrival rates.
	
	\begin{table}[!ht]
		\caption{\small 
			Optimality gap of the offering-all policy for random network instances.}
		\label{tab:random}
		\begin{center}
			{\small
				\begin{tabular}{cc | cc | ccc | ccc | ccc}
					\toprule
					\multirow{2}{*}{$J$}    & \multirow{2}{*}{$N$}     &                   {\footnotesize $\#$ of}  & {\footnotesize $\#$ of}   & \multicolumn{3}{c|}{\footnotesize $\veclambda = \veclambda^{(1)}$} & \multicolumn{3}{c|}{\footnotesize $\veclambda = \veclambda^{(2)}$} & \multicolumn{3}{c}{\footnotesize $\veclambda = \veclambda^{(3)}$} \\
					&          &           {\footnotesize Instances}   &   {\footnotesize Scenarios}             & \footnotesize Max & \footnotesize Average  & \footnotesize Median   & \footnotesize Max   &  \footnotesize Average  & \footnotesize Median   & \footnotesize Max &  \footnotesize Average  & \footnotesize Median  \\
					\otoprule
					3    &   10    &  100  &	36    &	3.9\% &   	0.2\%  &   0.0\%   &	3.9\% &	0.2\% &	0.0\% &	3.7\% &	0.3\% &	0.0\% \\
					&   20	&  80  &	120   &	4.6\% &	0.3\% &	0.1\% &	4.6\% &	0.3\% &	0.1\% &	5.1\% &	0.4\% &	0.1\%  \\
					&   30	&  40  &	253   &	5.2\% &	0.3\% &	0.1\% &	5.1\% &	0.4\% &	0.1\% &	3.6\% &	0.3\% &	0.1\%  \\
					&   40	&  10  &	435   &	3.0\% &	0.3\% &	0.0\% &	4.0\% &	0.5\% &	0.1\% &	4.4\% &	0.5\% &	0.2\%  \\
					\hline
					4   &   10  &	100   &	84    &	5.0\% &	0.3\% &	0.1\% &	4.7\% &	0.4\% &	0.2\% &	4.0\% &	0.3\% &	0.1\% \\
					&   20  &	10    &	455   &	3.7\% &	0.5\% &	0.3\% &	2.8\% &	0.5\% &	0.3\% &	4.7\% &	0.5\% &	0.3\% \\
					&   30  &	10    &	83    &	3.6\% &	0.6\% &	0.3\% &	3.6\% &	0.7\% &	0.2\% &	3.0\% &	0.6\% &	0.4\% \\
					\hline
					5   &   10  &	100   &	126   &	3.9\% &	0.4\% &	0.2\% &	3.7\% &	0.4\% &	0.3\% &	4.2\% &	0.5\% &	0.3\% \\
					&   20  &	10    &	126   &	2.1\% &	0.3\% &	0.2\% &	3.9\% &	1.0\% &	0.6\% &	4.8\% &	0.6\% &	0.3\% \\
					\bottomrule
				\end{tabular}
			}
		\end{center}
	\end{table}


\nan{Before proceeding to the next section, we briefly discuss the
	performance of  $\pi^{p^*}$ (i.e., the static randomized policy
	arising from the fluid model in Section \ref{sec.asymp}). We focus
	on the M model and vary $N$, the arrival probabilities and the
	initial capacity vectors. We report the optimality gap statistics
	for $\pi^{p^*}$ in Table \ref{tab:M.random} in the Online Appendix.
	We observe that the average optimality gap decreases from about 8\%
	to 5\% when $N$ increases from 20 to 50. This is consistent with our
	theory above that $\pi^{p^*}$ is asymptotically optimal when the
	demand and capacity increase simultaneously. Due to space constraint, we
	shall refrain us from further exploring the computational issues of
	$\pi^{p^*}$ and leave those for future research.}


\subsection{Performance of the ``Drain'' Heuristic}\label{sec:performance_drain}

%
%

In this section, we evaluate the performance of our ``drain''
heuristics developed in Section \ref{sec:drain}. We focus on the N,
M and W model instances. As in Section~\ref{sec.offer.all.num}, we vary the mix of
customer types, the total number of periods and the initial capacity
vectors. The performance of the optimal sequential offering policy
is evaluated by backward induction. The performances of the drain
heuristic are evaluated by running a discrete event simulation with
1000 days replication and then computing the average fill count per
day. We present the statistics on the percentage optimality gaps of
drain in Tables \ref{tab:drain.N}, \ref{tab:drain.M} and
\ref{tab:drain.W} for the N, M and W instances, respectively. In
particular, for the N an W model instances, the optimality gap
statistics are taken over all initial capacity vectors $(b_1,b_2)$
such that $(b_1,b_2) \in \{(x,y)\in \mathbb{Z}_+^2: x,y \geq 0.2 N,
x+y = N\}$. The second column of each table shows the number of
initial capacity vectors consider for each $N$.


%

\begin{table}[h!]
	\caption{Optimality gap of the Drain Heuristic in the N model
		instance. \label{tab:drain.N}} {
		\begin{tabular}{ c | c | ccc | ccc | ccc}
			\toprule
			\multirow{2}{*}{$N$}     &                    {\footnotesize $\#$ of}   & \multicolumn{3}{c|}{\footnotesize $(\lambda_1, \lambda_2) = (1/2, 1/2)$} & \multicolumn{3}{c|}{\footnotesize $(\lambda_1, \lambda_2) = (1/3, 2/3)$} & \multicolumn{3}{c}{\footnotesize $(\lambda_1, \lambda_2) = (1/4, 3/4)$} \\
			&           {\footnotesize Scenarios}            & \footnotesize Max & \footnotesize Average  & \footnotesize Median   & \footnotesize Max   &  \footnotesize Average  & \footnotesize Median   & \footnotesize Max &  \footnotesize Average  & \footnotesize Median  \\
			\otoprule
			$20$                                          & 13                     &  -0.8\% &  -0.4\%  &  -0.4\% &  -0.6\% &  -0.1\%  &  -0.2\%     & -0.7\% &  -0.2\%  &  -0.3\%                    \\
			$30$                                          & 19                     &  -0.6\% &  -0.2\%  &  -0.4\% &  -0.8\% &  -0.2\%  &  -0.1\%     & -0.5\% &  -0.0\%  &  -0.1\%            \\
			$40$                                          & 25                    &  -0.6\%  &  -0.2\%  &  -0.2\% &  -0.8\%  &  -0.1\%  &  -0.1\%     & -0.5\%  &  -0.0\%  &  -0.0\%        \\
			$50$                                           & 31                   &  -0.5\%  &  -0.1\%  &  -0.2\%  &  -0.5\%  &  -0.2\%  &  -0.2\%     & -0.6\%  &  -0.0\%  &  -0.1\%      \\
			\bottomrule
		\end{tabular}
	}{} 
\end{table}

\begin{table}[h!]
	\caption{Optimality gap of the Drain Heuristic in the M model
		instance. \label{tab:drain.M}} {
		\begin{tabular}{ c | c | ccc | ccc | ccc}
			\toprule
			\multirow{2}{*}{$N$}     &                    {\footnotesize $\#$ of}   & \multicolumn{3}{c|}{\footnotesize $(\lambda_1, \lambda_2) = (1/2, 1/2)$} & \multicolumn{3}{c|}{\footnotesize $(\lambda_1, \lambda_2) = (1/3, 2/3)$} & \multicolumn{3}{c}{\footnotesize $(\lambda_1, \lambda_2) = (1/4, 3/4)$} \\
			&           {\footnotesize Scenarios}            & \footnotesize Max & \footnotesize Average  & \footnotesize Median   & \footnotesize Max   &  \footnotesize Average  & \footnotesize Median   & \footnotesize Max &  \footnotesize Average  & \footnotesize Median  \\
			\otoprule
			$20$                                          & 45                     &  -1.4\% &  -0.7\%  &  -0.8\% &  -1.1\% &  -0.4\%  &  -0.4\%     & -0.9\% &  -0.2\%  &  -0.2\%                    \\
			$30$                                          & 91                     &  -0.9\% &  -0.6\%  &  -0.6\% &  -0.8\% &  -0.3\%  &  -0.3\%     & -0.7\% &  -0.2\%  &  -0.2\%            \\
			$40$                                          & 153                    &  -0.7\%  &  -0.5\%  &  -0.5\% &  -0.7\%  &  -0.3\%  &  -0.3\%     & -0.9\%  &  -0.2\%  &  -0.1\%        \\
			$50$                                           & 231                   &  -0.6\%  &  -0.4\%  &  -0.5\%  &  -0.6\%  &  -0.2\%  &  -0.3\%     & -0.6\%  &  -0.1\%  &  -0.2\%      \\
			\bottomrule
		\end{tabular}
	}{} 
\end{table}

\begin{table}[h!]
	\caption{Optimality gap of the Drain Heuristic in the W model
		instance. \label{tab:drain.W}} {
		\begin{tabular}{c | c | ccc | ccc | ccc}
			\toprule
			\multirow{2}{*}{$N$}                   & {\footnotesize $\#$ of}   & \multicolumn{3}{c|}{\footnotesize $(\lambda_1, \lambda_2, \lambda_3) = (1/3, 1/3, 1/3)$} & \multicolumn{3}{c|}{\footnotesize $(\lambda_1, \lambda_2, \lambda_3) = (1/5, 1/2, 3/10)$} & \multicolumn{3}{c}{\footnotesize $(\lambda_1, \lambda_2, \lambda_3) = (1/10, 3/10, 3/5)$} \\
			&           {\footnotesize Scenarios}            & \footnotesize Max & \footnotesize Average  & \footnotesize Median   & \footnotesize Max   &  \footnotesize Average  & \footnotesize Median   & \footnotesize Max &  \footnotesize Average  & \footnotesize Median  \\
			\otoprule
			$20$                                          & 13                     &  -0.2\% &  0.0\%  &  0.1\% &  -0.1\% &  0.1\%  &  0.1\%     & -0.7\% &  -0.2\%  &  -0.1\%                     \\
			$30$                                          & 19                     &  -0.7\% &  0.0\%  &  0.0\% &  -0.4\% &  0.0\%  &  0.0\%    & -0.6\% &  -0.1\%  &  -0.1\%            \\
			$40$                                          & 25                    &  -0.2\%  &  0.0\%  &  0.0\% &  -0.2\%  & 0.0\%  & 0.0\%     & -0.5\%  &  0.1\%  &  0.0\%        \\
			$50$                                           & 31                   &  -0.2\%  &  0.0\%  &  0.0\%  &  -0.2\%  &  0.0\%  &  0.0\%     & -0.4\%  &  -0.1\%  &  -0.1\%      \\
			\bottomrule
		\end{tabular}
	}{} 
\end{table}

In the N model, the optimality gap of drain is on average within
0.4\% (max 0.8\%) in all 264 scenarios we tested. The performances
of drain in the W instance is slightly better than those in the N
model. For the M model instance, the optimality gap of drain is on
average within 0.7\% (max 1.4\%) across all 1560 scenarios we
tested.  These observations suggest that the drain heuristic has a
remarkable performance. Given its simplicity, it can serve as an
effective scheduling rule for practitioners.


\subsection{Value of Sequential Offering} \label{sec.seq.num}

\subsubsection{Comparison of Optimal Policies}


In this section, we investigate the value of sequential scheduling
by comparing the optimal sequential policy to the optimal
non-sequential policy. We focus on the N, M and W model instances.
To provide a robust performance evaluation, we vary a range of model
parameters, including the mix of customer types, the total number of
periods and the initial capacity vectors like in earlier sections.
%
%
%
Table~\ref{tab:Nseq} presents the maximum, average and median
percentage improvement in fill count by following an optimal
sequential offering policy compared to the optimal non-sequential
policy in the N model instance. Tables~\ref{tab:Mseq}
and~\ref{tab:Wseq} present the similar information for the M and W
model instances, respectively.


\begin{table}[h!]
	\caption{Fill Count Improvement in the N Model instance (Opt
		Sequential vs. Opt Non-sequential). \label{tab:Nseq}} {
		\begin{tabular}{c | c | ccc | ccc | ccc}
			\toprule
			\multirow{2}{*}{$N$}     &                    {\footnotesize $\#$ of}   & \multicolumn{3}{c|}{\footnotesize $(\lambda_1,\lambda_2) = (1/2, 1/2)$} & \multicolumn{3}{c|}{\footnotesize $(\lambda_1,\lambda_2) = (1/4, 3/4)$} & \multicolumn{3}{c}{\footnotesize $(\lambda_1,\lambda_2) = (3/4, 1/4)$} \\
			&           {\footnotesize Scenarios}            & \footnotesize Max & \footnotesize Average  & \footnotesize Median   & \footnotesize Max   &  \footnotesize Average  & \footnotesize Median   & \footnotesize Max &  \footnotesize Average  & \footnotesize Median  \\
			\otoprule
			$20$                                          & 13                     &  16.0\% &  10.6\%  &  12.4\% &  10.8\%  &  9.0\% &  9.6\%     & 13.2\%  &  6.2\%  & 5.5\%                     \\
			$30$                                          & 19                     &  16.8\% &  10.9\%  &  12.6\% &  11.1\%   &  9.3\%  & 9.6\%     & 14.0\%  &  6.3\%  & 5.4\%             \\
			$40$                                          & 25                    &  17.2\% &  11.1\%  &  12.5\% &  11.2\%   &  9.5\%   & 9.9\%     & 14.5\%  & 6.4\%  & 5.3\%         \\
			$50$                                           & 31                   &  17.5\% &  11.2\%  & 12.9\%  &  11.3\%   &  9.6\%   & 9.8\%     & 14.8\%  & 6.4\%  & 5.3\%       \\
			\bottomrule
		\end{tabular}
	}{}
\end{table}

\begin{table}[h!]
	\caption{Fill Count Improvement in the M Model instance (Opt
		Sequential vs. Opt Non-sequential). \label{tab:Mseq}} {
		\begin{tabular}{c | c | ccc | ccc | ccc}
			\toprule
			\multirow{2}{*}{$N$}     &                    {\footnotesize $\#$ of}   & \multicolumn{3}{c|}{\footnotesize $(\lambda_1,\lambda_2) = (1/2, 1/2)$} & \multicolumn{3}{c|}{\footnotesize $(\lambda_1,\lambda_2) = (1/3, 2/3)$} & \multicolumn{3}{c}{\footnotesize $(\lambda_1,\lambda_2) = (1/4, 3/4)$} \\
			&           {\footnotesize Scenarios}            & \footnotesize Max & \footnotesize Average  & \footnotesize Median   & \footnotesize Max   &  \footnotesize Average  & \footnotesize Median   & \footnotesize Max &  \footnotesize Average  & \footnotesize Median  \\
			\otoprule
			$20$                                          & 45                     &  7.4\% &  4.0\%  &  3.4\% &  7.8\%  &  4.2\% &  3.8\%     & 6.9\%  &  4.1\%  & 4.0\%                     \\
			$30$                                          & 91                     &  7.9\% &  3.7\%  &  3.3\% &  8.3\%   &  4.1\%  & 3.8\%     & 7.2\%  &  4.2\%  & 4.2\%             \\
			$40$                                          & 153                    &  8.3\% &  3.5\%  &  2.9\% &  8.5\%   & 4.1\%   & 3.9\%     & 7.4\%  & 4.3\%  & 4.3\%         \\
			$50$                                           & 231                   &  8.5\% &  3.3\%  &  2.6\% &  8.7\%   &  4.1\%   & 4.0\%     & 7.5\%  & 4.4\%  & 4.4\%       \\
			\bottomrule
		\end{tabular}
	}{} 
\end{table}

\begin{table}[h!]
	\caption{Fill Count Improvement in the W Model instance (Opt
		Sequential vs. Opt Non-sequential). \label{tab:Wseq}} {
		\begin{tabular}{c | c | ccc | ccc | ccc}
			\toprule
			\multirow{2}{*}{$N$}     &                    {\footnotesize $\#$ of}   & \multicolumn{3}{c|}{\footnotesize $(\lambda_1,\lambda_2,\lambda_3) = (1/3, 1/3, 1/3)$} & \multicolumn{3}{c|}{\footnotesize $(\lambda_1,\lambda_2,\lambda_3) = (1/5, 1/2, 3/10)$} & \multicolumn{3}{c}{\footnotesize $(\lambda_1,\lambda_2,\lambda_3) = (1/10, 3/10,3/5)$} \\
			&           {\footnotesize Scenarios}            & \footnotesize Max & \footnotesize Average  & \footnotesize Median   & \footnotesize Max   &  \footnotesize Average  & \footnotesize Median   & \footnotesize Max &  \footnotesize Average  & \footnotesize Median  \\
			\otoprule
			$20$                                          & 13                     &  8.2\% &  6.1\%  &  6.5\% &  10.8\%  &  6.6\% &  7.0\%     & 11.2\%  &  7.7\%  & 8.3\%                     \\
			$30$                                          & 19                     &  9.0\% &  6.6\%  &  7.9\% &  11.8\%   &  7.0\%  & 7.2\%     & 11.6\%  &  8.1\%  & 9.2\%             \\
			$40$                                          & 25                    &  9.5\% &  6.9\%  &  7.9\% &  12.3\%   &  7.2\%   & 7.3\%     & 12.0\%  & 8.3\%  & 9.4\%         \\
			$50$                                           & 31                   &  9.8\% &  7.1\%  & 8.3\%  &  12.7\%   &  7.3\%   & 7.3\%     & 12.2\%  & 8.4\%  & 9.4\%       \\
			\bottomrule
		\end{tabular}
	}{} 
\end{table}

We observe that the efficiency gains in the M and W model instances
are robust to the initial customer type mix. The efficiency gain in
the W model instance is about 6-7\% on average, and can be as high
as 13\%. The efficiency gain in the M model instance is slightly
lower. For the N model instance, the gain is relatively more
sensitive to customer type mix, and ranges between 6-11\% on average.
In certain cases, the efficiency gain in the N model can be as high
as 18\%. These numerical findings show that sequential offering
holds strong potentials to improve the operational efficiency in
appointment scheduling systems.

\subsubsection{Comparison of Heuristics}

In this section, we compare the performances of two heuristic
scheduling policies discussed above: the offering-all policy and the
``drain'' heuristic developed in Section \ref{sec:drain}. We also
consider another policy called the \emph{random sequential offering
	policy}, which offers available slot types one at a time in a
permutation chosen uniformly at random. This policy mimics the
existing practice of telephone scheduling, which is often done
without careful planning. These three policies are used or can be
easily used by practice, and therefore the comparison results in
this section reveal the value of sequential scheduling that may be
realized by adopting these policies in practice.

We focus on the N, M and W model instances, and use the combinations
of parameters as in earlier sections. The performance of these three
policies are evaluated by running a discrete event simulation with
1000 days replication and then computing the average fill count per
day for each policy. We present the percentage improvement in the
fill count of drain over the other two policies. Detailed results
are shown in Tables \ref{tab:drain.N2}, \ref{tab:drain.M2} and
\ref{tab:drain.W2}.

For the N model instance, we see an average 9-11\% improvement (max
18\%) if using drain compared to using random sequential or
offering-all. In the M model, the average improvement is around
7-8\% with max 14\%. For the W model, drain makes on average 6-8\%
improvement over random sequential or offering-all with the maximum
improvement up to 13\%. It is worth remarking upon that in all model
instances the random sequential policy has about the same
performance as offering-all. So although the former is a sequential
policy and the latter is not, the potential of sequential offering
is not exploited due to the careless choice of the offered slots.

\begin{table}[h!]
	\caption{Comparison of the Drain Heuristic with Other Scheduling
		Policies (The N Model instance). \label{tab:drain.N2}} {
		\begin{tabular}{c | c | c | ccc | ccc | ccc}
			\toprule
			& \multirow{2}{*}{$N$}     &                    {\footnotesize $\#$ of}   & \multicolumn{3}{c|}{\footnotesize $(\lambda_1,\lambda_2) = (1/2, 1/2)$} & \multicolumn{3}{c|}{\footnotesize $(\lambda_1,\lambda_2) = (1/3, 2/3)$} & \multicolumn{3}{c}{\footnotesize $(\lambda_1,\lambda_2) = (1/4, 3/4)$} \\
			&                         &           {\footnotesize Scenarios}            & \footnotesize Max & \footnotesize Average  & \footnotesize Median   & \footnotesize Max   &  \footnotesize Average  & \footnotesize Median   & \footnotesize Max &  \footnotesize Average  & \footnotesize Median  \\
			\otoprule
			&  $20$                                          & 13                     &  16.5\%  &  10.2\%  & 11.6\%  &  14.0\%  &  10.3\%  & 11.6\%     & 11.0\%  &  9.0\%  & 9.6\%                      \\
			{\footnotesize \% Imp. over} &  $30$              & 19                     &  17.1\%  &  10.8\%  & 12.4\%  &  14.0\%  &  10.6\%  & 11.4\%     & 11.7\%  &  9.1\%  & 9.7\%            \\
			{\footnotesize Offering-all} &  $40$              & 25                    &  17.8\%  & 10.9\%  & 12.5\%   &  13.9\%  & 10.8\%  & 11.4\%     &  11.1\%  & 9.3\%  & 9.6\%         \\
			&  $50$                                          & 31                    &  17.8\%  & 11.2\%  & 13.3\%   &  14.5\%  & 10.9\%  & 11.5\%      & 11.7\%  & 9.6\%  & 9.8\%      \\
			\otoprule
			&  $20$                                          & 13                     &  16.6\%  &  10.2\% &  11.9\%  &  14.0\%  &  10.3\% &  11.2\%     & 11.2\%  &  8.7\% &  9.0\%                      \\
			{\footnotesize \% Imp. over} &  $30$              & 19                     &  16.7\%   &  10.8\%  & 12.5\% &  13.9\%   &  10.5\%  & 11.4\%     & 11.8\%   &  9.3\%  & 9.5\%             \\
			{\footnotesize Random Sequential} &  $40$         & 25                    &  17.4\%   &  11.1\%   & 12.7\% &  14.0\%   &  11.0\%   & 11.9\%     & 11.9\%   &  9.5\%   & 9.6\%         \\
			&  $50$                                          & 31                    &  18.0\%   &  11.1\%   & 13.2\%  &  14.5\%   &  10.9\%   & 11.2\%     & 11.6\%   &  9.5\%   & 10.0\%       \\
			\bottomrule
		\end{tabular}
	}{} 
\end{table}

\begin{table}[h!]
	\caption{Comparison of the Drain Heuristic with Other Scheduling
		Policies (The M Model instance). \label{tab:drain.M2}} {
		\begin{tabular}{c | c | c | ccc | ccc | ccc}
			\toprule
			& \multirow{2}{*}{$N$}     &                    {\footnotesize $\#$ of}   & \multicolumn{3}{c|}{\footnotesize $(\lambda_1,\lambda_2) = (1/2, 1/2)$} & \multicolumn{3}{c|}{\footnotesize $(\lambda_1,\lambda_2) = (1/3, 2/3)$} & \multicolumn{3}{c}{\footnotesize $(\lambda_1,\lambda_2) = (1/4, 3/4)$} \\
			&                         &           {\footnotesize Scenarios}            & \footnotesize Max & \footnotesize Average  & \footnotesize Median   & \footnotesize Max   &  \footnotesize Average  & \footnotesize Median   & \footnotesize Max &  \footnotesize Average  & \footnotesize Median  \\
			\otoprule
			&  $20$                                          & 45                     &  12.1\%  &  7.0\%  & 6.8\%  &  11.8\%  &  7.4\%  & 7.1\%     & 10.8\%  &  6.9\%  & 6.9\%                      \\
			{\footnotesize \% Imp. over} &  $30$              & 91                     &  13.7\%  &  7.1\%  & 6.5\%  &  13.4\%  &  7.6\%  & 7.1\%     & 11.3\%  &  7.4\%  & 7.4\%            \\
			{\footnotesize Offering-all} &  $40$              & 153                    &  14.0\%  & 7.0\%  & 6.3\%   &  13.8\%  & 7.7\%  & 7.7\%     &  11.5\%  & 7.6\%  & 7.8\%         \\
			&  $50$                                          & 231                    &  13.9\%  & 7.0\%  & 6.4\%   &  14.0\%  & 7.8\%  & 7.9\%      & 11.6\%  & 7.9\%  & 8.0\%      \\
			\otoprule
			&  $20$                                          & 45                     &  12.3\%  &  7.0\% &  6.9\%  &  13.2\%  &  7.4\% &  6.8\%     & 11.1\%  &  6.9\% &  6.7\%                      \\
			{\footnotesize \% Imp. over} &  $30$              & 91                     &  13.5\%   &  7.1\%  & 6.7\% &  13.7\%   &  7.6\%  & 7.1\%     & 11.0\%   &  7.4\%  & 7.4\%             \\
			{\footnotesize Random Sequential} &  $40$         & 153                    &  13.8\%   &  7.0\%   & 6.4\% &  13.7\%   &  7.7\%   & 7.4\%     & 11.4\%   &  7.7\%   & 7.9\%         \\
			&  $50$                                          & 231                    &  14.2\%   &  7.0\%   & 6.4\%  &  14.2\%   &  7.9\%   & 7.7\%     & 11.8\%   &  7.8\%   & 8.0\%       \\
			\bottomrule
		\end{tabular}
	}{} 
\end{table}

\begin{table}[h!]
	\caption{Comparison of the Drain Heuristic with Other Scheduling
		Policies (The W Model instance). \label{tab:drain.W2}} {
		\begin{tabular}{c | c | c | ccc | ccc | ccc}
			\toprule
			& \multirow{2}{*}{$N$}     &                    {\footnotesize $\#$ of}   & \multicolumn{3}{c|}{\footnotesize $(\lambda_1,\lambda_2,\lambda_3) = (1/3, 1/3, 1/3)$} & \multicolumn{3}{c|}{\footnotesize $(\lambda_1,\lambda_2,\lambda_3) = (1/5, 1/2, 3/10)$} & \multicolumn{3}{c}{\footnotesize $(\lambda_1,\lambda_2,\lambda_3) = (1/10, 3/10,3/5)$} \\
			&                          &           {\footnotesize Scenarios}            & \footnotesize Max & \footnotesize Average  & \footnotesize Median   & \footnotesize Max   &  \footnotesize Average  & \footnotesize Median   & \footnotesize Max &  \footnotesize Average  & \footnotesize Median  \\
			\otoprule
			& $20$                                          & 13                     &  8.0\%  &  6.1\%  & 6.9\% &  10.8\%  &  6.6\%  & 6.7\%     & 11.6\%  &  7.8\%  & 8.5\%                     \\
			{\footnotesize \% Imp. over} & $30$             & 19                     &  9.1\%  &  6.5\%  & 7.6\% &  11.7\%  &  6.9\%  & 6.9\%     & 11.9\%  &  8.1\%  & 9.1\%             \\
			{\footnotesize Offering-all} & $40$             & 25                    &  9.7\%  & 6.9\%  & 7.8\%   &  12.5\%  & 7.2\%  & 7.2\%      & 12.3\%  & 8.4\%  & 9.5\%         \\
			& $50$                                          & 31                   &  10.3\%  & 7.0\%  & 7.9\%   &  12.6\%  & 7.3\%  & 7.3\%     & 12.4\%  & 8.4\%  & 9.3\%       \\
			\otoprule
			& $20$                                          & 13                     &  8.5\%  &  6.1\% &  6.8\% &  10.9\%  &  6.7\% &  6.9\%     & 10.9\%  &  7.5\% &  8.1\%                     \\
			{\footnotesize \% Imp. over} & $30$             & 19                     &  9.6\%   &  6.5\%  & 7.9\%  &  11.9\%   &  6.9\%  & 7.3\%     & 11.5\%   &  8.0\%  & 9.3\%             \\
			{\footnotesize Random Sequential} & $40$        & 25                    &  9.8\%   &  7.0\%   & 8.1\% &  12.5\%   &  7.2\%   & 7.5\%     & 12.3\%   &  8.4\%   & 9.4\%         \\
			& $50$                                          & 31                   &  10.2\%   &  7.1\%   & 7.9\%  &  12.8\%   &  7.3\%   & 7.4\%     & 12.2\%   &  8.3\%   & 9.4\%      \\
			\bottomrule
		\end{tabular}
	}{} 
\end{table}

\subsection{Simulation of a Multi-day Setting}
\label{sec.multi.day.num}

Our scheduling policy is based on a model that looks at how
appointment slots are depleted in a single day, and implicitly assumes that
customer demand to a single day is independent from other days. In
practice, customer demand for different days may be correlated
because customers who do not find an acceptable slot in one day may
opt for another day. To incorporate this effect, we develop a
simulation model to evaluate the potential benefits of using our
scheduling policies in a multi-day rolling-horizon setting.

Specifically, we assume that the number of daily customer arrivals
is either deterministic $N$ or a Poisson random variable with mean
$N$. Daily capacity of the service provider is $N$ slots. We
consider an M model for within-day preferences.  That is, each
customer will either be type 1 or 2, and there are three slot types
in each day. Suppose that the scheduling window is $T$ days, i.e.,
customers are allowed to make appointments $T$ days ahead. Upon each
customer's arrival, she has $D \leq T$ acceptable days, and these
$D$ acceptable days are randomly generated within the scheduling
window. (For example, if $T=10$ and $D=3$ then one customer may accept day 3, 5
and 6 from her arrival day.) The customer will then ask the provider
for potential slots in each of these $D$ days (one day at a time in
a random order). The provider offers slots following one of
the three scheduling policies discussed above: offering-all,
non-sequential optimal (blocking type 2 if available), and
sequential optimal in the M model instance. If the customer finds an
acceptable slot in a day, the customer will take it and the
scheduling is done for this customer; if the customer cannot find
acceptable slots in all $D$ acceptable days, she will leave and not book the
appointment.

In our experiments, we fix $T=15$ and $N=30$;\footnote{\nan{We also try
	other $N$'s in our numerical experiments, and observe that the value
	of $N$ has only marginal impact on the numerical results.}} we vary
the arrival probabilities and the initial capacity vectors in each
day (similar to Table \ref{tab:drain.M2}). We also vary $D=1,2,3,4$
to study the impact of customer flexibility in their choices of days
(a larger $D$ implies that customers are more flexible in their
choices). We run simulations for 1200 days, and use the first 200
days as warm-up periods.  Based on the results of the last 1000
days, we calculate the percentage improvement, if any, in the slot
fill count for non-sequential optimal and sequential optimal against
offering-all for each combination of parameters. For each $D$ and
the arrival probability vector, we report the max, mean and median
percentage improvement over the initial capacity vectors we
consider.

Table \ref{tab:multiday} shows the comparison results with deterministic daily
arrivals. (Results when daily arrivals are Poisson
random variables are similar; see Table \ref{tab:multiday.poi}
in the Online Appendix.) We observe that the optimal non-sequential
and sequential scheduling policies obtained in our single-day model
still bring sizable benefits to the multi-day scheduling setting we
consider. Consistent with earlier findings, sequential
offering brings much higher efficiency gains compared to non-sequential
offering. The maximum improvement in fill count by sequential
offering compared to offering-all can be as high as 11\%.
We also observe that when customers become more flexible in their
day choices (i.e., when $D$ increases), the benefits due to
``smart'' scheduling decrease. This can be explained by that when customers
are more flexible, their preferences are immaterial and
thus taking customers' preferences into account when making scheduling decisions
becomes less valuable.



\begin{table}[h!]
	\caption{Policy Comparison in a Multi-day Scheduling Setting.
		\label{tab:multiday}} {
		\begin{tabular}{c | c | c | ccc | ccc | ccc}
			\toprule
			& \multirow{2}{*}{$D$}     &                    {\footnotesize $\#$ of}   & \multicolumn{3}{c|}{\footnotesize $(\lambda_1,\lambda_2) = (1/2, 1/2)$} & \multicolumn{3}{c|}{\footnotesize $(\lambda_1,\lambda_2) = (1/3, 2/3)$} & \multicolumn{3}{c}{\footnotesize $(\lambda_1,\lambda_2) = (1/4, 3/4)$} \\
			&                         &           {\footnotesize Scenarios}            & \footnotesize Max & \footnotesize Average  & \footnotesize Median   & \footnotesize Max   &  \footnotesize Average  & \footnotesize Median   & \footnotesize Max &  \footnotesize Average  & \footnotesize Median  \\
			\otoprule
			&  $1$                                          & 45                     &  3.9\%  &  2.1\%  & 2.0\%  & 4.3\%  &  2.4\%  & 2.2\%     & 3.7\%  &  2.4\%  & 2.5\%                      \\
			{\footnotesize Non-sequential Optimal} &  $2$     & 91                    &  3.7\%  &  1.9\%  & 1.7\%  &  3.3\%  &  2.0\%  & 2.2\%     & 2.8\%  &  2.0\%  & 2.1\%            \\
			{\footnotesize vs. Offering-all}  & $3$              & 153                &  3.3\%  &  1.6\%  & 1.7\%  &  2.7\%  &  1.7\%  & 1.9\%     & 2.5\%  &  1.6\%  & 1.7\%          \\
			&  $4$                                          & 231                    &  2.7\%  & 1.3\%  & 1.3\%   &  2.3\%  & 1.4\%  & 1.5\%      & 2.0\%  & 1.4\%  & 1.4\%      \\
			\otoprule
			&  $1$                                          & 45                     &  9.5\%  &  4.0\% &  3.3\%  &  11.0\%  &  5.1\% &  4.0\%     & 9.8\%  &  5.5\% &  5.6\%                      \\
			{\footnotesize Sequential Optimal} &  $2$         & 91                    &  9.8\%   &  3.8\%  & 2.7\% &  9.5\%   &  4.5\%  & 4.6\%     & 7.6\%   &  4.9\%  & 5.2\%             \\
			{\footnotesize vs. Offering-all} &  $3$         & 153                     &  8.9\%   &  3.4\%  & 2.4\% &  7.7\%   &  3.8\%  & 4.4\%     & 6.1\%   &  4.0\%  & 4.2\%         \\
			&  $4$                                          & 231                    &  7.6\%   &  2.8\%   & 1.8\%  &  6.6\%   &  3.1\%   & 3.5\%     & 5.2\%   &  3.3\%   & 3.4\%       \\
			\bottomrule
		\end{tabular}
	}{} 
\end{table}

\section{Conclusion}
\label{sec.conclusion}

\nan{Motivated by the increasing popularity of online appointment
booking platforms, we study how to offer appointment slots to
customers in order to maximize the total number of slots filled. \nanr{We consider two models, non-sequential offering and sequential offering, for different customer-scheduler interactions in the appointment
booking process. For each model, we develop optimal or near-optimal booking policies.}} 


\peterr{In our numerical experiments, we find that sequential offering in a proper manner makes a significant improvement over the two benchmark policies: random sequential offering policy (which mimics the existing practice of telephone scheduling) and the offering-all policy (that resembles many of the current online appointment booking systems). This finding suggests substantial potentials for improving the current appointment scheduling practice.}


\nan{Another notable observation from our numerical study is
that the two benchmark policies have quite similar performances, which indicates that current online scheduling (that often offers all
available slots) and traditional telephone scheduling (without a
careful offer sequence) would result in similar fill rates. Thus,
one should \emph{not} expect that implementing an online scheduling
system in place of traditional telephone scheduling can
automatically lead to more appointments booked. However, as our
research suggests, one may improve the performance of online
scheduling by designing an interface that uses the idea of sequential offering, collecting information on customer choice
behavior and then making offers in a smarter way.} 



In summary, our work provides the first analytical framework to model,
compare and improve the appointment booking process. Our study also suggests many possible directions for future research. To name a few, 
\nan{first,  we assume a specific model for customer choice, and 
%
future research may consider scheduling decisions under different choice
models. Second, it would be interesting to
consider other customer behaviors (e.g., cancellations, no-shows,
recall, renege after a few trials) in the scheduling models. Third,
our numerical study of the multi-day scheduling is by no means
exhaustive and it would be a fruitful direction to investigate
the (optimal) joint offering policy for both day and slot choices. Last but not least, asymptotic regimes with different scalings of model parameters may be interesting objects of study both from a stochastic model theoretical perspective and for informing more efficient operations in practical settings.}

\bibliographystyle{ormsv080} 
\bibliography{myrefs} 

\begin{thebibliography}{19}
\expandafter\ifx\csname natexlab\endcsname\relax\def\natexlab#1{#1}\fi
\expandafter\ifx\csname url\endcsname\relax
  \def\url#1{{\tt #1}}\fi
\expandafter\ifx\csname urlprefix\endcsname\relax\def\urlprefix{URL }\fi
\expandafter\ifx\csname urlstyle\endcsname\relax
  \expandafter\ifx\csname doi\endcsname\relax
  \def\doi#1{doi:\discretionary{}{}{}#1}\fi \else
  \expandafter\ifx\csname doi\endcsname\relax
  \def\doi{doi:\discretionary{}{}{}\begingroup \urlstyle{rm}\Url}\fi \fi

\bibitem[{Ak{\c{c}}ay et~al.(2010)Ak{\c{c}}ay, Balakrishnan, and
  Xu}]{akccay2010dynamic}
Ak{\c{c}}ay, Y, A~Balakrishnan, SH~Xu. 2010.
\newblock Dynamic assignment of flexible service resources.
\newblock {\it Production and Operations Management\/} {\bf 19}(3) 279--304.

\bibitem[{Bernstein et~al.(2015)Bernstein, K\"ok, and Xie}]{bkx2015}
Bernstein, F, AG~K\"ok, L~Xie. 2015.
\newblock Dynamic assortment customization with limited inventories.
\newblock {\it Manufacturing \& Service Operations Management\/} {\bf 17}(4)
  538--553.

\bibitem[{Billingsley(1968)}]{billingsley1968convergence}
Billingsley, Patrick. 1968.
\newblock {\it Convergence of probability measures\/}.
\newblock John Wiley \& Sons.

\bibitem[{Cayirli and Veral(2003)}]{CayirliVeral03}
Cayirli, T, E~Veral. 2003.
\newblock Outpatient scheduling in health care: A review of literature.
\newblock {\it Production and Operations Management\/} {\bf 12}(4) 519--549.

\bibitem[{Chan and Farias(2009)}]{ChanFarias}
Chan, CW, VF~Farias. 2009.
\newblock Stochastic depletion problems: Effective myopic policies for a class
  of dynamic optimization problems.
\newblock {\it Mathematics of Operations Research\/} {\bf 34}(2) 333--350.

\bibitem[{Cooper(2002)}]{cooper2002asymptotic}
Cooper, William~L. 2002.
\newblock Asymptotic behavior of an allocation policy for revenue management.
\newblock {\it Operations Research\/} {\bf 50}(4) 720--727.

\bibitem[{Feldman et~al.(2014)Feldman, Liu, Topaloglu, and
  Ziya}]{feldman2014appointment}
Feldman, J, N~Liu, H~Topaloglu, S~Ziya. 2014.
\newblock Appointment scheduling under patient preference and no-show behavior.
\newblock {\it Operations Research\/} {\bf 62}(4) 794--811.

\bibitem[{Gallego et~al.(2016)Gallego, Li, Truong, and
  Wang}]{gallego2016online}
Gallego, G, A~Li, V-A Truong, X~Wang. 2016.
\newblock Online resrouce allocation with customer choice. Working paper,
  Department of Industrial Engineering and Operations Research, Columbia
  University.

\bibitem[{Golovin and Krause(2011)}]{golovin}
Golovin, D, A~Krause. 2011.
\newblock Adaptive submodularity: Theory and applications in active learning
  and stochastic optimization.
\newblock {\it Journal of Artificial Intelligence Research\/} {\bf 42}(1)
  427--486.

\bibitem[{Golrezaei et~al.(2014)Golrezaei, Nazerzadeh, and
  Rusmevichientong}]{golrezaei2014real}
Golrezaei, N, H~Nazerzadeh, P~Rusmevichientong. 2014.
\newblock Real-time optimization of personalized assortments.
\newblock {\it Management Science\/} {\bf 60}(6) 1532--1551.

\bibitem[{Green et~al.(2006)Green, Savin, and Wang}]{green2006managing}
Green, L, S~Savin, B~Wang. 2006.
\newblock Managing patient service in a diagnostic medical facility.
\newblock {\it Operations Research\/} {\bf 54}(1) 11--25.

\bibitem[{Gupta and Denton(2008)}]{GuptaDenton08}
Gupta, D, B~Denton. 2008.
\newblock {Appointment scheduling in health care: Challenges and
  opportunities}.
\newblock {\it IIE Transactions\/} {\bf 40}(9) 800--819.

\bibitem[{Gupta and Wang(2008)}]{GuptaWang08}
Gupta, D, L~Wang. 2008.
\newblock Revenue management for a primary-care clinic in the presence of
  patient choice.
\newblock {\it Operations Research\/} {\bf 56}(3) 576--592.

\bibitem[{Liu et~al.(2010)Liu, Ziya, and Kulkarni}]{LiuZiyaKulkarni08}
Liu, N, S.~Ziya, V.~G. Kulkarni. 2010.
\newblock Dynamic scheduling of outpatient appointments under patient no-shows
  and cancellations.
\newblock {\it Manufacturing \& Service Operations Management\/} {\bf 12}(2)
  347--364.

\bibitem[{Mehta(2013)}]{mehtareview}
Mehta, A. 2013.
\newblock Online matching and ad allocation.
\newblock {\it Foundations and Trends in Theoretical Computer Science\/} {\bf
  8}(4) 265--368.

\bibitem[{Subramanian et~al.(1999)Subramanian, Stidham~Jr, and
  Lautenbacher}]{subramanian1999airline}
Subramanian, J, S~Stidham~Jr, CJ~Lautenbacher. 1999.
\newblock Airline yield management with overbooking, cancellations, and
  no-shows.
\newblock {\it Transportation Science\/} {\bf 33}(2) 147--167.

\bibitem[{Talluri and Van~Ryzin(2004)}]{talluri2004revenue}
Talluri, KT, GJ~Van~Ryzin. 2004.
\newblock Revenue management under a general discrete choice model of consumer
  behavior.
\newblock {\it Management Science\/} {\bf 50}(1) 15--33.

\bibitem[{Wang and Gupta(2011)}]{wang2011adaptive}
Wang, WY, D~Gupta. 2011.
\newblock Adaptive appointment systems with patient preferences.
\newblock {\it Manufacturing \& Service Operations Management\/} {\bf 13}(3)
  373--389.

\bibitem[{Zhang and Cooper(2005)}]{zhang2005revenue}
Zhang, D, WL~Cooper. 2005.
\newblock Revenue management for parallel flights with customer-choice
  behavior.
\newblock {\it Operations Research\/} {\bf 53}(3) 415--431.

\end{thebibliography}

\newpage
\appendix
\section*{Proofs and Additional Numerical Results}

	\section{Proof of the Results in Section \ref{sec.non.seq.model}}
	\subsection{Preliminarily results} We first state and prove an
	auxiliary lemma on the structural results of the value function for
	the non-sequential offering model. This lemma will be used in
	proving other results in the paper.
	
	\begin{lemma}
		\label{lemma.Vn.monotone}
		Let $\Omega$ be a preference matrix, $\mbf \ge 0$, $j=1,\dots,J$ and $n \in \{1,\dots,N\}$, then the value function $V_n(\mathbf{m})$ satisfies \\
		(i) $0 \leq V_{n+1}(\mathbf{m}) - V_n(\mathbf{m}) \leq 1;\quad \forall n=0,1,2,\dots;$ \\
		(ii) $0 \leq V_n(\mathbf{m} + \mathbf{e}_j) - V_n(\mathbf{m}) \leq
		1;\quad \forall n=0,1,2,\dots;$\\
		(iii) if $\lambda_0 > 0$, then $V_n(\mathbf{m} + \mathbf{e}_j) -
		V_n(\mathbf{m}) < 1;\quad \forall n=1,2,\dots.$.
	\end{lemma}
	These monotonicity results are quite intuitive. Properties (i) and
	(ii) state that the optimal expected reward is increasing in the
	number of customers and the number of slots left and  the changes in
	the optimal expected reward are bounded by the changes in the number
	of customers to go and the number of slots available. Property (iii)
	suggests that if there is a strictly positive probability that no
	customers would come in each period, then the increase of the
	optimal expected reward is strictly smaller than that of the
	available slots.
	
	\proof{Proof.} We use induction to prove this lemma. We first
	prove the first two properties. For $n=0$, these two properties
	hold trivially. Suppose that they also hold up to $n=t$. Consider
	$n=t+1$. Let $\mathbf{g}^*_t(\mathbf{m})$ represent the optimal
	decision rule in period $t$ when the system state is $\mathbf{m}$.
	Let $V_s^\mathbf{f}(\mathbf{m})$ be the expected number of slots
	filled given that the decision rule $\mathbf{f}$ is taken at stage
	$s$ and from stage $s-1$ onwards the optimal decision rule is
	used. Let $p_k(\mathbf{m}, \mathbf{f})$ be the probability that a
	type $k$ slot is booked at state $\bfm$ if action $\mathbf{f}$ is
	taken. It follows that
	\begin{align*}
	V_{t+1}(\mathbf{m})  \geq  V_{t+1}^{\mathbf{g}^*_t(\mathbf{m})}(\mathbf{m})  =&    \sum_{k=0}^J p_k(\mathbf{m}, \mathbf{g}^*_t(\mathbf{m})) [\mathbbm{1}_{\{k>0\}} + V_{t}(\mathbf{m} - \mathbf{e}_{k})] \\
	\geq & \sum_{k=0}^J  p_k(\mathbf{m}, \mathbf{g}^*_t(\mathbf{m})) [\mathbbm{1}_{\{k>0\}} + V_{t-1}(\mathbf{m} - \mathbf{e}_{k})]  \\
	=    &  V_{t}(\mathbf{m}),
	\end{align*}
	where the first inequality is due to the definition of
	$V_{t+1}(\mathbf{m})$ and the second inequality follows from the
	induction hypothesis. Following a similar argument and fixing $j
	\in \{1,2,\dots, J\}$, we have
	\begin{align*}
	V_{t+1}(\mathbf{m}+e_j) \geq &  V_{t+1}^{\mathbf{g}^*_{t+1}(\mathbf{m})} (\mathbf{m}+e_j) \\
	=    & \sum_{k=0}^J p_k(\mathbf{m}+\mathbf{e}_j, \mathbf{g}^*_{t+1}(\mathbf{m})) [\mathbbm{1}_{\{k>0\}} + V_{t}(\mathbf{m} + \mathbf{e}_j - \mathbf{e}_{k})] \\
	=    & \sum_{k=0}^J p_k(\mathbf{m}, \mathbf{g}^*_{t+1}(\mathbf{m})) [\mathbbm{1}_{\{k>0\}} + V_{t}(\mathbf{m} + \mathbf{e}_j - \mathbf{e}_{k})] \\
	\geq & \sum_{k=0}^J p_k(\mathbf{m}, \mathbf{g}^*_{t+1}(\mathbf{m})) [\mathbbm{1}_{\{k>0\}} + V_{t}(\mathbf{m} - \mathbf{e}_{k})]  \\
	=    &  V_{t+1}(\mathbf{m}),
	\end{align*}
	where the second equality results from the decision rules and the
	state transition probability \eqref{eqn.dist.ln}.

	To show the RHS of the inequality in (i) for $n=t+1$, note that
	\begin{align*}
	& V_{t+1}(\mathbf{m}) - V_{t}(\mathbf{m}) \\
	= &    \sum_{k=0}^J p_k(\mathbf{m}, \mathbf{g}^*_{t+1}(\mathbf{m})) [\mathbbm{1}_{\{k>0\}} + V_{t}(\mathbf{m} - \mathbf{e}_{k})] - V_t(\mathbf{m}) \\
	= & \sum_{k=1}^J p_k(\mathbf{m}, \mathbf{g}^*_{t+1}(\mathbf{m})) + \sum_{k=0}^J p_k(\mathbf{m}, \mathbf{g}^*_{t+1}(\mathbf{m})) [ V_{t}(\mathbf{m} - \mathbf{e}_{k}) - V_{t}(\mathbf{m}) ]  \\
	\leq &  \sum_{k=1}^J p_k(\mathbf{m}, \mathbf{g}^*_{t+1}(\mathbf{m})) \\
	\leq & 1,
	\end{align*}
	where the first inequality follows from that $ V_{t}(\mathbf{m} -
	\mathbf{e}_{k}) \leq V_{t}(\mathbf{m})$, which has been shown
	above.
	
	To show the RHS of the inequality in (ii) for $n=t+1$, we define a
	decision rule ${\mathbf{h}}$ in period $t+1$ such that $\mathbf{h}
	= \mathbf{g}^*_{t+1}(\mathbf{m}+\mathbf{e}_j)$ except ${h}_j=0$.
	It follows that
	\begin{equation}
	\label{eqn.lemma.monotonicity.bd.slot}
	V_{t+1}(\mathbf{m}+\mathbf{e}_j) - V_{t+1}(\mathbf{m})  \leq
	V_{t+1}(\mathbf{m}+\mathbf{e}_j)
	-V_{t+1}^{{\mathbf{h}}}(\mathbf{m}),
	\end{equation}
	because ${\mathbf{h}}$ may not be the optimal given system state
	$\mathbf{m}$ at period $t+1$. For $u =1, 2, \dots, J$, let
	\begin{align*}
	q_u  = p_u (\mathbf{m}+\mathbf{e}_j,
	\mathbf{g}^*_{t+1}(\mathbf{m}+\mathbf{e}_j) )
	\end{align*}
	and
	\begin{align*}
	q'_u  & = p_u(\mathbf{m}, {\mathbf h} ). 
	\end{align*}
	It is easy to check that $q_u \leq q'_u, \forall u \neq 0, j$ and
	$q'_j=0$. Now, let $\Omega_{i \cdot} = (\Omega_{i1}, \Omega_{i2},
	\dots, \Omega_{iJ})$ and use $\langle \cdot,\cdot \rangle$ to
	represent the inner product. We have that
	\begin{align*}
	\sum_{u=1}^J q_u
	%
	=  \sum_{i=1}^I \lambda_i \mathbbm{1}_{ \{ \langle \Omega_{i \cdot} , \mathbf{g}_{t+1}(\mathbf{m}+\mathbf{e}_j) \rangle > 0\} } 
	\geq \sum_{i=1}^I \lambda_i \mathbbm{1}_{ \{ \{ \langle \Omega_{i
			\cdot} , \mathbf{h} \rangle  > 0\} } = \sum_{u=1}^J q'_u,
	\end{align*}
	because $\mathbf{h} = \mathbf{g}^*_{t+1}(\mathbf{m}+e_j)$ except
	${h}_j=0$.
	Therefore, $q_0 = 1- \sum_{u=1}^J q_u \leq 1- \sum_{u=1}^J q'_u =
	q'_0$. Define $\delta_u = q'_u - q_u$ for $u \neq j$. It is clear
	that $\delta_u \geq 0, \forall u \neq j$, and we note the
	following relationship.
	\[
	q_j = 1 - \sum_{u \neq j} q_u = 1 - \sum_{u \neq j} (q'_u - \delta_u) 
	= \sum_{u \neq j} \delta_u.
	\]
	
	Now, we can continue the inequality
	\eqref{eqn.lemma.monotonicity.bd.slot} as follows.
	\begin{align*}
	&    V_{t+1}(\mathbf{m}+\textbf{e}_j) -V_{t+1}^{{\mathbf{h}}}(\mathbf{m}) \\
	=  & \sum_{u=0}^J q_u [\mathbbm{1}_{\{u>0\}} + V_{t}(\mathbf{m} + \mathbf{e}_{j} - \mathbf{e}_u)] - \sum_{u=0}^J q'_u [\mathbbm{1}_{\{u>0\}} + V_{t}(\mathbf{m}  - \mathbf{e}_u)] \\
	=  & (1-q_0) - (1-q'_0) + \sum_{u=0}^J q_u V_{t}(\mathbf{m} + \mathbf{e}_{j} - \mathbf{e}_u) - \sum_{u=0}^J q'_u V_{t}(\mathbf{m}  - \mathbf{e}_u)] \\
	=  & \delta_0 + \sum_{u \neq j} q_u (V_{t}[\mathbf{m} + \mathbf{e}_{j} - \mathbf{e}_u) - V_{t}(\mathbf{m}  - \mathbf{e}_u)] + q_j V_t(\mathbf{m}) - \sum_{u \neq j} \delta_u V_{t}(\mathbf{m}  - \mathbf{e}_u) \\
	=  & \delta_0 + \sum_{u \neq j} q_u (V_{t}[\mathbf{m} + \mathbf{e}_{j} - \mathbf{e}_u) - V_{t}(\mathbf{m}  - \mathbf{e}_u)] + \sum_{u \neq j} \delta_u V_t(\mathbf{m}) - \sum_{u \neq j} \delta_u V_{t}(\mathbf{m}  - \mathbf{e}_u) \\
	=  & \delta_0 + \sum_{u \neq j} q_u (V_{t}[\mathbf{m} + \mathbf{e}_{j} - \mathbf{e}_u) - V_{t}(\mathbf{m}  - \mathbf{e}_u)] + \sum_{u \neq 0,j} \delta_u [V_t(\mathbf{m}) - V_{t}(\mathbf{m}  - \mathbf{e}_u)] \\
	\leq & \delta_0 + \sum_{u \neq j} q_u + \sum_{u \neq 0,j} \delta_u \\
	=    & 1,
	\end{align*}
	where the last inequality comes from the induction hypothesis for
	property (ii).
	
	As for property (iii), first note that it trivially holds for
	$n=1$. We can then follow similar induction steps as those used to
	prove the RHS of the inequality in property (ii) to complete the proof.
	
	\endproof

	\subsection{Proof of Proposition~\ref{lemma.W.opt.policy}}
	\proof{Proof.} We focus on the W model instance here, as the N Model
	instance is a special case of this. For $n \geq 1$ and any system
	state $(x, y) \geq (1, 1)$, the optimality equation for the W model
	instance reads.
	\begin{equation}
	\label{eqn.W.bellman} {V_n(x,y) = \max } \left\{
	\begin{array}{lll}
	1 -\lambda_0 + (\lambda_1 + \frac{1}{2}\lambda_2) V_{n-1}(x-1,y) + (\frac{1}{2}\lambda_2 + \lambda_3) V_{n-1}(x, y-1) + \lambda_0 V_{n-1}(x,y), \\ 
	(1-\lambda_3 - \lambda_0) + (\lambda_3 + \lambda_0) V_{n-1}(x,y) + (\lambda_1 + \lambda_2) V_{n-1}(x-1,y), \\ 
	(1-\lambda_1 - \lambda_0) + (\lambda_1 + \lambda_0) V_{n-1}(x,y) + (\lambda_2 + \lambda_3) V_{n-1}(x,y-1) 
	\end{array}
	\right\},
	\end{equation}
	where the three terms in the max operator correspond to the action
	of offering slot types $\{1,2\}$, $\{1\}$ and $\{2\}$, respectively.
	For the boundary conditions, it is easy to see that $V_0(x,y)=0$
	regardless of $x$ and $y$.  When one type of the slots are depleted,
	it is optimal to offer the other type of the slots. To calculate
	$V_n(x,0)$, note that type 1 slots are accepted only by type 1 and
	type 2 customers and the number of type 1 and type 2 customers in
	the last $n$ customers yet to come has a binomial distribution with
	parameters $n$ and $\lambda_1+\lambda_2$. Denote this random
	variable by $X_1 \sim \mathrm{Bin}(n,\lambda_1+\lambda_2)$.
	It follows that
	\begin{equation}
	\label{eq:Vnx0} V_n(x,0) = \mathbf{E}(\min\{x,
	X_1\}) = \sum_{k=0}^n \min(x,k) {n
		\choose k} (\lambda_1+\lambda_2)^k (1-\lambda_1-\lambda_2)^{n-k}.
	\end{equation}
	Similarly, with $X_2 \sim \mathrm{Bin}(n,\lambda_2+\lambda_3)$
	\begin{equation}
	\label{eq:Vn0y} V_n(0,y) = \mathbf{E}(\min\{y,
	X_2\}) = \sum_{k=0}^n \min(y,k) {n
		\choose k} (\lambda_2+\lambda_3)^k (1-\lambda_2-\lambda_3)^{n-k}.
	\end{equation}
	For ease of presentation, we define $\Delta_n^{ij}(x,y)$ to be the
	difference of the $i$th and $j$th terms in the max operator
	\eqref{eqn.W.bellman} above, $i,j \in \{1,2,3\}$.  In particular, we
	have
	\begin{equation}
	\label{eq:Delta12} \Delta^{12}_n(x,y) = \lambda_3 - \frac{1}{2}
	\lambda_2 V_{n-1}(x-1,y) + (\frac{1}{2}\lambda_2 + \lambda_3)
	V_{n-1}(x, y-1) -\lambda_3 V_{n-1}(x,y),
	\end{equation}
	and
	\begin{equation}
	\label{eq:Delta13} \Delta^{13}_n(x,y) = \lambda_1 - \frac{1}{2}
	\lambda_2 V_{n-1}(x,y-1) + (\frac{1}{2}\lambda_2 + \lambda_1)
	V_{n-1}(x-1, y) -\lambda_1 V_{n-1}(x,y).
	\end{equation}
	It suffices to show that $\Delta^{12}_n(x,y), \Delta^{13}_n(x,y)
	\geq 0$ for any $x,y \geq 1$ (the case when $x$ or $y$ equals 0 is
	trivial as it meets the boundary conditions  discussed above; see
	\eqref{eq:Vnx0} and \eqref{eq:Vn0y}). We use induction below to
	prove this. When $n=1$, it is a trivial proof as it is optimal to
	offer all available slots with one period left. Suppose that
	\eqref{eq:Delta12} and \eqref{eq:Delta13} hold up to $n=k$ and
	for any $x, y \geq 1$. 
	Now, consider $n=k+1$ and $x,y\geq 1$ . We have four cases to check:
	(1) $x=y=1$; (2) $y=1$ and $x \geq 2$; (3) $x=1$ and $y \geq 2$; and
	(4) $x,y \geq 2$. We start with case (1) and evaluate the term
	$\Delta^{13}_{k+1}(x,1)$ below.
	\begin{align*}
	\Delta^{13}_{k+1}(1,1) =  & \lambda_1  + (\lambda_1 + \frac{1}{2}\lambda_2) V_{k}(0, 1) - \frac{1}{2} \lambda_2 V_{k}(1,0) -\lambda_1 V_{k}(1,1)  \\
	=  & \lambda_1  + (\lambda_1 + \frac{1}{2}\lambda_2) [1-(\lambda_1 + \lambda_0) ^k] - \frac{1}{2} \lambda_2 [1- (\lambda_3 + \lambda_0)^k] \\
	&- \lambda_1 [1 - \lambda_0 + (\lambda_1 + \frac{1}{2} \lambda_2)V_{k-1}(0,1) + (\lambda_3 + \frac{1}{2} \lambda_2)V_{k-1}(1,0) + \lambda_0 V_{k-1}(1,1)]  \\
	=  & \lambda_0 \Delta^{13}_k(1,1)  + (\lambda_1 + \frac{1}{2}\lambda_2) [1-(\lambda_1 + \lambda_0) ^k] - \frac{1}{2} \lambda_2 [1- (\lambda_3 + \lambda_0)^k] \\
	& - \lambda_1 [(\lambda_1 + \frac{1}{2} \lambda_2)V_{k-1}(0,1) + (\lambda_3 + \frac{1}{2} \lambda_2)V_{k-1}(1,0)] \\
	& - \lambda_0 (\lambda_1 + \frac{1}{2}\lambda_2) V_{k-1}(0, 1) + \frac{1}{2} \lambda_0 \lambda_2 V_{k-1}(1,0)   \\
	=  & \lambda_0 \Delta^{13}_k(1,1)  + (\lambda_1 + \frac{1}{2} \lambda_2) [1 - (\lambda_1 + \lambda_0)] - \frac{1}{2} \lambda_2 [1- (\lambda_3 + \lambda_0)^k] \\
	& + [\frac{1}{2} \lambda_0 \lambda_2 - \lambda_1 (\lambda_3 + \frac{1}{2}\lambda_2) ] [1- (\lambda_3 + \lambda_0)^{k-1}] \\
	= & \lambda_0 \Delta^{13}_k(1,1)  +  [\frac{1}{2}\lambda_2 (\lambda_3 + \lambda_0) - \frac{1}{2} \lambda_0 \lambda_2 + \lambda_1 (\lambda_3 + \frac{1}{2}\lambda_2) ] (\lambda_3 + \lambda_0)^{k-1} \\
	= & \lambda_0 \Delta^{13}_k(1,1)  +  [\frac{1}{2}\lambda_2 (\lambda_1 + \lambda_3) + \lambda_1 \lambda_3 ] (\lambda_3 + \lambda_0)^{k-1} \geq 0
	\end{align*}
	where the second equality follow from \eqref{eq:Vnx0},
	\eqref{eq:Vn0y}  and the induction hypothesis. Observing the
	symmetry, we can show $\Delta^{12}_{k+1}(1,1) \geq 0$.
	
	We now study case (2). We can evaluate the term
	$\Delta^{13}_{k+1}(x,1)$ as below.
	\begin{align*}
	\Delta^{13}_{k+1}(x,1)=  & \lambda_1  + \lambda_1 V_{k}(x-1, 1) -\lambda_1 V_{k}(x,1) + \frac{1}{2}\lambda_2 V_{k}(x-1, 1) - \frac{1}{2} \lambda_2 V_{k}(x,0)   \\
	=  & \lambda_1  + \lambda_1                            [1 - \lambda_0 + (\lambda_1 + \frac{1}{2}\lambda_2) V_{k-1}(x-2,1) + (\frac{1}{2}\lambda_2 + \lambda_3) V_{k-1}(x-1, 0) + \lambda_0 V_{k-1}(x-1,1)]  \\
	& - \lambda_1                            [1 - \lambda_0 + (\lambda_1 + \frac{1}{2}\lambda_2) V_{k-1}(x-1,1) + (\frac{1}{2}\lambda_2 + \lambda_3) V_{k-1}(x, 0) + \lambda_0 V_{k-1}(x,1)]    \\
	& + \frac{1}{2}\lambda_2                 [1 - \lambda_0 + (\lambda_1 + \frac{1}{2}\lambda_2) V_{k-1}(x-2,1) + (\frac{1}{2}\lambda_2 + \lambda_3) V_{k-1}(x-1, 0) + \lambda_0 V_{k-1}(x-1,1)]   \\
	& - \frac{1}{2}\lambda_2                 [1-\lambda_3 - \lambda_0 + (\lambda_1 + \lambda_2) V_{k-1}(x-1,0) + \lambda_3 V_{k-1}(x, 0) + \lambda_0 V_{k-1}(x, 0)],
	\end{align*}
	where the second equality follows from the induction hypothesis.
	Note that
	$\lambda_1=\lambda_1(\lambda_0+\lambda_1+\lambda_2+\lambda_3)$. We
	can continue the equality chain above as follows.
	\begin{align*}
	& \Delta^{13}_{k+1}(x,1) \\
	= &   (\lambda_1 + \frac{1}{2}\lambda_2)   [\lambda_1 + \lambda_1 V_{k-1}(x-2,1) -  \lambda_1 V_{k-1}(x-1,1) +\frac{1}{2} \lambda_2 V_{k-1}(x-2,1) -  \frac{1}{2} \lambda_2 V_{k-1}(x-1,0)]\\
	&  +   \lambda_0[\lambda_1 +\lambda_1 V_{k-1}(x-1,1) -  \lambda_1 V_{k-1}(x,1) + \frac{1}{2} \lambda_2 V_{k-1}(x-1,1) - \frac{1}{2} \lambda_2 V_{k-1}(x,0)] \\
	& + (\frac{1}{2}\lambda_2 + \lambda_3)   [\lambda_1 +  \lambda_1 V_{k-1}(x-1,0) -  \lambda_1 V_{k-1}(x,0) + \frac{1}{2} \lambda_2 V_{k-1}(x-1,0)]  \\
	& +\frac{1}{2}\lambda_2 \lambda_3 - \frac{1}{2}\lambda_2 \frac{1}{2}\lambda_2 V_{k-1}(x-1,0) - \frac{1}{2}\lambda_2 \lambda_3V_{k-1}(x,0) \\
	= & (\lambda_1 + \frac{1}{2}\lambda_2) \Delta_k^{13}(x-1,1) + \lambda_0 \Delta_k^{13}(x,1) + (\frac{1}{2}\lambda_2 \lambda_1 + \lambda_3 \lambda_1 + \frac{1}{2}\lambda_2 \lambda_3 )[1+V_{k-1}(x-1,0) - V_{k-1}(x,0)] \geq 0,
	\end{align*}
	where the last inequality follows from the induction hypothesis
	\eqref{eq:Delta13} and Lemma \ref{lemma.Vn.monotone}. Following a
	similar proof, we can show that $\Delta^{12}_{k+1}(x,1),
	\Delta^{12}_{k+1}(1,y), \Delta^{13}_{k+1}(1,y) \geq 0$ for $x,y \geq
	2$.

	Finally, we consider case (4) and evaluate the term
	$\Delta^{13}_{k+1}(x,y)$ below.
	\begin{align*}
	& \Delta^{13}_{k+1}(x,y) \\
	=  & \lambda_1  + \lambda_1 V_{k}(x-1, y) -\lambda_1 V_{k}(x,y) + \frac{1}{2}\lambda_2 V_{k}(x-1, y) - \frac{1}{2} \lambda_2 V_{k}(x,y-1)   \\
	=  & \lambda_1  + \lambda_1                            [1-\lambda_0 + (\lambda_1 + \frac{1}{2}\lambda_2) V_{k-1}(x-2,y) + (\frac{1}{2}\lambda_2 + \lambda_3) V_{k-1}(x-1, y-1) + \lambda_0 V_{k-1}(x-1,y)]  \\
	&            - \lambda_1                            [1-\lambda_0 + (\lambda_1 + \frac{1}{2}\lambda_2) V_{k-1}(x-1,y) + (\frac{1}{2}\lambda_2 + \lambda_3) V_{k-1}(x, y-1) + \lambda_0 V_{k-1}(x,y)]    \\
	&            + \frac{1}{2}\lambda_2                 [1-\lambda_0 + (\lambda_1 + \frac{1}{2}\lambda_2) V_{k-1}(x-2,y) + (\frac{1}{2}\lambda_2 + \lambda_3) V_{k-1}(x-1, y-1) + \lambda_0 V_{k-1}(x-1,y)]   \\
	&            - \frac{1}{2}\lambda2                 [1-\lambda_0 + (\lambda_1 + \frac{1}{2}\lambda_2) V_{k-1}(x-1,y-1) + (\frac{1}{2}\lambda_2 + \lambda_3) V_{k-1}(x, y-2) + \lambda_0 V_{k-1}(x,y-1)],
	\end{align*}
	where the second equality follows from the induction hypothesis.
	Recall that $\sum_{i=0}^3 \lambda_i=1$ and thus
	$\lambda_1=\lambda_1(\lambda_0+\lambda_1+\lambda_2+\lambda_3)$. We
	can continue the equality chain above as follows.
	\begin{align*}
	\Delta^{13}_{k+1}(x,y) =             &   (\lambda_1 + \frac{1}{2}\lambda_2)   [\lambda_1 +                  \lambda_1 V_{k-1}(x-2,y) -                            \lambda_1 V_{k-1}(x-1,y) \\
	&+                                                        \frac{1}{2} \lambda_2 V_{k-1}(x-2,y) -                \frac{1}{2} \lambda_2 V_{k-1}(x-1,y-1)]\\
	&+ (\frac{1}{2}\lambda_2 + \lambda_3)   [\lambda_1 +                   \lambda_1 V_{k-1}(x-1,y-1) -                            \lambda_1 V_{k-1}(x,y-1) \\
	&+                                                        \frac{1}{2} \lambda_2 V_{k-1}(x-1,y-1) -                \frac{1}{2} \lambda_2 V_{k-1}(x,y-2)]\\
	&+ \lambda_0                            [\lambda_1 -      \frac{1}{2} \lambda_2 V_{k-1}(x,y-1) + (\lambda_1 + \frac{1}{2}\lambda_2) V_{k-1}(x-1,y) - \lambda_1 V_{k-1}(x,y)] \\
	=             & (\lambda_1 + \frac{1}{2}\lambda_2) \Delta_k^{13}(x-1,y) + (\frac{1}{2}\lambda_2 + \lambda_3) \Delta_k^{12}(x,y-1) + \lambda_0 \Delta_k^{13}(x,y)                                     \geq 0,
	\end{align*}
	where the last inequality follows from the induction hypothesis.
	Using similar arguments, we can show that $\Delta^{12}_{k+1}(x,y)
	\geq 0$ for $x,y \geq 2$.  Combining the four cases above, we prove
	the desired result. 
	\endproof
	
	
	\subsection{Proof of Proposition~\ref{lemma.M.opt.policy}}
	
	Before we prove Proposition \ref{lemma.M.opt.policy}, we first
	present an auxiliary result.
	\begin{lemma}\label{lemma.M.boundary}
		Consider the ``M'' network and let $n \in \mathds{N}$. Then
		\begin{align}
		V_n(0,m_2,m_3-1) &\ge V_n(0,m_2-1,m_3), \quad m_2 \ge 1,\ m_3 \ge 1,\label{eqn:lemma_boundary1}\\
		V_n(m_1-1,m_2,0) &\ge V_n(m_1,m_2-1,0), \quad m_1 \ge 1,\ m_2 \ge
		1.\label{eqn:lemma_boundary2}
		\end{align}
	\end{lemma}
	
	\proof{Proof.}
	We will
	prove~\eqref{eqn:lemma_boundary1} by induction; this immediately
	implies~\eqref{eqn:lemma_boundary2} due to symmetry.
	
	First, we can see by inspection that
	\begin{equation*}
	V_1(0,m_2,m_3-1) = 1-\lambda_0 \ge V_1(0,m_2-1,m_3).
	\end{equation*}
	
	Now, let $t\in \mathds{N}$ and assume
	that~\eqref{eqn:lemma_boundary1} holds for all $n \le t$. In order
	to show that~\eqref{eqn:lemma_boundary1} holds for $n = t+1$ as
	well, first observe that for $m_1 = 0$, the $M$ model reduces to the
	$N$ model, and by Proposition~\ref{lemma.W.opt.policy} we know that it is
	optimal to offer all slots:
	\begin{align}
	\nonumber V_n(0,m_2,m_3) =& (1-\lambda_0) + (\lambda_1+\frac{1}{2}\lambda_2)V_{n-1}(0,m_2-1,m_3) + \frac{1}{2}\lambda_2 V_{n-1}(0,m_2,m_3-1) \\
	& + \lambda_0 V_{n-1}(0,m_2,m_3), \quad m_2 \ge 1,\ m_3 \ge 1,\ n \in \mathds{N}, \label{eqn.proof.lemma.M.1}
	\end{align}
	and
	\begin{equation}
	V_n(0,m_2,0) = (1-\lambda_0) + (1-\lambda_0)V_{n-1}(0,m_2-1,0) +
	\lambda_0 V_{n-1}(0,m_2,0),\quad m_2 \ge 1,\ n\in \mathds{N}.
	\label{eqn.proof.lemma.M.2}
	\end{equation}
	
	We first prove that~\eqref{eqn:lemma_boundary1} holds for $m_2 \ge
	2$ and $m_3 \ge 2$, and treat the boundary cases separately.
	Using~\eqref{eqn.proof.lemma.M.1} we can write
	\begin{align*}
	\nonumber &V_{t+1}(0,m_2,m_3-1)\\
	\nonumber ={}& (1-\lambda_0) + (\lambda_1+\frac{1}{2}\lambda_2)V_{t}(0,m_2-1,m_3-1) + \frac{1}{2}\lambda_2 V_{t}(0,m_2,m_3-2) + \lambda_0 V_{t}(0,m_2,m_3-1) \\
	\nonumber \ge{} & (1-\lambda_0) + (\lambda_1+\frac{1}{2}\lambda_2)V_{t}(0,m_2-2,m_3) + \frac{1}{2}\lambda_2 V_{t}(0,m_2-1,m_3-1) + \lambda_0 V_{t}(0,m_2-1,m_3) \\
	={}& V_{t+1}(0,m_2-1,m_3).
	\end{align*}
	Here we use the induction hypothesis~\eqref{eqn:lemma_boundary1}
	(with $n = t$) for the inequality, and
	use~\eqref{eqn.proof.lemma.M.1} for the second equality.
	
	For the case $m_2 \ge 2$ and $m_3 = 1$ we
	use~\eqref{eqn.proof.lemma.M.2} to obtain
	\begin{align*}
	\nonumber &V_{t+1}(0,m_2,0) \\
	= & (1-\lambda_0) + (1-\lambda_0)V_t(0,m_2-1,0) + \lambda_0 V_t(0,m_2,0)\\
	\geq & (1-\lambda_0) + (\lambda_1 + \frac{1}{2}\lambda_2) V_t(0,m_2-2,1) + \frac{1}{2} \lambda_2 V_t(0,m_2-1,0) + \lambda_0 V_t(0,m_2-1,1) \\
	= & V_{t+1}(0,m-1,1),
	\end{align*}
	where the inequality follows from 
	the induction hypothesis~\eqref{eqn:lemma_boundary1}, and the final
	equality from our knowledge on the optimal control for $n = t+1$,
	see~\eqref{eqn.proof.lemma.M.1}.
	
	For the case $m_2 = 1$ and $m_3 \ge 2$ we write,
	using~\eqref{eqn.proof.lemma.M.1},
	\begin{align*}
	\nonumber &V_{t+1}(0,1,m_3-1)\\
	\nonumber ={}& (1-\lambda_0) + (\lambda_1 + \frac{1}{2}\lambda_2) V_t(0,0,m_3-1) + \frac{1}{2} \lambda_2 V_t(0,1,m_3-2) + \lambda_0V_t(0,1,m_3-1) \\
	\nonumber = {}& (1-\lambda_0 - \lambda_1) + \frac{1}{2} \lambda_2 V_t(0,0,m_3-1) + \frac{1}{2} \lambda_2 V_t(0,1,m_3-2) +\lambda_1 \big( 1+ V_t(0,0,m_3-1) \big) \\
	\nonumber & + \lambda_0V_t(0,1,m_3-1) \\
	\nonumber \ge{}& (1-\lambda_0 - \lambda_1) +  \lambda_2 V_t(0,0,m_3-1) +  (\lambda_0 + \lambda_1) V_t(0,0,m_3) \\
	=& V_{t+1}(0,0,m_3).
	\end{align*}
	For the second inequality, we use the induction
	hypothesis~\eqref{eqn:lemma_boundary1} and apply
	Lemma~\ref{lemma.Vn.monotone}(ii) to show that $1 + V_t(0,0,m_3-1)
	\ge V_t(0,0,m_3)$.
	
	The case $m_2 = m_3 = 1$ we can do directly, by observing that
	\begin{equation*}
	V_{t+1}(0,1,0) = 1 -(1-\lambda_0)^{t+1}  \ge 1 -
	(1-\lambda_0-\lambda_2)^{t+1} = V_{t+1}(0,0,1),
	\end{equation*}
	completing the proof. 
	\endproof
	
	
	
	With Lemma~\ref{lemma.M.boundary}, we can now prove
	Proposition~\ref{lemma.M.opt.policy}.
	
	\proof{Proof of Proposition~\ref{lemma.M.opt.policy}.} From the
	boundary conditions, it is easy to see that $V_0(\bfm)=0$ regardless
	of $\bfm$.  When $m_2 = 0$ the problem degenerates into two separate
	problems with a single customer type and single slot type where the
	straightforward optimal decision is to offer all slots to customers.
	When either $m_1 = 0$ or $m_3 = 0$, the problem reduces to an ``N''
	model and it is optimal to offer all available slots (see
	Proposition \ref{lemma.W.opt.policy}). Thus, what remains to be
	shown is that when none of the slots are depleted, it is optimal to
	offer type-1 and type-3 slots, but block type-2 slots.
	
	Throughout this proof we assume that $\bfm \geq (1,1,1)$, unless
	stated otherwise. In this case, the Bellman equation can be written
	as
	\begin{equation}
	\label{eqn.M.bellman} {V_n(\bfm) = \max } \left\{
	\begin{array}{lll}
	1-\lambda_0 + \frac{1}{2}\lambda_1 V_{n-1}(\bfm-\bfe_1) + \frac{1}{2} (\lambda_1+\lambda_2) V_{n-1}(\bfm-\bfe_2) + \frac{1}{2}\lambda_2 V_{n-1}(\bfm-\bfe_3) \\
	+ \lambda_0 V_{n-1}(\bfm),  \\  
	1-\lambda_0 + \frac{1}{2}\lambda_1 V_{n-1}(\bfm-\bfe_1) + (\frac{1}{2}\lambda_1+\lambda_2) V_{n-1}(\bfm-\bfe_2)  + \lambda_0 V_{n-1}(\bfm),  \\ 
	1 - \lambda_0 + \lambda_1 V_{n-1}(\bfm-\bfe_1)+ \lambda_2 V_{n-1}(\bfm-\bfe_3)  + \lambda_0 V_{n-1}(\bfm),  \\ 
	1 - \lambda_0 + (\lambda_1+ \frac{1}{2}\lambda_2) V_{n-1}(\bfm-\bfe_2) + \frac{1}{2}\lambda_2 V_{n-1}(\bfm-\bfe_3)  + \lambda_0 V_{n-1}(\bfm),  \\ 
	\lambda_1 + \lambda_1 V_{n-1}(\bfm - \bfe_1) +  (\lambda_0+\lambda_2) V_{n-1}(\bfm),  \\ 
	1 -\lambda_0 +  (\lambda_1+\lambda_2) V_{n-1}(\bfm-\bfe_2) +  \lambda_0 V_{n-1}(\bfm), \\ 
	\lambda_2 +  \lambda_2 V_{n-1}(\bfm-\bfe_3) + (\lambda_0 + \lambda_1) V_{n-1}(\bfm)  
	\end{array}
	\right\},
	\end{equation}
	where the seven terms in the max operator correspond to the action
	of offering slot types $\{1,2,3\}$, $\{1,2\}$, $\{1,3\}$, $\{2,3\}$,
	$\{1\}$, $\{2\}$ and $\{3\}$, respectively.
	
	For ease of notation we define $\Delta_n^{ij}(\bfm)$ to be the
	difference of the $i$th and $j$th terms in the max operator
	\eqref{eqn.M.bellman} above, $i,j \in \{1,2,\dots,7\}$. To prove the
	desired result it suffices to show for any $n \in \mathds{n}$ that
	$\Delta_n^{3,j} \ge 0$, $j \neq 3$.
	
	First, by writing out the definition,
	\begin{align}
	\Delta^{35}_n(\bfm) &= \lambda_2 [1 + (V_{n-1}(\bfm - \bfe_3) - V_{n-1}(\bfm)) ] \geq 0, \label{eqn.compare35.M}\\
	\Delta^{37}_n(\bfm) &= \lambda_1 [1 + (V_{n-1}(\bfm - \bfe_1) -
	V_{n-1}(\bfm)) ]\geq 0. \label{eqn.compare.37.M}
	\end{align}
	The equalities follow from the fact that $V_{n-1}(\bfm) - V_{n-1}(\bfm -
	\bfe_3) \leq 1$ (for~\eqref{eqn.compare35.M}) and $V_{n-1}(\bfm) -
	V_{n-1}(\bfm - \bfe_1) \leq 1$ (for~\eqref{eqn.compare.37.M}), see
	Lemma \ref{lemma.Vn.monotone}.(i).
	
	The other four inequalities can be written as
	\begin{align}
	\Delta_{n+1}^{31} \ge 0 &\Leftrightarrow \lambda_1 V_n(\bfm-\bfe_1) + \lambda_2 V_n(\bfm-\bfe_3) \ge (\lambda_1+\lambda_2) V_n(\bfm-\bfe_2), \label{eqn:proof_M_1}\\
	\Delta_{n+1}^{32} \ge 0 &\Leftrightarrow \frac{1}{2}\lambda_1 V_n(\bfm-\bfe_1) + \lambda_2 V_n(\bfm-\bfe_3) \ge \big(\frac{1}{2} \lambda_1 + \lambda_2 \big)V_n(\bfm-\bfe_2),\label{eqn:proof_M_2}\\
	\Delta_{n+1}^{34} \ge 0 &\Leftrightarrow \lambda_1 V_n(\bfm-\bfe_1) + \frac{1}{2}\lambda_2 V_n(\bfm-\bfe_3) \ge (\lambda_1+ \frac{1}{2} \lambda_2)V_n(\bfm-\bfe_2),\label{eqn:proof_M_3}\\
	\Delta_{n+1}^{36} \ge 0 &\Leftrightarrow \lambda_1 V_n(\bfm-\bfe_1)
	+ \lambda_2 V_n(\bfm-\bfe_3) \ge (\lambda_1+\lambda_2)
	V_n(\bfm-\bfe_2). \label{eqn:proof_M_4}
	\end{align}
	Note that~\eqref{eqn:proof_M_1} and~\eqref{eqn:proof_M_4} are
	equivalent, as are~\eqref{eqn:proof_M_2} and~\eqref{eqn:proof_M_3},
	due to symmetry. Thus, we limit ourselves to showing
	that~\eqref{eqn:proof_M_1} and~\eqref{eqn:proof_M_2} hold, which we
	will do by induction.
	
	Let $n = 1$, then it is readily seen that for~\eqref{eqn:proof_M_1},
	\begin{equation*}
	\lambda_1 V_1(\bfm-\bfe_1) + \lambda_2 V_1(\bfm-\bfe_3) =
	(\lambda_1+\lambda_2)(1-\lambda_0) =
	(\lambda_1+\lambda_2)V_1(\bfm-\bfe_2),
	\end{equation*}
	and for~\eqref{eqn:proof_M_2},
	\begin{equation*}
	\frac{1}{2}\lambda_1 V_1(\bfm-\bfe_1) +  \lambda_2 V_1(\bfm-\bfe_3)
	= (\frac{1}{2}\lambda_1 + \lambda_2) (1-\lambda_0) =
	(\frac{1}{2}\lambda_1 + \lambda_2) V_1(\bfm-\bfe_2),
	\end{equation*}
	so both hold.
	
	Next we let $t \in \mathds{N}$ and assume
	that~\eqref{eqn:proof_M_1}-\eqref{eqn:proof_M_4} hold for all $n\le
	t - 1$, i.e.,
	\begin{align}
	\lambda_1 V_n(\bfm-\bfe_1) + \lambda_2 V_n(\bfm-\bfe_3) &\ge (\lambda_1+\lambda_2) V_n(\bfm-\bfe_2), \quad n\le t-1,\label{eqn:proof_M_induction_1}\\
	\frac{1}{2}\lambda_1 V_n(\bfm-\bfe_1) +  \lambda_2 V_n(\bfm-\bfe_3)
	&\ge (\frac{1}{2}\lambda_1 + \lambda_2) V_n(\bfm-\bfe_2), \quad n\le
	t-1. \label{eqn:proof_M_induction_2}
	\end{align}
	In this case we know that $g_n$ in~\eqref{eqn.opt.M} provides an
	optimal policy for all $n\le t$. We shall now demonstrate
	that~\eqref{eqn:proof_M_induction_1}
	and~\eqref{eqn:proof_M_induction_2} hold for $n = t$ as well, which
	implies that $g_n$ is also optimal for $n = t+1$. Since we know an
	optimal control policy for $n \le t$, we also know the transition
	probabilities given that we use optimal control.
	\begin{align*}
	p_0(\bfm) &= \lambda_0 + \lambda_1 \indi{m_1 = m_2 = 0} + \lambda_2 \indi{m_2 = m_3 = 0},\\
	p_1(\bfm) &= \lambda_1 \indi{m_1 \ge 1},\\
	p_2(\bfm) &= \lambda_1 \indi{m_1 = 0,\ m_2 \ge 1} + \lambda_2 \indi{m_2 \ge 1,\ m_3 = 0},\\
	p_3(\bfm) &= \lambda_2 \indi{m_3 \ge 1}.
	\end{align*}
	
	
	
	
	Using the above transition probabilities we can compute
	\begin{align}
	V_t(\bfm - \bfe_1) = 1-\lambda_0 + \lambda_1 V_{t-1}(\bfm-2\bfe_1) + \lambda_2 V_{t-1}(\bfm-\bfe_1-\bfe_3) + \lambda_0 V_{t-1}(\bfm-\bfe_1),\label{eqn:proof_M_5}\\
	V_t(\bfm - \bfe_3) = 1-\lambda_0 + \lambda_1 V_{t-1}(\bfm-\bfe_1 -
	\bfe_3) + \lambda_2 V_{t-1}(\bfm-2\bfe_3) + \lambda_0 V_{t-1}(\bfm -
	\bfe_3). \label{eqn:proof_M_6}
	\end{align}
	Moreover, we know from the induction
	hypothesis~\eqref{eqn:proof_M_induction_1} that
	\begin{align}
	\lambda_1 V_{t-1}(\bfm- 2\bfe_1) + \lambda_2 V_{t-1}(\bfm-\bfe_1-\bfe_3) &\ge (\lambda_1+\lambda_2) V_{t-1}(\bfm-\bfe_1 - \bfe_2), \label{eqn:proof_M_7}\\
	\lambda_1 V_{t-1}(\bfm- \bfe_1 - \bfe_3) + \lambda_2
	V_{t-1}(\bfm-2\bfe_3) &\ge (\lambda_1+\lambda_2) V_{t-1}(\bfm-\bfe_2
	- \bfe_3). \label{eqn:proof_M_8}
	\end{align}
	Using~\eqref{eqn:proof_M_5}-\eqref{eqn:proof_M_8}, we can write
	\begin{align*}
	\nonumber & \lambda_1 V_t(\bfm-\bfe_1) + \lambda_2 V_t(\bfm-\bfe_3) \\
	\nonumber&\ge (\lambda_1+\lambda_2)(1-\lambda_0) + \lambda_1 (\lambda_1+\lambda_2) V_{t-1}(\bfm-\bfe_1 - \bfe_2) + \lambda_2 (\lambda_1+\lambda_2) V_{t-1}(\bfm-\bfe_2 - \bfe_3) \\
	\nonumber& + \lambda_0 \big( \lambda_1  V_{t-1}(\bfm - \bfe_1) + \lambda_2  V_{t-1}(\bfm - \bfe_3)\big)\\
	\nonumber&\ge (\lambda_1+\lambda_2)(1-\lambda_0) + \lambda_1 (\lambda_1+\lambda_2) V_{t-1}(\bfm-\bfe_1 - \bfe_2) + \lambda_2 (\lambda_1+\lambda_2) V_{t-1}(\bfm-\bfe_2 - \bfe_3) \\
	\nonumber& + \lambda_0 ( \lambda_1 + \lambda_2)  V_{t-1}(\bfm - \bfe_2)\\
	&= (\lambda_1+\lambda_2) V_t(\bfm-\bfe_2),
	\end{align*}
	where the second inequality follows from the induction hypothesis
	\eqref{eqn:proof_M_induction_1}. This proves the desired inequality.
	
	Similarly, to verify~\eqref{eqn:proof_M_2} we use
	\eqref{eqn:proof_M_5} and \eqref{eqn:proof_M_6} and apply the
	induction hypothesis \eqref{eqn:proof_M_induction_2} to obtain,
	after some rearranging,
	\begin{align}
	\nonumber &\frac{1}{2}\lambda_1 V_t(\bfm-\bfe_1) + \lambda_2V_t(\bfm - \bfe_3)\\
	\nonumber \ge& (\frac{1}{2} \lambda_1 +\lambda_2) (1-\lambda_0) + (\frac{1}{2} \lambda_1 +\lambda_2) \lambda_1 V_{t-1}(\bfm -\bfe_1-\bfe_2) + (\frac{1}{2} \lambda_1 +\lambda_2) \lambda_2 V_{t-1}(\bfm -\bfe_2-\bfe_3) \\
	& + (\frac{1}{2} \lambda_1 +\lambda_2) \lambda_0 (\bfm -\bfe_2) \\
	&= (\frac{1}{2} \lambda_1 +\lambda_2) V_t(\bfm - \bfe_2).
	\label{eqn.proof.M.9}
	\end{align}
	
	
	Next, we verify the induction hypotheses for the various boundary
	cases. First, it is readily verified, using our knowledge of the
	optimal control for $n = t$, that for $m_1 = 1$
	\begin{align}
	\nonumber V_t(\bfm-\bfe_1) & = (1-\lambda_0) + (\lambda_1 + \frac{1}{2}\lambda_2) V_{t-1}(0,m_2-1,m_3) + \frac{1}{2} \lambda_2 V_{t-1}(0,m_2,m_3-1) + \lambda_0 V_{t-1}(0,m_2,m_3) \\
	&\ge 1-\lambda_0 + (\lambda_1 + \lambda_2) V_{t-1}(0,m_2-1,m_3)  +
	\lambda_0 V_{t-1}(0,m_2,m_3), \label{eqn.proof.M.10}
	\end{align}
	where the inequality follows from Lemma~\ref{lemma.M.boundary}.
	Analogously, we derive
	\begin{equation}\label{eqn.proof.M.11}
	V_t(\bfm-\bfe_3) \ge 1 -\lambda_0 + (\lambda_1 + \lambda_2)
	V_{t-1}(m_1,m_2-1,0)  + \lambda_0 V_{t-1}(m_1,m_2,0), \quad m_3 = 1.
	\end{equation}
	
	First we treat the case $m_1 = 1$ and $m_3 \ge 2$.
	Combining~\eqref{eqn:proof_M_6} and~\eqref{eqn.proof.M.10} yields
	\begin{align*}
	\nonumber &\lambda_1 V_t(\bfm-\bfe_1) + \lambda_2 V_t(\bfm-\bfe_3)\\
	\nonumber \ge& \lambda_1 [(1-\lambda_0)+  (\lambda_1 + \lambda_2) V_{t-1}(\bfm - \bfe_1-\bfe_2)  + \lambda_0 V_{t-1}(\bfm-\bfe_1) ] \\
	\nonumber& + \lambda_2 [(1-\lambda_0) + \lambda_1 V_{t-1}(\bfm - \bfe_1 - \bfe_3) + \lambda_2 V_{t-1}(\bfm - 2 \bfe_3) + \lambda_0 V_{t-1}(\bfm - \bfe_3)]\\
	\nonumber\ge& (\lambda_1+\lambda_2)(1-\lambda_0) + (\lambda_1+\lambda_2) \lambda_1 V_{t-1}(\bfm-\bfe_1-\bfe_2) + (\lambda_1+\lambda_2) \lambda_2 V_{t-1}(\bfm-\bfe_2-\bfe_3) \\
	\nonumber&+ (\lambda_1+\lambda_2) \lambda_0 V_{t-1}(\bfm-\bfe_2) \\
	&= (\lambda_1+\lambda_2) V_t(\bfm - \bfe_2),
	\end{align*}
	with the second inequality due to the induction
	hypothesis~\eqref{eqn:proof_M_induction_1}.
	
	In order to show~\eqref{eqn:proof_M_induction_2} we can again
	use~\eqref{eqn:proof_M_6} and~\eqref{eqn.proof.M.10}, and do some
	rearranging to show that
	\begin{align}
	\nonumber &\frac{1}{2}\lambda_1 V_t(\bfm-\bfe_1) + \lambda_2 V_t(\bfm-\bfe_3)\\
	\nonumber \ge& \frac{1}{2} \lambda_1 [(1-\lambda_0)+  (\lambda_1 + \lambda_2) V_{t-1}(\bfm - \bfe_1-\bfe_2)  + \lambda_0 V_{t-1}(\bfm-\bfe_1) ] \\
	\nonumber& + \lambda_2 [(1-\lambda_0) + \lambda_1 V_{t-1}(\bfm - \bfe_1 - \bfe_3) + \lambda_2 V_{t-1}(\bfm - 2 \bfe_3) + \lambda_0 V_{t-1}(\bfm - \bfe_3)]\\
	\nonumber\ge& (\frac{1}{2}\lambda_1+\lambda_2)(1-\lambda_0) + (\lambda_1+\lambda_2) \frac{1}{2}\lambda_1 V_{t-1}(\bfm-\bfe_1-\bfe_2) + \lambda_2 [\frac{1}{2}\lambda_1 V_{t-1}(\bfm-\bfe_1-\bfe_3) \\
	\nonumber &+ (\frac{1}{2}\lambda_1 +\lambda_2) V_{t-1}(\bfm-\bfe_2-\bfe_3)] + (\frac{1}{2}\lambda_1+\lambda_2) \lambda_0 V_{t-1}(\bfm-\bfe_2) \\
	\nonumber  \ge& (\frac{1}{2}\lambda_1+\lambda_2)(1-\lambda_0) + (\lambda_1+\lambda_2) \frac{1}{2}\lambda_1 V_{t-1}(\bfm-\bfe_1-\bfe_2) + \lambda_2 [\frac{1}{2}\lambda_1 V_{t-1}(\bfm-\bfe_1-\bfe_2) \\
	\nonumber &+ (\frac{1}{2}\lambda_1 +\lambda_2) V_{t-1}(\bfm-\bfe_2-\bfe_3)] + (\frac{1}{2}\lambda_1+\lambda_2) \lambda_0 V_{t-1}(\bfm-\bfe_2) \\
	=& (\frac{1}{2}\lambda_1 + \lambda_2) V_{t-1}(\bfm - \bfe_2),
	\label{eqn.proof.M.12}
	\end{align}
	where the second and third equalities follows from the induction
	hypothesis~\eqref{eqn:proof_M_induction_2} and
	Lemma~\ref{lemma.M.boundary}, respectively. This shows that the
	\eqref{eqn:proof_M_induction_2} holds for $m_1 = 1$, $m_3 \ge 2$.
	
	
	The proof for the case $m_1 \ge 2$, $m_3 = 1$ follows from symmetry.
	Finally, we verify the case $m_1 = m_3 = 1$. We first bound,
	using~\eqref{eqn.proof.M.10} and~\eqref{eqn.proof.M.11},
	\begin{align}
	\nonumber & \lambda_1 V_t(\bfm - \bfe_1) + \lambda_2 V_t(\bfm - \bfe_3) \ge (\lambda_1+\lambda_2) (1-\lambda_0) + (\lambda_1+\lambda_2) \lambda_1 V_{t-1}(\bfm - \bfe_1 - \bfe_2)\\
	\nonumber  &  + (\lambda_1+\lambda_2) \lambda_2 V_{t-1}(\bfm - \bfe_2 - \bfe_3) + (\lambda_1+\lambda_2) \lambda_0 V_{t-1}(\bfm - \bfe_2)\\
	\nonumber &= (\lambda_1+\lambda_2) V_t(\bfm - \bfe_2).
	\end{align}
	Using these same inequalities we can show
	\begin{align}
	\nonumber &\frac{1}{2}\lambda_1 V_t(\bfm - \bfe_1) +  \lambda_2 V_t(\bfm - \bfe_3)\\
	\nonumber &\ge \frac{1}{2}\lambda_1 [(1-\lambda_0) + (\lambda_1 + \lambda_2) V_{t-1}(0,m_2-1,1) + \lambda_0 V_{t-1}(0,m_2,1)] \\
	\nonumber &+  \lambda_2 [(1-\lambda_0) + (\lambda_2 + \frac{1}{2}\lambda_1) V_{t-1}(1,m_2-1,0) + \frac{1}{2} \lambda_1 V_{t-1}(0,m_2,0) + \lambda_0 V_{t-1}(1,m_2,0)]\\
	\nonumber & \geq (\frac{1}{2}\lambda_1 + \lambda_2) (1-\lambda_0) + (\frac{1}{2}\lambda_1 + \lambda_2) \lambda_2 V_{t-1}(1,m_2-1,0) + (\frac{1}{2}\lambda_1 + \lambda_2) \lambda_1 V_{t-1}(0,m_2-1,1) \\
	\nonumber &+ \lambda_0 [\frac{1}{2}\lambda_1  V_{t-1}(0,m_2,1) + \lambda_2 V_{t-1}(1,m_2,0)] \\
	\nonumber & \geq (\frac{1}{2}\lambda_1 + \lambda_2) (1-\lambda_0) + (\frac{1}{2}\lambda_1 + \lambda_2) \lambda_2 V_{t-1}(1,m_2-1,0) + (\frac{1}{2}\lambda_1 + \lambda_2) \lambda_1 V_{t-1}(0,m_2-1,1) \\
	\nonumber &+ \lambda_0 (\frac{1}{2}\lambda_1 + \lambda_2)  V_{t-1}(1,m_2-1,1)  \\
	\nonumber & =  \Big(\frac{1}{2}\lambda_1 + \lambda_2 \Big) V_t(\bfm
	- \bfe_2),
	\end{align}
	with the second inequality using Lemma~\ref{lemma.M.boundary}.(i).
	and the third inequality due to the induction
	hypothesis~\eqref{eqn:proof_M_induction_2}. This completes the
	proof. 
	\endproof

	
	\subsection{Proof of Corollary~\ref{mcorollary}}
	
	We prove by contradiction. Suppose \eqref{m.morevaluable} does not
	hold and thus
	\begin{equation}\label{lessval}
	V_n(\mathbf{m}-e_2)>V_n(\mathbf{m}-e_1)~\mbox{and}~V_n(\mathbf{m}-e_2)>V_n(\mathbf{m}-e_3).
	\end{equation}
	In period $n+1$ and at state $\mathbf{m}$, action $\mathbf{d}_1:=(1,0,1)$ yields the value-to-go of
	\begin{equation*}
	p_1(\mathbf{m},\mathbf{d}_1) V_n(\mathbf{m}-e_1)+p_3(\mathbf{m},\mathbf{d}_1) V_n(\mathbf{m}-e_3)+\lambda_0 V_n(\mathbf{m}),
	\end{equation*}
	which is strictly less than the value-to-go under action $\mathbf{d}_2=(0,1,0)$ given by
	\begin{equation*}
	p_2(\mathbf{m},\mathbf{d}_2) V_n(\mathbf{m}-e_2)+\lambda_0 V_n(\mathbf{m}),
	\end{equation*}
	by using \eqref{lessval} and
	$p_1(\mathbf{m},\mathbf{d}_1)+p_3(\mathbf{m},\mathbf{d}_1)=p_2(\mathbf{m},\mathbf{d}_2)=1-\lambda_0$.
	This contradicts the result in Proposition~\ref{lemma.M.opt.policy} on
	the optimality of $\mathbf{d}_1$. \qed
	
	
	%

	\subsection{Proof of Theorem \ref{prop.asym}}
	This proof entails a few key steps. First, we show that the optimal
	amount of the customers scheduled in the fluid model is an upper
	bound to that in the corresponding stochastic model (see Proposition
	\ref{thm.fluid.bound} below). Then, we construct a lower bound for
	the objective value of the stochastic model under any static
	randomized policy (see Lemma \ref{lemma.vn.bound.pip} below).
	Finally, we show that under the static randomized policy $\pi^{p^*}$
	this lower bound, after normalization (i.e., divided by the scaling
	factor $K$), converges to the optimal objective value of the fluid
	model, which is a constant upper bound for the stochastic model. To
	economize our notation in the proof below, we let
	$\Ical=\{1,2,\dots,I\}$ be the set of customer types and
	$\Jcal=\{1,2,\dots,J\}$ be the set of slot types.

	\begin{lemma} \label{thm.fluid.bound}
		\[Z_n(\bfm) \geq V_n(\bfm), ~ \forall n=1,2,\dots, N, ~ \bfm \in \mathbb{Z}_+^J.\]
	\end{lemma}
	
	\proof{Proof.} We first show
	that Problem (P1) has an equivalent dynamic programming (DP)
	formulation. This DP formulation will facilitate our proof that the
	fluid model provides an upper bound for the stochastic model. To
	differentiate from the stochastic model, we let $\tildeV_n(\bfm)$ be
	the maximum amount of fluid that can be served given $n$ periods to
	go and the capacity vector $\bfm$. Consider the following DP
	formulation.
	\begin{align}
	& \tildeV_{n}(\mathbf{m}) =   \max \{\sum_{j \in \Jcal}  y_{j}(n) + \tildeV_{n-1}(\mathbf{m} -  \mathbf{y}(n)) \}, &  \\
	\mbox{subject to:}~~~
	& y_{j}(n) =  \sum_{i\in\Ical}y_{i,j}(n), ~ j \in \Jcal,    \\
	& \mathbf{m} -  \mathbf{y}(n) \geq 0, \\
	& \tildeV_0(\mathbf{x}) = 0, ~\forall \mathbf{x} \geq 0, \\
	&\mbox{and \eqref{zzeroone}, \eqref{zsum}, \eqref{opentime}, 
		\eqref{yexp1}, \eqref{yexp2} defined for $n$ only.}
	\end{align}
	
	Recall that $Z_n(\bfm)$ is the optimal objective value to Problem
	(P1) with $M_j(n) = m_j$ and $n$ periods left to go. We claim that
	\begin{equation}
	\label{eqn.fluid.eq.dp} Z_n(\bfm) = \tildeV_n(\bfm),~ \forall
	n=1,2,\dots, N, ~ \bfm \in \mathbb{Z}_+^J.
	\end{equation}
	
	We use induction to prove this claim. It is easy to check the cases
	for $n=1$. Now suppose that $Z_n(\bfm) = \tildeV_n(\bfm)$ holds for
	$n=2,3,\dots, N-1$, and consider that $n=N$. Consider an optimal
	solution $f^*$ under the LP formulation. Following the decision at
	period $N$ specified by $f^*$ in both the LP and DP formulations. We
	see that the amount of fluid served in period $N$ is the same under
	both formulations, and that the capacity left for period $N-1$ is
	also the same for both formulations. Following the induction
	hypothesis, we know that the total amount of fluid served from $N-1$
	periods onward is the same under both formulations. Now, the optimal
	action for period $N$ under the LP formulation is clearly feasible
	for the DP formulation, but not necessarily optimal. Thus we have
	$Z_N(\bfm) \leq \tildeV_N(\bfm)$.
	
	Taking the optimal action in period $N$ under the DP formulation,
	and apply it to both the DP and LP formulations. Following a similar
	argument above, we can show that $Z_N(\bfm) \geq \tildeV_N(\bfm)$.
	It thus follows that $Z_N(\bfm) = \tildeV_N(\bfm)$, as desired.
	
	Now, to prove Lemma~\ref{thm.fluid.bound}, it suffices to show
	that
	\[\tildeV_n(\bfm) \geq V_n(\bfm), ~ \forall n=1,2,\dots, N, ~ \bfm
	\in \mathbb{Z}_+^J.\]
	We use induction to show this result. We first check the case when
	$n=1$. Suppose that action $k$ corresponds to the optimal action
	taken in the stochastic model at $n=1$. In the fluid model, we use
	the same action.  That is, we set $z_k(n)=1$ and set $z_s(n)=0$ for
	$s \neq k$. The feasibility of the optimal action in the stochastic
	model implies that for each type of the slots opened, there is at
	least 1 unit of capacity. Thus, in the fluid model, we can set
	$\tau_{k,j}(n) = \mathbf{w}_j^k$ for all $j \in \Jcal$ as the
	draining speed for each type of slots is bounded by 1 implied by
	constraint \eqref{yexp1}. Then, one can algebraically check that the
	expected number of customers served in the stochastic model is the
	same as the amount of the fluid served in the fluid model.
	
	Now suppose $\tildeV_n(\bfm) \geq V_n(\bfm)$ holds up to
	$n=2,3,\dots, N-1$ and consider $n=N$. Again, we apply the optimal
	action in the stochastic model at period $N$, say $\dbf$, to the
	fluid model at period $N$. Here we let $p_k(\mathbf{m}, \mathbf{d})$
	be the probability that a type $k$ slot is booked at state $\bfm$ if
	action $\mathbf{d}$ is taken. Using \eqref{eqn.dist.ln}, one can
	check that
	\begin{equation}
	\label{stofluidrel1} \sum_{j \in \Jcal} p_j(\mbf,\dbf)  = \sum_{j
		\in \Jcal}  y_{j}(n)
	\end{equation}
	and
	\begin{equation}
	\label{stofluidrel2} \sum_{j =0}^J p_j(\mbf,\dbf) (\mbf - \bfe_j) =
	\mbf - \ybf(N).
	\end{equation}
	The first equation \eqref{stofluidrel1} above suggests that the
	amount of customers served in both models are the same.  The second
	equation \eqref{stofluidrel2} implies that the system state at
	period $N-1$ in the fluid model is a convex combination of the
	possible states that a stochastic model may reach, in which the
	weights are the associated state transition probabilities. To
	simplify notation, we let  $p_j=p_j(\mbf,\dbf), ~j=0,1,\dots,J$. We
	claim that, for the fluid model,
	\begin{equation}
	\label{bound.claim} \tildeV_{N-1}(\mathbf{m} -  \mathbf{y}(n)) \geq
	\sum_{j=0}^J p_j \tildeV_{N-1}(\mathbf{m} -  \bfe_j),
	\end{equation}
	which will be proved at the end. By the induction hypothesis, we
	have that
	\begin{equation}
	\label{bound.induction} \tildeV_{N-1}(\mathbf{m} -  \bfe_j) \geq
	V_{N-1}(\mathbf{m} - \bfe_j), ~ \forall j=0,1,\dots,J.
	\end{equation}
	It follows that
	\begin{align}
	\tildeV_{N}(\mathbf{m})  & \geq   \sum_{j \in \Jcal} y_{j}(N) + \tildeV_{N-1}(\mathbf{m} -  \mathbf{y}(N)) \label{bound.inequ1} \\
	& =      \sum_{j \in \Jcal} p_j  + \tildeV_{N-1}(\mathbf{m} -  \mathbf{y}(N)) \label{bound.inequ2} \\
	& \geq   \sum_{j =0}^J p_j [\mathbbm{1}_{\{j>0\}}  + \tildeV_{N-1}(\mathbf{m} -   \bfe_j)] \label{bound.inequ3} \\
	& \geq   \sum_{j =0}^J p_j [\mathbbm{1}_{\{j>0\}} + V_{N-1}(\mathbf{m} -  \bfe_j) ] \label{bound.inequ4} \\
	& = V_{N}(\mbf) \label{bound.inequ5}.
	\end{align}
	Inequality \eqref{bound.inequ1} holds as the optimal action $\dbf$
	for the stochastic model may not be optimal for the fluid model;
	\eqref{bound.inequ2} holds because of \eqref{stofluidrel1};
	inequalities \eqref{bound.inequ3} and \eqref{bound.inequ4} follow
	from \eqref{bound.claim} and \eqref{bound.induction}, respectively;
	equality \eqref{bound.inequ5} holds by definition.
	
	Finally, we prove our claim \eqref{bound.claim} for $\mathbf{y}(N)$
	that satisfies \eqref{stofluidrel2}. To do this, we turn to the LP
	formulation (P1) for the fluid model. So $\tildeV_{N-1}(\cdot)$ in
	\eqref{bound.claim} is equal to the optimal objective value of the
	corresponding LP formulation by claim \eqref{eqn.fluid.eq.dp}. To
	simplify the notation, we can imagine that this fluid model can be
	written into the following standard form of LP:
	\begin{align*}
	& \tildeV_{N-1}(\mathbf{h}) =   \max ~\mathbf{c} \mathbf{x}, &  \\
	\mbox{subject to:}~~~
	& A \mathbf{x} = \mathbf{h} \\
	& B \mathbf{x} = \mathbf{b} \\
	& \mathbf{x} \geq 0,
	\end{align*}
	in which $\mathbf{x}$ is the vector of decision variables,
	$\mathbf{h}$ is the vector for slots capacity, $\mathbf{b}$ is the
	vector for other right-hand-side coefficients, and $A$, $B$ and $C$
	are properly constructed matrices representing the coefficients for
	$\mathbf{x}$ in the constraint sets. Denote the optimal decision to
	this LP formulation when $\mathbf{h} = \bfm - \bfe_j$ as
	$\mathbf{x}_j$, $j=0,1,\dots,J$. It is easy to check that a solution
	$\sum_{j =0}^J p_j \mathbf{x}_j$ is feasible (but not necessarily
	optimal) to the LP when $\mathbf{h}$ is replaced by $\bfm -
	\mathbf{y}(N)$ and other coefficients are fixed, due to
	\eqref{stofluidrel2} and that $\sum_{j =0}^J p_j \mathbf{x}_j$ is a
	convex combination of $\mathbf{x}_j$'s. Thus we have
	\[
	\tildeV_{N-1}(\mathbf{m} - \mathbf{y}(N) ) \geq \mathbf{c} \sum_{j
		=0}^J p_j \mathbf{x}_j =  \sum_{j =0}^J p_j \mathbf{c} \mathbf{x}_j
	= \sum_{j =0}^J p_j \tildeV_{N-1}(\mathbf{m} - \bfe_j),
	\]
	proving the claim \eqref{bound.claim} and completing the whole
	proof. \endproof

	Before presenting Lemma \ref{lemma.vn.bound.pip}, we introduce a few
	ancillary notations first. Recall that a static randomized policy
	$\pi^p$ offers $\textbf{w}^k$ with probability $p_k$. Define
	$\Ical_j = \{i: \Omega_{ij}=1\}$ be the set of customer types who
	accept type $j$ slots. Recall that $\Kcal_j = \{s: w^s_j=1, s \in
	\Kcal \}$ be the index set of actions that offer type $j$ slots. Let
	$$\Upsilon_j = \sum_{k \in \Kcal_j} p_k (\sum_{i \in \Ical_j}
	\frac{\lambda_i}{\sum_{l, l \in \Jcal} \Omega_{il} w^k_l})$$ be the
	probability that a type $j$ slot will be taken under policy $\pi^p$
	when $\bfm > 0$. To simplify notations below, let $\Ical_j(\bfm) =
	\Ical_j$ if $m_j > 0$ and $\Ical_j(\bfm)=\emptyset$ if $m_j=0$.
	Define
	\begin{equation} \label{ineq.Upsilon}
	\Upsilon_j(\bfm) = \sum_{k \in \Kcal_j} p_k (\sum_{i \in
		\Ical_j(\bfm)} \frac{\lambda_i}{\sum_{l: m_l
			>0, l \in \Jcal} \Omega_{il} w^k_l})
	\end{equation}
	be the probability that a type $j$ slot will be taken under policy
	$\pi^p$ when some of $m_j$'s are zeros. Note that $\Upsilon_j$ is a
	constant while $\Upsilon_j(\bfm)$ depends on $\bfm$. Also,
	$\Upsilon_j(\bfm) \geq \Upsilon_j$ for $j$ such that $m_j > 0$.

	\begin{lemma} \label{lemma.vn.bound.pip}
		For any static randomized policy $\pi^p$, $V^{\pi^p}_n (\bfm) \geq
		\sum_{j \in \Jcal} \mathbb{E} [Bin(n, \Upsilon_j) \wedge m_j ]$,
		$\forall \bfm \geq 0$ and $n=1,2,\dots,N$.
	\end{lemma}
	
	\proof{Proof.}  We prove this
	result by induction. Consider the case when $n=1$. The above
	inequality holds as equality if $\bfm
	> 0$. If there are some $m_j=0$, then
	\begin{eqnarray*}
		V^{\pi^p}_1 (\bfm) &=& \sum_{j: m_j > 0} \mathbb{E} [Bin(1, \Upsilon_j(\bfm ))] \\
		&=& \sum_{j: m_j > 0} \mathbb{E} [Bin(1, \Upsilon_j(\bfm )) \wedge m_j] \\
		&\geq& \sum_{j: m_j > 0} \mathbb{E} [Bin(1, \Upsilon_j ) \wedge m_j ] \\
		&=& \sum_{j} \mathbb{E} [Bin(n, \Upsilon_j) \wedge m_j ],
	\end{eqnarray*}
	where the first inequality above follows from \eqref{ineq.Upsilon}.
	
	Now, assume that the desired inequality holds up to $n-1$ and
	consider the case of $n$. If $\bfm > 0$, then
	%
	%
	\begin{eqnarray*}
		V^{\pi^p}_n (\bfm)  &=& \sum_{j}  \Upsilon_j (1 + V^{\pi^p}_{n-1}
		(\bfm - e_j)) + (1-\sum_j \Upsilon_j) V^{\pi^p}_{n-1} (\bfm) \\
		&\geq& \sum_j \Upsilon_j [1+ \sum_{s \neq j} E(Bin(n-1, \Upsilon_s)
		\wedge m_s) + E(Bin(n-1, \Upsilon_j) \wedge (m_j-1))] \\
		&& + (1-\sum_j \Upsilon_j) \sum_t E(Bin(n-1,\Upsilon_t) \wedge
		m_t)\\
		&=& \sum_j \Upsilon_j [1 + E(Bin(n-1, \Upsilon_j) \wedge (m_j-1))] +
		\sum_t E(Bin(n-1,\Upsilon_t) \wedge m_t)\\
		&+& \sum_j \Upsilon_j [\sum_{s \neq j} E(Bin(n-1, \Upsilon_s) \wedge
		m_s)  - \sum_t
		E(Bin(n-1,\Upsilon_t) \wedge m_t)] \\
		&=& \sum_j \Upsilon_j [1 + E(Bin(n-1, \Upsilon_j) \wedge (m_j-1))] \\
		&+& \sum_j E(Bin(n-1,\Upsilon_j) \wedge m_j) - \sum_j \Upsilon_j
		E(Bin(n-1, \Upsilon_j) \wedge m_j) \\
		&=& \sum_j \Upsilon_j [1 + E(Bin(n-1, \Upsilon_j) \wedge (m_j-1))] +
		\sum_j (1- \Upsilon_j) E(Bin(n-1,\Upsilon_j) \wedge m_j) \\
		&=& \sum_j E(Bin(n, \Upsilon_j) \wedge m_j).
	\end{eqnarray*}
	
	If there are some $m_j = 0$, then
	\begin{eqnarray*}
		V^{\pi^p}_n (\bfm)  &=& \sum_{j: m_j > 0}  \Upsilon_j(\bfm) (1 +
		V^{\pi^p}_{n-1}
		(\bfm - e_j)) + (1-\sum_j \Upsilon_j(\bfm)) V^{\pi^p}_{n-1} (\bfm) \\
		&\geq& \sum_{j: m_j>0} \Upsilon_j(\bfm) [1+ \sum_{s \neq j}
		E(Bin(n-1, \Upsilon_s)
		\wedge m_s) + E(Bin(n-1, \Upsilon_j) \wedge (m_j-1))] \\
		&& + (1-\sum_{j: m_j>0} \Upsilon_j(\bfm)) \sum_t
		E(Bin(n-1,\Upsilon_t) \wedge
		m_t)\\
		&=& \sum_{j: m_j>0} \Upsilon_j(\bfm) [1 + E(Bin(n-1, \Upsilon_j)
		\wedge (m_j-1))] +
		\sum_t E(Bin(n-1,\Upsilon_t) \wedge m_t)\\
		&+& \sum_{j: m_j>0} \Upsilon_j(\bfm) [\sum_{s \neq j} E(Bin(n-1,
		\Upsilon_s) \wedge m_s)  - \sum_t
		E(Bin(n-1,\Upsilon_t) \wedge m_t)] \\
		&=& \sum_{j: m_j>0} \Upsilon_j(\bfm) [1 + E(Bin(n-1, \Upsilon_j)
		\wedge (m_j-1))] +
		\sum_{t: m_t>0} E(Bin(n-1,\Upsilon_t) \wedge m_t)\\
		&+& \sum_{j: m_j>0} \Upsilon_j(\bfm) [\sum_{s \neq j} E(Bin(n-1,
		\Upsilon_s) \wedge m_s)  - \sum_{t:m_t>0}
		E(Bin(n-1,\Upsilon_t) \wedge m_t)] \\
		&=& \sum_{j:m_j>0} \Upsilon_j(\bfm) [1 + E(Bin(n-1, \Upsilon_j) \wedge (m_j-1))] \\
		&+& \sum_{j:m_j>0} E(Bin(n-1,\Upsilon_j) \wedge m_j) -
		\sum_{j:m_j>0} \Upsilon_j(\bfm)
		E(Bin(n-1, \Upsilon_j) \wedge m_j) \\
		&=& \sum_{j:m_j>0} \Upsilon_j(\bfm) [1 + E(Bin(n-1, \Upsilon_j)
		\wedge (m_j-1))] +
		\sum_{j:m_j>0} (1- \Upsilon_j(\bfm)) E(Bin(n-1,\Upsilon_j) \wedge m_j) \\
		&=& \sum_{j:m_j>0} E\{[Bin(1,\Upsilon_j(\bfm))+Bin(n-1,\Upsilon_j)]
		\wedge m_j\} \\
		&\geq& \sum_{j:m_j>0} E[Bin(Bin(n,\Upsilon_j) \wedge m_j] \\
		&=&\sum_j E(Bin(n, \Upsilon_j) \wedge m_j),
	\end{eqnarray*}
	where the last inequality results from \eqref{ineq.Upsilon}. This
	completes the proof. \endproof
	
	Before presenting the proof of Theorem \ref{prop.asym}, we need two
	more ancillary results. The first result states that for the $K$th
	problem, its objective value of the fluid model is $K$ times that of
	the base model with $K=1$. The second is a convergence result, and
	we let $\stackrel{D}{\longrightarrow}$ denote convergence in
	distribution.

	\begin{lemma}
		\label{lemma.fluid.scaling} $Z_{NK}(\bfm K) = K Z_N(\bfm),~\forall
		\bfm \geq 0, ~ K=1,2,3,\dots,$.
	\end{lemma}
	
	\proof{Proof.} It suffices to
	show that
	\begin{equation}
	\label{eq:fluid.ineq1} K^{-1} Z_{NK}(\bfm K) \leq Z_N(\bfm),~\forall
	\bfm \geq 0, ~ K=1,2,3,\dots,
	\end{equation}
	and
	\begin{equation} \label{eq:fluid.ineq2}
	Z_{NK}(\bfm K) \geq K Z_N(\bfm),~\forall \bfm \geq 0, ~
	K=1,2,3,\dots.
	\end{equation}
	
	To show \eqref{eq:fluid.ineq1}, we let $z^*_k(i,K), ~\forall
	i=1,2,\dots,NK, \forall k=1,2,\dots,2^J$ be the optimal solution for
	the $K$th fluid model. Define
	\[z_k(n,1) = \frac{\sum_{i=(n-1)K+1}^{nK} z^*_k(i,K)}{K},~\forall n=1,2,\dots,N,~k=1,2,\dots,2^J. \]
	It suffices to show that $z_k(n,1)$ is a feasible solution for the
	base fluid model with $K=1$, and gives an objective value of $K^{-1}
	Z_{NK}(\bfm K)$. It is easy to check that $z_k(n,1)$ satisfies
	\eqref{zzeroone}, because that $0 \leq z^*_k(i,K) \leq 1$ by
	definition. To check that $z_k(n,1)$ satisfies \eqref{zsum}, we have
	that
	\begin{align*}
	\sum_k z_k(n,1) &= \frac{1}{K} \sum_k \sum_{i=(n-1)K+1}^{nK}
	z^*_k(i,K) \\
	&= \frac{1}{K}  \sum_{i=(n-1)K+1}^{nK} \underbrace{\sum_k z^*_k(i,K)}_\text{$=1$ by definition of $z^*_k(i,K)$} \\
	&= \frac{1}{K} \cdot K = 1.
	\end{align*}
	Constraints \eqref{opentime}-\eqref{yexp2} hold as they are simply
	definitions of $\tau_{k,j}(n)$ and $y_{i,j}(n)$.
	
	Now, let $M_j(n,K)$ be the capacity left for slot type $j$ with $n$
	periods to go in the $K$th fluid model under its respective solution
	under consideration. To show that $z_k(n,1)$ gives an objective
	value of $K^{-1} Z_{NK}(\bfm K)$, it suffices to show
	\begin{equation} \label{eq:Mjn1}
	M_j(n,1) = \frac{1}{K} M_j(nK,K),~\forall n=1,2,\dots,N.
	\end{equation}
	We prove \eqref{eq:Mjn1} by induction. First consider $n=N$. By
	definition, we have
	$$M_j(N,1) = m_j = \frac{1}{K}(m_j K) = \frac{1}{K} M_j(NK,K).$$
	Assume that \eqref{eq:Mjn1} holds for $N-1$, $N-2$, ..., $n$.
	Consider the case of $n-1$.
	\begin{align*}
	M_j(n-1,1)  &= M_j(n,1) - \sum_i \sum_{k \in \Kcal_j} z_k(n,1)
	\mathbf{w}^k_j
	\frac{\lambda_i}{\sum_{l\in\Jcal}\min\{\Omega_{i,l},\mathbf{w}^k_l\}}
	\\         &= M_j(n,1) - \sum_i \sum_{k \in \Kcal_j}
	\frac{\sum_{s=(n-1)K+1}^{nK} z^*_k(s,K)}{K} \mathbf{w}^k_j
	\frac{\lambda_i}{\sum_{l\in\Jcal}\min\{\Omega_{i,l},\mathbf{w}^k_l\}}
	\\
	&= \frac{1}{K} (\underbrace{K M_j(n,1)}_\text{$=M_j(nK, K)$ by induction}
	- \underbrace{\sum_{s=(n-1)K+1}^{nK}  \sum_i \sum_{k \in \Kcal_j} z^*_k(s,K)
		\mathbf{w}^k_j \frac{\lambda_i}{\sum_{l\in\Jcal}\min\{\Omega_{i,l},\mathbf{w}^k_l\}}}_\text{fluid taking type $j$ slots from periods $nK$ to $(n-1)K+1$ in
		the $K$th fluid
		model}) \\
	&= \frac{1}{K} M_j(nK-K, K) \\
	&= \frac{1}{K} M_j((n-1)K, K),
	\end{align*}
	which proves \eqref{eq:Mjn1}. And thus \eqref{eq:fluid.ineq1} holds.
	
	Next, we prove \eqref{eq:fluid.ineq2}. Let $z^*_k(i,1)$ be the
	optimal solution to the base fluid model with $K=1$. For
	$i=1,2,\dots,N$, define
	\[
	z_k(n,K) = z^*_k(i,1),~\textrm{if $(i-1)K+1 \leq n \leq iK$.}
	\]
	It suffices to show that $z_k(n,K)$ is a feasible solution to the
	$K$th fluid model and gives rise to an objective value of $K
	Z_N(\bfm)$. It is easy to check that $z_k(n,K)$ satisfies
	\eqref{zzeroone}, \eqref{opentime}, \eqref{yexp1} and \eqref{yexp2}.
	To check \eqref{zsum}, note that
	\[
	\sum_k z_k(n,K) = \sum_k z_k^*(i,1) = 1, ~\textrm{for
		$i=1,2,\dots,N$ and $(i-1)K+1 \leq n \leq iK$.}
	\]
	
	To show that $z_k(n,K)$ gives rise to an objective value of $K
	Z_N(\bfm)$, it suffices to show that
	\begin{equation} \label{eq:Mjn2}
	M_j(nK,K) = K M_j(n,1),~\forall n=1,2,\dots,N.
	\end{equation}
	We prove \eqref{eq:Mjn2} by induction. For $n=N$, \eqref{eq:Mjn2}
	holds by definition. Assume that \eqref{eq:Mjn2} holds for $N-1$,
	$N-2$, ..., $n$. Consider the case of $n-1$.
	\begin{align*}
	M_j((n-1)K,K)  &= M_j(nK,K) - \sum_{s=(n-1)K+1}^{nK} \sum_i \sum_{k
		\in \Kcal_j} z_k(s,K) \mathbf{w}^k_j
	\frac{\lambda_i}{\sum_{l\in\Jcal}\min\{\Omega_{i,l},\mathbf{w}^k_l\}}
	\\
	&= \underbrace{M_j(nK,K)}_\text{$=K M_j(n, 1)$ by induction} - \sum_{s=(n-1)K+1}^{nK} \sum_i \sum_{k
		\in \Kcal_j} z^*_k(n,1) \mathbf{w}^k_j
	\frac{\lambda_i}{\sum_{l\in\Jcal}\min\{\Omega_{i,l},\mathbf{w}^k_l\}}\\
	&= K M_j(n, 1) - K \underbrace{\sum_i \sum_{k
			\in \Kcal_j} z^*_k(n,1) \mathbf{w}^k_j
		\frac{\lambda_i}{\sum_{l\in\Jcal}\min\{\Omega_{i,l},\mathbf{w}^k_l\}}}_\text{fluid taking type $j$ slots in period $n$ for model with $K=1$}\\
	&= K M_j(n-1,1),
	\end{align*}
	which proves \eqref{eq:Mjn2}. Thus \eqref{eq:fluid.ineq2} holds. This completes the proof. 
	
	

	\endproof
	
	\begin{lemma}\emph{\cite[p.\ 34]{billingsley1968convergence}}
		\label{lemma.fun.converge} Suppose that $X$ and $\{X_k\}$ are
		$\mathbb{R}^n$-valued random variables such that $X_k
		\stackrel{D}{\longrightarrow} X$, and suppose that the functions
		$h_k: \mathbb{R}^n \rightarrow \mathbb{R}$ converge uniformly on
		compact sets to a continuous function $h:\mathbb{R}^n \rightarrow
		\mathbb{R}$. Then $h_k(X_k) \stackrel{D}{\longrightarrow} h(X)$.
	\end{lemma}

	We are now in a position to prove Theorem \ref{prop.asym}, the main
	result in this section.
	
	\proof{Proof of Theorem~\ref{prop.asym}.} Consider the stochastic
	scheduling policy $\pi^{p^*}$ defined above. To simplify notations
	below, we define
	\begin{equation} \label{eqn.Upsilon.star}
	\Upsilon_j^* = \sum_{k \in \Kcal_j} p^*_k (\sum_{i \in \Ical_j}
	\frac{\lambda_i}{\sum_{l, l \in \Jcal} \Omega_{il} \textbf{w}^k_l})
	\end{equation} be the probability that a type $j$ slot will be taken under
	policy $\pi^{p^*}$ when $\bfm > 0$. We have
	\begin{equation}
	\label{eq:asym.chain.bound} K^{-1} \sum_{j} \mathbb{E}[Bin(NK,
	\Upsilon^*_j) \wedge m_j K ] \leq K^{-1} V^{\pi^{p^*}}_{NK} (\bfm K)
	\leq K^{-1} Z_{NK}(\bfm K) = Z_N(\bfm),
	\end{equation} where the first inequality follows from Lemma
	\ref{lemma.vn.bound.pip}, the second inequality follows from
	Lemma~\ref{thm.fluid.bound}, and the last equality follows
	from Lemma \ref{lemma.fluid.scaling}. The LHS of
	\eqref{eq:asym.chain.bound} can be rewritten as
	\begin{equation}
	\label{eq:asym.chain.bound2}
	K^{-1} \sum_{j} \mathbb{E}[Bin(NK, \Upsilon^*_j) \wedge m_j K ]
	=  \sum_{j} \mathbb{E}[K^{-1} Bin(NK, \Upsilon^*_j) \wedge m_j ]
	\end{equation}
	
	The strong law of large numbers implies that
	$$K^{-1} Bin(NK, \Upsilon^*_j) \stackrel{D}{\longrightarrow}
	N\Upsilon^*_j~~\textrm{as $K \rightarrow \infty$}.$$ Applying Lemma
	\ref{lemma.fun.converge}, we conclude that
	$$\sum_{j} [K^{-1} Bin(NK, \Upsilon^*_j) \wedge m_j ] \stackrel{D}{\longrightarrow} \sum_{j} [N\Upsilon^*_j \wedge m_j
	] ~~\textrm{as $K \rightarrow \infty$}.$$ Because that the random
	variable $\sum_{j} [K^{-1} Bin(NK, \Upsilon^*_j) \wedge m_j ]$ is
	uniformly bounded by $\sum_j m_j$, we have that
	\begin{equation}
	\label{eq:asym.chain.bound3}\lim_{K \rightarrow \infty}\mathds{E}
	\sum_{j} [K^{-1} Bin(NK, \Upsilon^*_j) \wedge m_j ] = \sum_{j} (N
	\Upsilon^*_j \wedge m_j )
	= Z_N(\bfm),
	\end{equation}
	where the last equality follows from the definition of $p^*$ and
	$\Upsilon^*_j$. Specifically, $\Upsilon_j^*$ is defined based on the
	fluid model, via the use of $p_k^*$ which is the proportion of time
	in which slot type $k$ is offered in the fluid model; see equations
	\eqref{eq:p*def} and \eqref{eqn.Upsilon.star}. Note that in the
	fluid model, we have constraints ensuring that a slot type can only
	be offered if it is still available. Thus $\Upsilon_j^*$ matches
	exactly the proportion of time in which slot type $j$ is being
	drawn. This quantity times $N$ is exactly the amount of type $j$
	slots being taken in the fluid model. In fact, $(N \Upsilon^*_j
	\wedge m_j ) = N\Upsilon_j^*$ and $\sum_{j} N \Upsilon^*_j =
	Z_N(\bfm)$. Combining \eqref{eq:asym.chain.bound},
	\eqref{eq:asym.chain.bound2} and \eqref{eq:asym.chain.bound3} gives
	the desired result and completes the whole proof. 
	\endproof


	\subsection{Proof of Theorem~\ref{2approximation}} \proof{Proof.} We
	prove this by induction. Let $\Omega$ be any preference matrix. For
	$n=1$, $\pi_0$ is optimal and thus
	$V_1(\mathbf{m})=V_{1,\pi_0}(\mathbf{m})\leq
	2V_{1,\pi_0}(\mathbf{m})$. Suppose the desired result holds for any
	$n\leq k-1$ and state $\mathbf{m}$.
	
	Now we consider two systems, one under an optimal policy and the
	other under $\pi_0$, both starting from state $\mathbf{m}$ in period
	$n=k$ and operating independently from each other. We denote by
	$L^*_k(\mathbf{m})$ the slot type filled in period $k$ in the first
	system (i.e., using an optimal policy), and by
	$L^{\pi_0}_k(\mathbf{m})$ the slot type filled in period $k$ in the
	second system (i.e., using $\pi_0$). These two random variables are
	independent and we shall next condition on them. Specifically, let
	$V_k(\mathbf{m}|L^*_k(\mathbf{m}),L^{\pi_0}_k(\mathbf{m}))$ denote
	the value attained in the first system conditioning on these two
	random variables. 
	types.
	We have that
	\begin{align}
	V_k(\mathbf{m}|L^*_k(\mathbf{m}),L^{\pi_0}_k(\mathbf{m}))&=E[\mathbbm{1}_{\{L^*_k(\mathbf{m})>0\}}|L^*_k(\mathbf{m}),L^{\pi_0}_k(\mathbf{m})]+V_{k-1}(\mathbf{m}-\mathbf{e}_{L^*_k(\mathbf{m})})\notag\\
	&=\mathbbm{1}_{\{L^*_k(\mathbf{m})>0\}}+V_{k-1}(\mathbf{m}-\mathbf{e}_{L^*_k(\mathbf{m})})\notag\\
	&\leq \mathbbm{1}_{\{L^*_k(\mathbf{m})>0\}}+V_{k-1}(\mathbf{m})\label{valueine1}\\
	&\leq
	\mathbbm{1}_{\{L^*_k(\mathbf{m})>0\}}+\mathbbm{1}_{\{l>0\}}+V_{k-1}(\mathbf{m}-\mathbf{e}_{l}),~\forall
	l\in\bar{S}(\bfm)\cup\{0\},\label{valueine2}
	\end{align}
	where inequality \eqref{valueine1} follows from the left inequality
	of Lemma \ref{lemma.Vn.monotone} (ii) and inequality
	\eqref{valueine2} holds due to the right inequality of Lemma
	\ref{lemma.Vn.monotone} (ii). We now let $l=L^{\pi_0}_k(\mathbf{m})$
	in \eqref{valueine2} and in turn have
	\begin{equation*}
	V_k(\mathbf{m}|L^*_k(\mathbf{m}),L^{\pi_0}_k(\mathbf{m}))\leq
	\mathbbm{1}_{\{L^*_k(\mathbf{m})>0\}}+\mathbbm{1}_{\{L^{\pi_0}_k(\mathbf{m})>0\}}+V_{k-1}(\mathbf{m}-\mathbf{e}_{L^{\pi_0}_k(\mathbf{m})}).
	\end{equation*}
	Further applying the induction hypothesis to the above inequality,
	we obtain
	\begin{equation*}
	V_k(\mathbf{m}|L^*_k(\mathbf{m}),L^{\pi_0}_k(\mathbf{m}))\leq
	\mathbbm{1}_{\{L^*_k(\mathbf{m})>0\}}+\mathbbm{1}_{\{L^{\pi_0}_k(\mathbf{m})>0\}}+2V_{k-1,\pi_0}(\mathbf{m}-\mathbf{e}_{L^{\pi_0}_k(\mathbf{m})}).
	\end{equation*}
	Finally, taking expectations on both sides of the above inequality
	leads to
	\begin{equation*}
	V_k(\mathbf{m})\leq
	E[\mathbbm{1}_{\{L^*_k(\mathbf{m})>0\}}]+E[\mathbbm{1}_{\{L^{\pi_0}_k(\mathbf{m})>0\}}]+2E[V_{k-1,\pi_0}(\mathbf{m}-\mathbf{e}_{L^{\pi_0}_k(\mathbf{m})})].
	\end{equation*}
	Now note that $E[\mathbbm{1}_{\{L^{\pi_0}_k(\mathbf{m})>0\}}]\geq
	E[\mathbbm{1}_{\{L^*_k(\mathbf{m})>0\}}]$ by the definition of the
	greedy policy, and hence we arrive at
	\begin{equation*}
	V_k(\mathbf{m})\leq
	2E[\mathbbm{1}_{\{L^{\pi_0}_k(\mathbf{m})>0\}}]+2E[V_{k-1,\pi_0}(\mathbf{m}-\mathbf{e}_{L^{\pi_0}_k(\mathbf{m})})]=2V_{k,\pi_0}(\mathbf{m}).
	\end{equation*}
	\endproof

	\section{Proof of the Results in Section \ref{sec.seq.model}}
	
	\subsection{Proof of Lemma~\ref{lemma.Vn.monotone.seq}}
	
	\proof{Proof.}
	The proof
	follows that of Lemma~\ref{lemma.Vn.monotone} with some minor
	modifications.  In particular, to prove
	\eqref{eqn.lemma.monotonicity.bd.slot}, we define a decision rule
	${\mathbf{h}}$ in period $t+1$ which acts the same as
	$\mathbf{g}^*_{t+1}(\mathbf{m}+\mathbf{e}_j)$ regarding all slot
	types but type $j$. For type $j$, ${\mathbf{h}}$ does not offer it
	in any subsets it offers. That is, $\mathbf{h} =
	\mathbf{g}^*_{t+1}(\mathbf{m}+\mathbf{e}_j)$ except that we enforce
	${h}_{kj}=0,~,\forall k$.  All other parts of the proof readily follow.
	\endproof
	
	\subsection{Proof of Lemma~\ref{lemma.fullset.seq}}
	
	\proof{Proof.} For convenience, we introduce some notation here.
	We say the decision rule in period $n$ given the system occupies
	state $\mathbf{m} \in S$ can be described by a matrix-valued
	function: $\mathbf{g}_n: S \rightarrow \mathbf{d}$ in which
	$\mathbf{d}=\{d_{kj}\}$ is a $J$ by $J$ matrix,  $d_{kj} \in
	\{0,1\}$.
	If $d_{kj}=1$, type $j$ slots are offered in the $k$th subset.
	Since these subsets offered are mutually exclusive, $\sum_{k=1}^K
	d_{kj} \leq 1, ~\forall j$. As before, depleted slot types cannot
	be offered:  $d_{kj}\leq m_j$.

	Let $\hat{\mathbf{d}}$ denote an optimal decision rule. Without loss of generality, we assume that $m_j > 0, ~\forall j
	\in \Jcal$. Otherwise we would consider a network where the preference matrix has been modified by removing empty slots. Let $\hat{\Jcal}=\{j: ~\sum_{k=1}^{K-1} d_{kj} = 1, ~
	j \in \Jcal\}$ be the set of slot types offered by
	$\hat{\mathbf{d}}$ collectively in all subsets it offers. Assume
	that $\Jcal \setminus \hat{\Jcal} \neq \emptyset$. Consider
	another decision rule $\tilde{\mathbf{d}}$ which follows exactly
	the same sequential offering rule as $\hat{\mathbf{d}}$, and then
	offers all slots types in $\Jcal \setminus \hat{\Jcal}$ as the
	$K$th offer set. So $\tilde{\mathbf{d}}$ eventually offers all
	slot types.  To prove the desired result, it suffices to show that
	$\tilde{\mathbf{d}}$ is no worse than $\hat{\mathbf{d}}$, and thus
	must be optimal as well.
	
	First consider a policy that uses $\hat{\mathbf{d}}$ in the first slot, and then
	follows the optimal scheduling rule. Let ${V}^{\hat{d}}_n(\mbf)$
	denote the expected objective value following such a policy. Recall
	that $V_n(\mbf)$ is the optimal expected objective value, and that
	${p}_j(\mbf,\dbf)$ denotes the probability that a type $j$ slot will
	be booked in state $\mbf$ if decision rule $\dbf$ is used. It
	follows that
	\begin{equation}
	\label{eq:Vn.subset.seq} {V}^{\hat{d}}_n(\mbf) =   \sum_{j \in
		\hat{\Jcal}} {p}_j(\mbf,\hat{\dbf}) + \sum_{j \in \hat{\Jcal}}
	{p}_j(\mbf,\hat{\dbf}) V_{n-1}(\mbf - \mathbf{e}_{j}) + [1-\sum_{j
		\in \hat{\Jcal}}  {p}_j(\mbf,\hat{\dbf})] V_{n-1}(\mbf).
	\end{equation}
	Then, consider a a policy that uses $\tilde{\mathbf{d}}$ first, and
	then follow the optimal scheduling rule. The expected objective
	valuing of this policy is
	\begin{equation}
	\label{eq:Vn.fullset.seq} {V}^{\tilde{d}}_n(\mbf) =   \sum_{j \in
		{\Jcal}} {p}_j(\mbf,\tilde{\dbf}) + \sum_{j \in {\Jcal}}
	{p}_j(\mbf,\tilde{\dbf}) V_{n-1}(\mbf - \mathbf{e}_{j}) + [1-\sum_{j
		\in {\Jcal}}  {p}_j(\mbf,\tilde{\dbf})] V_{n-1}(\mbf)
	\end{equation}
	It is easy to check that ${p}_j(\mbf,\hat{\dbf}) =
	{p}_j(\mbf,\tilde{\dbf})$ for $j \in \hat{\Jcal}$, as $\hat{\dbf}$
	acts the same as $\tilde{\dbf}$ in the first $K-1$ offer sets that
	cover slots types in $\hat{\Jcal}$. Subtracting
	\eqref{eq:Vn.subset.seq} from \eqref{eq:Vn.fullset.seq} and
	simplifying, we arrive at
	\[
	{V}^{\tilde{d}}_n(\mbf) - {V}^{\hat{d}}_n(\mbf) = \sum_{j \in
		\Jcal \setminus \hat{\Jcal}} {p}_j(\mbf,\tilde{\dbf}) (1+
	V_{n-1}(\mbf - \mathbf{e}_{j}) - V_{n-1}(\mbf)) \geq 0,
	\]
	where the last inequality directly follows from Lemma
	\ref{lemma.Vn.monotone.seq}, proving the desired result. 
	\endproof
	
	\subsection{Proof of Theorem~\ref{lemma.one.by.one.seq} }
	\proof{Proof.}
	Lemma
	\ref{lemma.fullset.seq} suggests that there exists an optimal
	decision rule $\bfS^*=S^*_1$-$\dots$-$S^*_K$ such that $\cup_{i=1}^K S^*_i
	= J$. Suppose that $S^*_1$-$\dots$-$S^*_K$ does not take the form as
	desired, we will show below that the objective value obtained by partitioning $\bfS^*$ into singletons
	$\{j_1\}$-$\dots$-$\{j_J\}$ is no worse than that of
	$S^*_1$-$\dots$-$S^*_K$.
	
	If there exists some $k$ that $|S^*_k|>1$, let us consider an
	alternative decision rule
	\begin{equation}
	\label{eqn.alter1} \hat{\bfS}^* = S^*_1-\dots-S^*_{k-1}-\{t_1\}-S^*_k \setminus
	\{t_1\} -\dots-S^*_K,
	\end{equation}
	such that
	\begin{equation}
	\label{eqn.alter1.con} V_{n-1}(\mathbf{m} - \bfe_{t_1}) \geq
	V_{n-1}(\mathbf{m} - \bfe_{t}), ~ \forall t \in S^*_k \setminus
	\{t_1\}.
	\end{equation}
	This new decision rule follows the same offering sequence as the
	original rule, except that it splits the offer set $S^*_k$ into two
	sub-offer sets $S^*_{k-1}$ and $\{t_1\}$.
	
	Now, we will show that $\hat{\bfS}^*$ does no worse than $\bfS^*$.  To do that, let $V^1_n(\mathbf{m})$ be the expected
	number of slots filled at the end of the booking horizon by
	following decision rule $\hat{\bfS}^*$ at period $n$ and then
	following the optimal decision afterwards. Let $\Delta^1 =
	V_n(\mathbf{m}) - V^{\hat{\bfS}^*}_n(\mathbf{m})$. Let $I^* = \{i:
	\Omega_{it_1}=1, \sum_{j \in S^*_k \setminus t_1} \Omega_{ij} \geq
	1, \sum_{j \in \cup_{l=1}^{k-1} S^*_l}
	\Omega_{ij} =0\}$ be the set of customer types that accept type
	$t_1$ slots and also at least one slot type in the set of $S^*_k
	\setminus t_1$, but do not accept any slot type that has been
	offered so far in sets $S^*_1$ through $S^*_{k-1}$. Let $J^*(i) =
	\{j: j \in S^*_k, \Omega_{ij} =1\}$ be the subset of slots type in
	$S^*_k$ that are acceptable by customer type $i$, $i \in I^*$.
	Clearly, $t_1 \in J^*(i)$. One can find that
	\[
	\Delta^1 = \sum_{i \in I^*} \frac{\lambda_i}{|J^*(i)|} \sum_{j \in
		J^*(i)} V_{n-1} (\mathbf{m} - \bfe_{j}) - \sum_{i \in I^*} \lambda_i
	V_{n-1} (\mathbf{m} - \bfe_{t_1}) \leq 0,
	\]
	where the last equality follows from \eqref{eqn.alter1.con}, proving
	that $\hat{\bfS}^*$ does no worse than $\bfS^*$.
	
	Following the procedure above to keep splitting offer sets that
	contain more than one slot types, we can obtain an optimal action of
	form $\{j'_1\}$-$\dots$-$\{j'_J\}$ so that each sequential offer set
	contains exactly one slot type. Suppose that
	$\{j'_1\}$-$\dots$-$\{j'_J\}$ does not follow the order desired.
	That is, there exists  $1 \leq u \leq J+1$ such that
	$V_{n-1}(\mathbf{m} - \bfe_{j'_u}) < V_{n-1}(\mathbf{m} -
	\bfe_{j'_{u+1}})$. Consider another decision rule with only $j'_u$
	and $j'_{u+1}$ switched and others remained the same order.
	\begin{equation}
	\label{eqn.alter2}
	\{j'_1\}-\dots-\{j'_{u+1}\}-\{j'_u\}-\dots-\{j'_J\}.
	\end{equation}
	It suffices to show the claim that \eqref{eqn.alter2} either
	provides the same objective value as $\{j'_1\}$-$\dots$-$\{j'_J\}$,
	or strictly higher, which contradicts with the optimality of
	$\{j'_1\}$-$\dots$-$\{j'_J\}$, and thus for all $1 \leq u \leq J+1$,
	$V_{n-1}(\mathbf{m} - \bfe_{j'_u}) \geq V_{n-1}(\mathbf{m} -
	\bfe_{j'_{u+1}})$ as desired.
	
	To show the claim above, let $I' = \{i: \Omega_{ij'_u}=1,
	\Omega_{ij'_{u+1}} = 1, \sum_{v=1}^{u-1} \Omega_{ij_v} =0\}$ be the
	set of customer types that accept both types $j'_u$ and $j'_{u+1}$
	slots, but do not accept any slot type that has been offered so far.
	Let $V^2_n(\mathbf{m})$ be the expected number of slots filled at
	the end of the booking horizon by following decision rule
	\eqref{eqn.alter2} at period $n$ and then following the optimal
	decision afterwards. We consider
	\[
	\Delta^2 = V_n(\mathbf{m}) - V^2_n(\mathbf{m}) = \sum_{i \in I'}
	\lambda_i V_{n-1} (\mathbf{m} - \bfe_{j'_u}) - \sum_{i \in I^*}
	\lambda_i V_{n-1} (\mathbf{m} - \bfe_{j'_{u+1}}).
	\]
	If $\sum_{i \in I'} \lambda_i = 0$, then $\Delta^2=0$ and thus
	\eqref{eqn.alter2} is optimal. However, if $\sum_{i \in I'}
	\lambda_i > 0$, then $\Delta^2 < 0$ leading to the contradiction
	desired. This proves our claim and completes the proof. 
	\endproof

	\subsection{Proof of Theorem~\ref{theorem.equivalence}}
	\proof{Proof.} For notational convenience, here we consider the case
	when $\lambda_0=0$ and all customer types $i \in \mathcal{I}$ can be
	covered by at least one slot type left in $\mathbf{m}$. Proofs of
	other cases follow a similar procedure.
	
	It is trivial that $V^s_n(\mathbf{m}) = V^f_n(\mathbf{m})$, for
	$n=0,1$ and for all $\mathbf{m} \geq 0$. Assume the desired equality
	holds up to $n=t-1$, and consider $n=t$.
	Let $V^f_n(\mathbf{m} | i)$ be the optimal value function with
	system state $\mathbf{m}$, the current arrival being customer type
	$i \in \mathcal{I}$ and $n$ periods to go. Then $\forall i \in
	\mathcal{I}$,
	\[
	V^f_n(\mathbf{m} | i) = \max_{\mathbf{d}} \{\sum_{j=1}^J
	p_{ij}(\mathbf{m},\mathbf{d}) [1 + V^f_{n-1}(\mathbf{m} -
	\mathbf{e}_{j})],
	\]
	where $\mathbf{d}$ is the offered set and
	$p_{ij}(\mathbf{m},\mathbf{d})$ is the probability that slot type
	$j$ will be taken if $\mathbf{d}$ is offered and the arrival is type
	$i$ customer.  It is not difficult to see that the optimal offer set
	will be the slot type $j^*(i)$ (which is a function of $i$) such
	that
	\begin{equation}
	\label{eqn:j.star.i} j^*(i) = \argmax_{j \in \{k: \Omega_{ik}=1, k
		\in \mathcal{J}\}} V^f_{n-1}(\mathbf{m} - \mathbf{e}_{j}).
	\end{equation}
	That is, for any arriving customer type, the optimal action is to
	offer the slot type that is acceptable by this customer type and
	that leads to the largest value-to-go. It follows that
	\[
	V^f_n(\mathbf{m}) = \sum_{i \in \mathcal{I} }  \lambda_i
	V^f_n(\mathbf{m} | i) = 1+ \sum_{i \in \mathcal{I} } \lambda_i
	V^f_{n-1}(\mathbf{m} - \mathbf{e}_{j^*(i)}) = 1 + \sum_{i \in
		\mathcal{I}} \lambda_i V^s_{n-1}(\mathbf{m} - \mathbf{e}_{j^*(i)}) =
	V^s_n(\mathbf{m}).
	\]
	To see the last equality, note that the optimal action stipulated by
	Theorem \ref{lemma.one.by.one.seq} ensures that (i) for any arriving
	customer type, it will find an acceptable slot type and (ii) this
	accepted slot type leads to the largest value-to-go among all slot
	types accepted by this customer type. This is exactly enforced by
	\eqref{eqn:j.star.i}. \qed

	
	\subsection{Proof of Theorem~\ref{thm.vn.order}} \proof{Proof.} For
	notational convenience, here we consider the case when
	$\lambda_0=0$. The proof for the case when $\lambda_0 > 0$ follows
	similar steps. In light of Theorem \ref{lemma.one.by.one.seq}, it
	suffices to show that for any $j_1,~j_2$ such that $I(j_1) \subset
	I(j_2)$,
	\begin{equation}
	\label{eqn.vn.order.nested} V_{n-1}(\mathbf{m} - \bfe_{j_1}) \geq
	V_{n-1}(\mathbf{m} - \bfe_{j_{2}}), \forall~n=1,2,\dots, N.
	\end{equation}
	It is easy to see that \eqref{eqn.vn.order.nested} holds for
	$n=1$.
	%
	Assume it holds up to $n>1$ and consider $n+1$. Let us consider a
	few cases below.
	
	Case 1: $\bfm_{j_1} \geq 2, \bfm_{j_2} \geq 1$. Let $g^*$ be an
	optimal action for the state $\bfm-\bfe_{j_2}$. Let
	$B_{ij}^{g}(\bfm)$ denote the probability that a type $i$ customers
	will choose type $j$ slots in state $\bfm$ when action $g$ is taken.
	If $j=0$, then no slots are chosen. Note that $g^*$ is always
	feasible for state $\bfm-\bfe_{j_1}$, and that $B_{ij}^{g^*}(\bfm -
	\bfe_{j_1}) = B_{ij}^{g^*}(\bfm - \bfe_{j_2})$. Thus,
	\begin{align*}
	&V_n(\bfm - \bfe_{j_1}) \geq V_n^{g^*}(\bfm - \bfe_{j_1}) = \sum_{i=1}^I \lambda_i \sum_{j=0}^J B_{ij}^{g^*}(\bfm - \bfe_{j_1}) [\mathbbm{1}_{j>0} + V_{n-1}(\bfm - \bfe_{j_1} - \bfe_{j})]  \\
	\geq &\sum_{i=1}^I \lambda_i \sum_{j=0}^J B_{ij}^{g^*} (\bfm - \bfe_{j_2}) [\mathbbm{1}_{j>0} + V_{n-1}(\bfm -\bfe_{j_2}- \bfe_{j})] = V_n^{g^*}(\bfm - \bfe_{j_2}) = V_n(\bfm - \bfe_{j_2}),
	\end{align*}
	where the second inequality follows from the induction hypothesis.
	
	Case 2: $\bfm_{j_1} =1, \bfm_{j_2} \geq 2$. Again, let $g^*$ be an
	optimal action for the state $\bfm-\bfe_{j_2}$. Following induction
	hypothesis, we choose $g^*$ so that a slot type with a smaller set
	of covered customer types will be offered before any other slot type
	with a larger covered set of customer types. Thus, $g^*$ offers
	$j_1$ before offering $j_2$. Let $\tilde{g}$ be an action that
	follows exactly as $g^*$ except that $\tilde{g}$ does not offer type
	$j_1$ slots. It is clear that $\tilde{g}$ is feasible for state
	$\bfm-\bfe_{j_1}$. There are two subcases.
	
	Case 2a: None of the slot types if any offered between $j_1$ and
	$j_2$ by $g^*$ are acceptable by customer type $i_1$, $\forall i_1
	\in \underline{I}(j_1)$ where $\underline{I}(j_1)$ represent the set
	of customer types who would choose slot type $j_1$ when it is
	offered by $g^*$ at state $\bfm - \bfe_{j_2}$. The following
	inequalities hold.
	\begin{align*}
	B_{ij}^{\tilde{g}}(\bfm - \bfe_{j_1}) =     B_{ij}^{g^*}(\bfm - \bfe_{j_2}), & ~\forall i \notin \underline{I}(j_1), ~\forall j; \\
	B_{i j_1}^{g^*}(\bfm - \bfe_{j_2})  =B_{i j_2}^{\tilde{g}}(\bfm -
	\bfe_{j_1}) =1,& ~\forall i \in \underline{I}(j_1).
	\end{align*}
	It follows that
	\begin{align*}
	& V_n(\bfm - \bfe_{j_1}) \geq V_n^{\tilde{g}}(\bfm - \bfe_{j_1}) = \sum_{i=1}^I \lambda_i \sum_{j=1}^J B_{ij}^{\tilde{g}}(\bfm - \bfe_{j_1}) [1+V_{n-1}(\bfm - \bfe_{j_1} - \bfe_{j})]  \\
	= & \sum_{i \notin \underline{I}(j_1)} \lambda_i \sum_{j=1}^J B_{ij}^{\tilde{g}}(\bfm - \bfe_{j_1}) [1+ V_{n-1}(\bfm - \bfe_{j_1} - \bfe_{j})] + \sum_{i \in \underline{I}(j_1)} \lambda_{i} B_{i j_2}^{\tilde{g}}(\bfm - \bfe_{j_1})  [1+V_{n-1}(\bfm - \bfe_{j_1} - \bfe_{j_2})]  \\
	\geq & \sum_{i \notin \underline{I}(j_1)} \lambda_i \sum_{j=1}^J B_{ij}^{g^*}(\bfm - \bfe_{j_2}) [1+ V_{n-1}(\bfm - \bfe_{j_2} - \bfe_{j})] + \sum_{i \in \underline{I}(j_1)} \lambda_{i}  B_{i j_1}^{g^*}(\bfm - \bfe_{j_2}) [1+V_{n-1}(\bfm - \bfe_{j_2} - \bfe_{j_1})] \\
	= & V_{n-1}(\bfm - \bfe_{j_2}),
	\end{align*}
	where the last inequality follows from the induction hypothesis.
	
	Case 2b: Let $\underline{I}(j_k)$ be the set of the customer types
	who will actually choose slot type $j_k$ under $g^*$ when it is
	offered at state $\bfm - \bfe_{j_2}$. One slot type, say $j_3$,
	offered between $j_1$ and $j_2$ by $g^*$ is acceptable by some
	customer type $i_1 \in \underline{I}(j_1)$. Consider $j_3$ to be the
	only one of such slots (extension to multiple of such slots uses a
	similar but more tedious proof). Because slot types $j_1$, $j_2$ and
	$j_3$ can cover some same customer types, we know that either
	$I(j_3) \subset I(j_k)$ or $I(j_3) \supset I(j_k)$, $k=1,2$, by the
	preassumption of the theorem. Because we choose $g^*$ based on the
	size of covered set of customer types, we have that $I(j_1) \subset
	I(j_2) \subset I(j_3)$. Also note that $\cap_{k=1}^3
	\underline{I}(j_k) = \emptyset$ because any customer type can choose
	at most one slot type. Let $j(i)$ be the slot type chosen by
	customer type $i$, $i \notin \underline{I}(j_1) \cup
	\underline{I}(j_2) \cup \underline{I}(j_3)$. Note that $j(i)$ is the
	same under $g^*$ or $\tilde{g}$, $\forall i \notin
	\underline{I}(j_1) \cup \underline{I}(j_2) \cup \underline{I}(j_3)$.
	It follows that
	\begin{align*}
	& V_n(\bfm - \bfe_{j_1}) - V_n(\bfm - \bfe_{j_2}) \geq V_n^{\tilde{g}}(\bfm - \bfe_{j_1}) - V^{g^*}_{n-1}(\bfm - \bfe_{j_2}) \\
	= & \sum_{i \notin \underline{I}(j_1) \cup \underline{I}(j_2) \cup  \underline{I}(j_3)} \lambda_i [V_{n-1}(\bfm - \bfe_{j_1} - \bfe_{j(i)}) - V_{n-1}(\bfm - \bfe_{j_2} - \bfe_{j(i)}) ] \\
	& + \sum_{i \in \underline{I}(j_1) \cup \underline{I}(j_3)} \lambda_i V_{n-1}(\bfm - \bfe_{j_1} - \bfe_{j_3}) +  \sum_{i \in \underline{I}(j_2)} \lambda_i V_{n-1}(\bfm - \bfe_{j_1} - \bfe_{j_2}) \\
	& - \{ \sum_{i \in \underline{I}(j_1)} \lambda_i V_{n-1}(\bfm - \bfe_{j_2} - \bfe_{j_1}) + \sum_{i \in \underline{I}(j_3)} \lambda_i V_{n-1}(\bfm - \bfe_{j_2} - \bfe_{j_3}) + \sum_{i \in \underline{I}(j_2)} \lambda_i V_{n-1}(\bfm - \bfe_{j_2} - \bfe_{j_2}) \} \geq
	0,
	\end{align*}
	where the last inequality follows from the induction hypothesis.

	Case 3: $\bfm_{j_1} =1, \bfm_{j_2} = 1$. Let $g^*$ be an optimal
	action for the state $\bfm-\bfe_{j_2}$ ($g^*$ does not offer slot
	type $j_2$ because none is available). Let $\tilde{g}$ be an action
	that follows exactly as $g^*$ except that $\tilde{g}$ does not offer
	type $j_1$ slots but offers type $j_2$ at the end. It is clear that
	$\tilde{g}$ is feasible for state $\bfm-\bfe_{j_1}$. Let
	$\underline{I}(j_2)$ be the customer types who choose $j_2$ under
	${\tilde g}$; these customers do not book any appointments under
	$g^*$. Let $j(i)$ be the slot type actually chosen by customer type
	$i$, $i \notin \underline{I}(j_2)$. Note that $j(i)$ is the same
	under $g^*$ or $\tilde{g}$.  It follows that
	\begin{align*}
	& V_n(\bfm - \bfe_{j_1}) - V_{n}(\bfm - \bfe_{j_2}) \geq V_n^{\tilde{g}}(\bfm - \bfe_{j_1}) - V^{g^*}_{n}(\bfm - \bfe_{j_2}) \\
	= & \sum_{i \notin \underline{I}(j_2) } \lambda_i [V_{n-1}(\bfm - \bfe_{j_1} - \bfe_{j(i)}) + V_{n-1}(\bfm - \bfe_{j_2} - \bfe_{j(i)} ) ] + \sum_{i \in \underline{I}(j_2) } \lambda_i [1 + V_{n-1}(\bfm - \bfe_{j_1} - \bfe_{j_2}) - V_{n-1}(\bfm - \bfe_{j_2}) ] \\
	& \geq 0,
	\end{align*}
	where the last inequality follows from the induction hypothesis and
	Lemma \ref{lemma.Vn.monotone.seq}. This completes the whole proof.
	\qed

	\subsection{Proof of Proposition~\ref{theorem.W.seq.opt.structure}}
	
	\proof{Proof.}
	It
	suffices to show the following monotonic results for the ``W'' model
	with sequential offers: $V_n(\bfm - \bfe_1) - V_{n}(\bfm -\bfe_2)$
	increases as $m_1$ increases and that $ V_n(\bfm - \bfe_2) -
	V_{n}(\bfm  -\bfe_1)$ increases as $m_2$ increases, $\forall n \geq
	1, ~ \forall \bfm \geq (1,1)$.
	
	That is,
	\begin{equation}
	\label{eq:monotone1.W.seq} V_n(\bfm) - V_{n}(\bfm + \bfe_1 -\bfe_2)
	\geq  V_n(\bfm - \bfe_1) - V_{n}(\bfm  -\bfe_2),  ~\forall n \geq 1,
	~\forall \bfm \geq (1,1),
	\end{equation}
	and
	\begin{equation}
	\label{eq:monotone2.W.seq} V_n(\bfm) - V_{n}(\bfm + \bfe_2 -\bfe_1)
	\geq  V_n(\bfm - \bfe_2) - V_{n}(\bfm  -\bfe_1), ~\forall n \geq 1,
	~\forall \bfm \geq (1,1).
	\end{equation}
	
	To facilitate the proof of \eqref{eq:monotone1.W.seq} and
	\eqref{eq:monotone2.W.seq}, we introduce a few notations. Let
	$\Delta^A_n(\bfm) = V_n(\bfm - \bfe_1) - V_{n}(\bfm -\bfe_2)$ and
	$\Delta^B_n(\bfm) =V_n(\bfm - \bfe_2) - V_{n}(\bfm -\bfe_1)$. Note
	that \eqref{eq:monotone1.W.seq} and \eqref{eq:monotone2.W.seq} are
	symmetric, and thus we limit ourselves to just prove
	\eqref{eq:monotone1.W.seq}.
	
	Consider the case for $n=1$. At $\bfm=(1,1)$, $\Delta^A_1(\bfm) =
	V_1(0,1) - V_1(1,0) = (\lambda_2+\lambda_3) - (\lambda_1 +
	\lambda_2) = \lambda_3 - \lambda_1$. For $\bfm =(m_1, 1)$ and $m_1
	\geq 2$, we have $\Delta^A_1(\bfm) = V_1(m_1-1,1) - V_1(m_1,0) =
	(1-\lambda_0) - (\lambda_1+\lambda_2) = \lambda_3$. Thus,
	\eqref{eq:monotone1.W.seq} holds for $n=1$, $\bfm = (m_1, 1)$ and
	$m_1 \geq 1$. Now, for $n=1$ and $\bfm = (1, m_2)$ and $m_2 \geq 2$,
	we have $\Delta^A_1(\bfm) = V_1(0,m_2) - V_1(1,m_2-1) =
	(\lambda_2+\lambda_3) - (1-\lambda_0) = -\lambda_1$. Consider $n=1$,
	$\bfm = (m_1, m_2)$, and $m_1, m_2 \geq 2$. In this case, we have
	that $\Delta^A_1(\bfm) = V_1(m_1-1,m_2) - V_1(m_1,m_2-1) =
	(1-\lambda_0) - (1-\lambda_0) = 0$. Thus, \eqref{eq:monotone1.W.seq}
	holds for $n=1$, $\bfm = (m_1, m_2)$ and $m_1 \geq 1, m_2 \geq 2$.
	This completes the proof of \eqref{eq:monotone1.W.seq} for $n=1$.
	
	Assume that \eqref{eq:monotone1.W.seq} holds up to $n=k$ for $\bfm
	\geq (1,1)$. We will use induction below to show that this is also
	true for $n=k+1$. We start by witting the Bellman's equation below.
	\begin{align}
	\label{eqn.W.bellman.seq}
	& V_{k+1}(\bfm) = \nonumber \\
	& \max \left\{
	\begin{array}{lll}
	1-\lambda_0 + (\lambda_1 + \frac{1}{2} \lambda_2) V_k(\bfm - \bfe_1) + (\frac{1}{2 }\lambda_2 + \lambda_3) V_k(\bfm - \bfe_2) + \lambda_0 V_k(\bfm) , \\ 
	1-\lambda_0 + (\lambda_1 + \lambda_2) V_k(\bfm - \bfe_1) + \lambda_3 V_k(\bfm - \bfe_2) + \lambda_0 V_k(\bfm), \\ 
	1-\lambda_0 + \lambda_1  V_k(\bfm - \bfe_1) + (\lambda_2 + \lambda_3) V_k(\bfm - \bfe_2) + \lambda_0 V_k(\bfm). 
	\end{array}
	\right\},
	\end{align}
	where the three terms in the max operator correspond to actions
	$\{1,2\}$, $\{1\}$-$\{2\}$ and $\{2\}$-$\{1\}$, respectively. Action
	$\{S_1\}$-$\{S_2\}$ offers subset $S_1$ followed by subset $S_2$.
	For ease of notation, we define $\Delta_{k+1}^{ij}(\bfm)$ to be the
	difference of the $i$th and $j$th terms in the max operator
	\eqref{eqn.W.bellman.seq} above, $i,j \in \{1,2,3\}$. It follows
	that
	\[
	\Delta_{k+1}^{21}(\bfm) = \frac{1}{2} \lambda_2 [V_k(\bfm - \bfe_1)
	- V_k(\bfm - \bfe_2)] = \frac{1}{2} \lambda_2 \Delta^A_k(\bfm),
	\]
	and
	\[
	\Delta_{k+1}^{31}(\bfm) = \frac{1}{2} \lambda_2 [V_k(\bfm - \bfe_2)
	- V_k(\bfm - \bfe_1)] = \frac{1}{2} \lambda_2 \Delta^B_k(\bfm).
	\]
	Because $\Delta_{k+1}^{21}(\bfm) + \Delta_{k+1}^{31}(\bfm)=0$, one
	of these two terms must be non-negative suggesting one of the
	corresponding actions is optimal. In particular, if
	$\Delta_{k+1}^{21}(\bfm) \geq 0$, or equivalently, $\Delta^A_k(\bfm)
	\geq 0$, the the optimal action is  $\{1\}$-$\{2\}$; otherwise, it
	would be $\{2\}$-$\{1\}$.
	
	To prove the desired result, we need to consider the following
	cases. Case (1): $\bfm=(1,1)$;  case (2): $\bfm=(m_1,1), ~m_1 \geq
	2$; case (3): $\bfm=(1,m_2), ~ m_2 \geq 2$; and case (4), $\bfm \geq
	(2,2)$.

	For Case (1) with $\bfm = (1,1)$, we have
	\[
	\Delta^A_{k+1}(1,1) =  V_{k+1}(0,1) - V_{k+1}(1,0) =
	[(1-\lambda_1-\lambda_0) + (\lambda_1+\lambda_0) V_k(0,1)]
	-[(1-\lambda_3-\lambda_0) + (\lambda_3+\lambda_0) V_k(1,0) ].
	\]
	We consider two subcases. Case (1a): if at state $(1,1)$ the optimal
	action is $\{1\}$-$\{2\}$, then
	\begin{align*}
	& \Delta^A_{k+1}(2,1) = V_{k+1}(1,1) - V_{k+1}(2,0)\\
	= & [(1-\lambda_0) + (\lambda_1+\lambda_2) V_k(0,1) + \lambda_3 V_k(1,0) + \lambda_0 V_k(1,1)] \\
	& -[(1-\lambda_3-\lambda_0) + (\lambda_1+\lambda_2) V_k(1,0) + (\lambda_3+\lambda_0) V_k(2,0) ].
	\end{align*}
	It follows that
	\begin{align*}
	& \Delta^A_{k+1}(2,1) -  \Delta^A_{k+1}(1,1) \\
	= & \lambda_1 [1-V_k(1,0)] + \lambda_2 [V_k(0,1) - V_k(1,0)] + \lambda_3 [2 V_k(1,0) - V_k(2,0)] + \lambda_0 [\Delta^A_k(2,1) - \Delta^A_k(1,1)] \\
	= & \lambda_1 [1-V_k(1,0)] + \lambda_2 \Delta^A_k(1,1) + \lambda_3
	[2 V_k(1,0) - V_k(2,0)] + \lambda_0 [\Delta^A_k(2,1) -
	\Delta^A_k(1,1)].
	\end{align*}
	It is trivial that $1-V_k(1,0) \geq 0$. We also know that
	$\Delta^A_k(1,1) \geq 0$ in this case because the optimal action is
	$\{1\}$-$\{2\}$; and that $\Delta^A_k(2,1) - \Delta^A_k(1,1) \geq 0$
	by the induction hypothesis. Finally, we claim that
	\begin{equation}
	\label{eq:claim} 2 V_k(1,0) - V_k(2,0) \geq 0, ~\forall k \geq 1,
	\end{equation}
	which will be shown at the end of this proof. Thus,
	$\Delta^A_{k+1}(2,1) -  \Delta^A_{k+1}(1,1) \geq 0$ if the optimal
	action is $\{1\}$-$\{2\}$ at state $(1,1)$.

	Case (1b): if the optimal action at state $(m_1,1)$ is
	$\{2\}$-$\{1\}$, then
	\begin{align*}
	& \Delta^A_{k+1}(2,1) = V_{k+1}(1,1) - V_{k+1}(2,0)\\
	= & [(1-\lambda_0) + \lambda_1 V_k(0,1) + (\lambda_2 + \lambda_3) V_k(1,0) + \lambda_0 V_k(1,1)] \\
	& -[(1-\lambda_3-\lambda_0) + (\lambda_1+\lambda_2) V_k(1,0) + (\lambda_3+\lambda_0) V_k(2,0) ].
	\end{align*}
	It follows that
	\begin{align*}
	& \Delta^A_{k+1}(2,1) -  \Delta^A_{k+1}(1,1) \\
	= & \lambda_1 [1-V_k(1,0)]  + \lambda_3 [2 V_k(1,0) - V_k(2,0)] +
	\lambda_0 [\Delta^A_k(2,1) - \Delta^A_k(1,1)] \geq 0.
	\end{align*}
	In summary, cases (1a) and (1b) collectively show that
	$\Delta^A_{k+1}(2,1) -  \Delta^A_{k+1}(1,1) \geq 0$.
	
	For Case (2) with $\bfm =(m_1, 1), ~m_1 \geq 2$, we evaluate
	$V_{k+1}(m_1-1,1) - V_{k+1}(m_1,0)$ in the following two subcases.
	Case(2a): if the optimal action at state $(m_1,1)$ is
	$\{1\}$-$\{2\}$, then
	\begin{align*}
	& V_{k+1}(m_1-1,1) - V_{k+1}(m_1,0) \\
	= & [(1-\lambda_0) + (\lambda_1 + \lambda_2) V_k(m_1-2, 1) + \lambda_3 V_k(m_1-1, 0) + \lambda_0 V_k(m_1-1,1)] \\
	& -[(1-\lambda_3-\lambda_0) + (\lambda_1+\lambda_2) V_k(m_1-1,0) + (\lambda_3+\lambda_0) V_k(m_1,0) ] \\
	= & \lambda_3 + (\lambda_1 + \lambda_2) \Delta^A_k(m_1-1,1) + \lambda_0 \Delta^A_k(m_1,1) + \lambda_3 [(1-(\lambda_1+\lambda_2)^{m_1-1}) - (1-(\lambda_1+\lambda_2)^{m_1})] \\
	= & \lambda_3 + (\lambda_1 + \lambda_2) \Delta^A_k(m_1-1,1) +
	\lambda_0 \Delta^A_k(m_1,1) - \lambda_3 (\lambda_3+\lambda_0)
	(\lambda_1+\lambda_2)^{m_1-1},
	\end{align*}
	which increases as $m_1$ increases by the induction hypothesis.
	
	Case (2b): if the optimal action at state $(m_1,1)$ is
	$\{2\}$-$\{1\}$, then for $\bfm =(m_1, 1), ~m_1 \geq 2$,
	\begin{align*}
	& V_{k+1}(m_1-1,1) - V_{k+1}(m_1,0) \\
	= & [(1-\lambda_0) + \lambda_1 V_k(m_1-2, 1) + (\lambda_2+\lambda_3) V_k(m_1-1, 0) + \lambda_0 V_k(m_1-1,1)] \\
	& -[(1-\lambda_3-\lambda_0) + (\lambda_1+\lambda_2) V_k(m_1-1,0) + (\lambda_3+\lambda_0) V_k(m_1,0) ] \\
	= & \lambda_3 + \lambda_1 \Delta^A_k(m_1-1,1) + \lambda_0 \Delta^A_k(m_1,1) + \lambda_3 [(1-(\lambda_1+\lambda_2)^{m_1-1}) - (1-(\lambda_1+\lambda_2)^{m_1})] \\
	= & \lambda_3 + \lambda_1 \Delta^A_k(m_1-1,1) + \lambda_0
	\Delta^A_k(m_1,1) - \lambda_3 (\lambda_3+\lambda_0)
	(\lambda_1+\lambda_2)^{m_1-1},
	\end{align*}
	which also increases as $m_1$ increases by the induction hypothesis.
	Thus, cases (1a) though (1d) shows that $\Delta^A_n(\bfm)$ increases
	in $m_1$ for $n \geq 1$ and $\bfm=(m_1,1),m_1 \geq 1$.
	
	Case (3): $\bfm=(1,m_2), ~ m_2 \geq 2$. We want to show that
	\[
	\Delta^A_{k+1}(2,m_2) - \Delta^A_{k+1}(1,m_2) =  [V_{k+1}(1,m_2) -
	V_{k+1}(2,m_2-1)] - [V_{k+1}(0,m_2) - V_{k+1}(1,m_2-1)] \geq 0.
	\]
	Again, we separate into a few subcases. If the optimal action at
	state $(1,m_2-1)$ is $\{1\}$-$\{2\}$, then $\Delta_k^A(1,m_2-1) \geq
	0$. It follows that $\Delta_k^A (2, m_2-1) \geq 0$ by the induction
	hypothesis, and the optimal action at state $(2,m_2-1)$ is also
	$\{1\}$-$\{2\}$. But the optimal actions at state $(1,m_2)$ can
	still be either $\{1\}$-$\{2\}$ or $\{2\}$-$\{1\}$. Following this
	logic, we need to consider four subcases. Case (3a): the optimal
	actions at state $(1,m_2-1)$, $(2,m_2-1)$ and $(1,m_2)$ are all
	$\{1\}$-$\{2\}$. Case (3b): the optimal actions at state
	$(1,m_2-1)$, $(2,m_2-1)$ and $(1,m_2)$ are $\{1\}$-$\{2\}$,
	$\{1\}$-$\{2\}$ and $\{2\}$-$\{1\}$, respectively. Case (3c): the
	optimal actions at state $(1,m_2-1)$, $(2,m_2-1)$ and $(1,m_2)$ are
	all $\{2\}$-$\{1\}$. Case (3d): the optimal actions at state
	$(1,m_2-1)$, $(2,m_2-1)$ and $(1,m_2)$ are $\{2\}$-$\{1\}$,
	$\{1\}$-$\{2\}$ and $\{2\}$-$\{1\}$, respectively.
	
	For case (3a), we have
	\begin{align*}
	& \Delta^A_{k+1}(2,m_2) - \Delta^A_{k+1}(1,m_2) \\
	=& \lambda_1 [1 - V_k(0,m_2-1)] + (\lambda_1+\lambda_2) \Delta_k^A(1,m_2)\\
	& + \lambda_3 [\Delta_k^A(2,m_2-1) - \Delta_k^A(1,m_2-1)] +
	\lambda_0 [\Delta_k^A(2,m_2) - \Delta_k^A(1,m_2)] \geq 0,
	\end{align*}
	where the inequality follows from the fact that the first term is
	trivially nonnegative, the second term is positive as the optimal
	actions at state $(1,m_2)$ is $\{1\}$-$\{2\}$ and the last two terms
	are nonnegative by the induction hypothesis.
	
	For case (3b), we have
	\begin{align*}
	& \Delta^A_{k+1}(2,m_2) - \Delta^A_{k+1}(1,m_2) \\
	=& \lambda_1 [ - V_k(1,m_2-1) + 1 + V_k(0,m_2-1) ]\\
	& + \lambda_3 [\Delta_k^A(2,m_2-1) - \Delta_k^A(1,m_2-1)] +\lambda_0
	[\Delta_k^A(2,m_2) - \Delta_k^A(1,m_2)] \geq 0,
	\end{align*}
	where the first term is nonnegative following Lemma
	\ref{lemma.Vn.monotone.seq} and the other two terms are nonnegative
	following the induction hypothesis.
	
	For case (3c), we have
	\begin{align*}
	& \Delta^A_{k+1}(2,m_2) - \Delta^A_{k+1}(1,m_2) \\
	=& \lambda_1 [- V_k(1,m_2-1) + 1 + V_k(0,m_2-1) ]\\
	& + (\lambda_2 + \lambda_3) [\Delta_k^A(2,m_2-1) -
	\Delta_k^A(1,m_2-1)] +\lambda_0 [\Delta_k^A(2,m_2) -
	\Delta_k^A(1,m_2)] \geq 0,
	\end{align*}
	following a similar argument of case (3b).
	
	For case (3d), we have
	\begin{align*}
	& \Delta^A_{k+1}(2,m_2) - \Delta^A_{k+1}(1,m_2) \\
	=& \lambda_1 [ - V_k(1,m_2-1) + 1 + V_k(0,m_2-1) ]\\
	& - \lambda_2 \Delta^A_k(1,m_2-1) +  \lambda_3 [\Delta_k^A(2,m_2-1)
	- \Delta_k^A(1,m_2-1)] +\lambda_0 [\Delta_k^A(2,m_2) -
	\Delta_k^A(1,m_2)] \geq 0,
	\end{align*}
	following a similar logic of case (3b) and the fact that
	$\Delta^A_k(1,m_2-1) \leq 0$ (because the optimal action at state
	$(1,m_2-1)$ is $\{2\}$-$\{1\}$). This completes the proof of case
	(3).

	For case (4) $\bfm \geq (2,2)$, we evaluate $V_n(\bfm - \bfe_1) -
	V_{n}(\bfm  -\bfe_2)$ and need to consider four subcases.  Case
	(4a): if the optimal actions at states $(\bfm - \bfe_1)$ and $(\bfm
	- \bfe_2)$ are both $\{1\}$-$\{2\}$. Then
	\begin{align*}
	& \Delta_{k+1}^A(\bfm) = V_{k+1}(\bfm - \bfe_1) - V_{k+1}(\bfm  -\bfe_2) \\
	=& (\lambda_1+\lambda_2) \Delta_k^A(\bfm - \bfe_1) + \lambda_3
	\Delta_k^A(\bfm - \bfe_2) + \lambda_0 \Delta_k^A(\bfm),
	\end{align*}
	which increases in $m_1$ by the induction hypothesis. Case (4b): if
	the optimal actions at states $(\bfm - \bfe_1)$ and $(\bfm -
	\bfe_2)$ are both $\{2\}$-$\{1\}$. Then,
	\begin{align*}
	& \Delta_{k+1}^A(\bfm) = V_{k+1}(\bfm - \bfe_1) - V_{k+1}(\bfm  -\bfe_2) \\
	=& \lambda_1 \Delta_k^A(\bfm - \bfe_1) + (\lambda_2+\lambda_3)
	\Delta_k^A(\bfm - \bfe_2) + \lambda_0 \Delta_k^A(\bfm),
	\end{align*}
	which increases in $m_1$ by the induction hypothesis. Case (4c): if
	the optimal actions at states $(\bfm - \bfe_1)$ and $(\bfm -
	\bfe_2)$ are $\{1\}$-$\{2\}$ and $\{2\}$-$\{1\}$, respectively.
	Then,
	\begin{align*}
	& \Delta_{k+1}^A(\bfm) = V_{k+1}(\bfm - \bfe_1) - V_{k+1}(\bfm  -\bfe_2) \\
	=& \lambda_1 \Delta_k^A(\bfm - \bfe_1) + \lambda_2 [V_k(\bfm - 2\bfe_1) - V_k(\bfm - 2\bfe_2) ] + \lambda_3 \Delta_k^A(\bfm - \bfe_2) + \lambda_0 \Delta_k^A(\bfm) \\
	=& \lambda_1 \Delta_k^A(\bfm - \bfe_1) + \lambda_2 [\Delta_k^A(\bfm
	- \bfe_1) + \Delta_k^A(\bfm - \bfe_2) ] + \lambda_3 \Delta_k^A(\bfm
	- \bfe_2) + \lambda_0 \Delta_k^A(\bfm),
	\end{align*}
	which increases in $m_1$ by the induction hypothesis. Case (4d): if
	the optimal actions at states $(\bfm - \bfe_1)$ and $(\bfm -
	\bfe_2)$ are $\{2\}$-$\{1\}$ and $\{1\}$-$\{2\}$, respectively.
	Then,
	\begin{align*}
	\Delta_{k+1}^A(\bfm) = V_{k+1}(\bfm - \bfe_1) - V_{k+1}(\bfm  -\bfe_2)
	= \lambda_1 \Delta_k^A(\bfm - \bfe_1)  + \lambda_3 \Delta_k^A(\bfm -
	\bfe_2) + \lambda_0 \Delta_k^A(\bfm),
	\end{align*}
	which increases in $m_1$ by the induction hypothesis.
	
	Finally, we show our claim \eqref{eq:claim}, which can be easily
	done by induction. When $k=1$, $2 V_k(1,0) - V_k(2,0) =
	2(\lambda_1+\lambda_2) - [1-(1-\lambda_1-\lambda_2)^2] =
	(\lambda_1+\lambda_2)^2 \geq 0.$ Assume this holds up to $k=u$.
	Consider $k=u+1$. We have that
	\begin{align*}
	& 2 V_{u+1}(1,0) - V_{u+1}(2,0) \\
	= & 2 [(\lambda_1+\lambda_2) + (\lambda_3+\lambda_0) V_{u}(1,0) ] - [(\lambda_1+\lambda_2) + (\lambda_1+\lambda_2)V_u(1,0) +(\lambda_3+\lambda_0) V_{u}(2,0)] \\
	= & (\lambda_1+\lambda_2)[1 - V_u(1,0)] + (\lambda_3+\lambda_0) [2
	V_{u}(1,0) - V_{u}(2,0)] \geq 0,
	\end{align*}
	proving our claim \eqref{eq:claim} and completing the whole proof.
	\endproof

	\newpage
	
	\section{Additional Numerical Results}
	
	\begin{table}[!ht]
		\caption{\small 
			Optimality gap of the static randomized policy $\pi^{p^*}$ in the M
			model instance.} \label{tab:M.random}
		\begin{center}
			{\small
				\begin{tabular}{c | c | ccc | ccc | ccc}
					\toprule
					\multirow{2}{*}{$N$}     &                   {\footnotesize $\#$ of}   & \multicolumn{3}{c|}{\footnotesize $(\lambda_1,\lambda_2) = (1/2, 1/2)$} & \multicolumn{3}{c|}{\footnotesize $(\lambda_1,\lambda_2) = (1/3, 2/3)$} & \multicolumn{3}{c}{\footnotesize $(\lambda_1,\lambda_2) = (1/4, 3/4)$} \\
					&           {\footnotesize Scenarios}            & \footnotesize Max & \footnotesize Average  & \footnotesize Median   & \footnotesize Max   &  \footnotesize Average  & \footnotesize Median   & \footnotesize Max &  \footnotesize Average  & \footnotesize Median  \\
					\otoprule
					$20$                                          & 45                     &  -10.7\% &  -7.7\%  &  -7.2\% &   -9.8\%  &  -7.9\% &  -8.0\%     & -10.8\%  &  -8.7\%  & -8.7\%                     \\
					$30$                                          & 91                     &  -9.1\% &  -6.4\%  &  -5.9\% &  -8.7\%   &  -6.7\%  & -6.7\%     & -8.8\%  &  -7.2\%  & -7.1\%             \\
					$40$                                          & 153                    &  -8.1\% &  -5.7\%  &  -5.2\% &  -7.9\%   &  -5.8\%   & -5.8\%     & -7.8\%  & -6.2\%  & -6.2\%         \\
					$50$                                           & 231                   &  -7.5\% &  -5.2\%  & -4.6\%  &  -7.1\%   &  -5.3\%   & -5.2\%     & -7.0\%  & -5.6\%  & -5.5\%       \\
					\bottomrule
				\end{tabular}
			}
		\end{center}
	\end{table}

	\begin{table}[h!]
		\caption{Policy Comparison in a Multi-day Scheduling Setting (with
			Poisson Arrivals). \label{tab:multiday.poi}} {
			\begin{tabular}{c | c | c | ccc | ccc | ccc}
				\toprule
				& \multirow{2}{*}{$D$}     &                    {\footnotesize $\#$ of}   & \multicolumn{3}{c|}{\footnotesize $(\lambda_1,\lambda_2) = (1/2, 1/2)$} & \multicolumn{3}{c|}{\footnotesize $(\lambda_1,\lambda_2) = (1/3, 2/3)$} & \multicolumn{3}{c}{\footnotesize $(\lambda_1,\lambda_2) = (1/4, 3/4)$} \\
				&                         &           {\footnotesize Scenarios}            & \footnotesize Max & \footnotesize Average  & \footnotesize Median   & \footnotesize Max   &  \footnotesize Average  & \footnotesize Median   & \footnotesize Max &  \footnotesize Average  & \footnotesize Median  \\
				\otoprule
				&  $1$                                          & 45                     &  4.1\%  &  2.1\%  & 2.0\%  & 4.2\%  &  2.4\%  & 2.1\%     & 3.7\%  &  2.4\%  & 2.5\%                      \\
				{\footnotesize Non-sequential Optimal} &  $2$     & 91                    &  3.7\%  &  1.9\%  & 1.8\%  &  3.4\%  &  2.0\%  & 2.2\%     & 2.9\%  &  2.0\%  & 2.1\%            \\
				{\footnotesize vs. Offering-all}  & $3$              & 153                &  3.0\%  &  1.5\%  & 1.5\%  &  1.6\%  &  1.7\%  & 2.3\%     & 2.3\%  &  1.6\%  & 1.6\%          \\
				&  $4$                                          & 231                    &  2.7\%  & 1.2\%  & 1.2\%   &  2.3\%  & 1.3\%  & 1.6\%      & 2.1\%  & 1.4\%  & 1.4\%      \\
				\otoprule
				&  $1$                                          & 45                     &  9.6\%  &  4.0\% &  3.4\%  &  11.2\%  &  5.1\% &  4.2\%     & 10.0\%  &  5.6\% &  5.8\%                      \\
				{\footnotesize Sequential Optimal} &  $2$         & 91                    &  9.8\%   &  3.7\%  & 2.8\% &  9.0\%   &  4.5\%  & 4.4\%     & 7.6\%   &  4.8\%  & 5.1\%             \\
				{\footnotesize vs. Offering-all} &  $3$         & 153                     &  8.6\%   &  3.2\%  & 2.2\% &  7.5\%   &  3.7\%  & 4.2\%     & 6.1\%   &  4.0\%  & 4.1\%         \\
				&  $4$                                          & 231                    &  7.5\%   &  2.6\%   & 1.7\%  &  6.4\%   &  3.1\%   & 3.5\%     & 5.2\%   &  3.3\%   & 3.3\%       \\
				\bottomrule
			\end{tabular}
		}{} 
	\end{table}






\end{document}